\input ./style/arxiv-general.cfg
\documentclass[aap,MSNbibl,seceqn,dvips]{arximspdf}
\makeatletter
   \@ifpackageloaded{graphicx}{}{\usepackage{graphicx}}
\makeatother
\usepackage{mathbh}
\usepackage{url,breakurl}

\doi{10.1214/14-AAP1035} 
\volume{25}
\issue{4}
\pubyear{2015}
\firstpage{1729}
\lastpage{1779}

\makeatletter
\newcommand{\TV}{\mathrm{TV}}
\renewcommand{\mid}{|}
\newcommand{\rrvert}{\vert}
\newcommand{\llvert}{\vert}
\newtheorem{theorem}{Theorem}[section]
\newtheorem{corollary}[theorem]{Corollary}
\newtheorem{proposition}[theorem]{Proposition}
\newproclaim{definition}[theorem]{Definition}
\newproclaim{remark}[theorem]{Remark}

\newcommand{\calA}{\mathcal{A}}
\newcommand{\calC}{C}
\newcommand{\calB}{\mathcal{B}}
\newcommand{\bM}{\mathbf{M}}
\newcommand{\bzero}{\mathbf{0}}
\newcommand{\bP}{\mathbf{P}}
\def\bone{\mathbf{1}}
\def\ba{\mathbf{a}}
\newcommand{\bfa}{\mathbf{a}}
\newcommand{\bfb}{\mathbf{b}}
\newcommand{\bfh}{\mathbf{h}}
\newcommand{\bfg}{\mathbf{g}}
\makeatother

\begin{document}
\begin{frontmatter}

\title{On meteors, earthworms and WIMPs}
\runtitle{On meteors, earthworms and WIMPs}

\begin{aug}
\author[A]{\fnms{Sara}~\snm{Billey}\thanksref{T1}\ead[label=e1]{billey@math.washington.edu}},
\author[A]{\fnms{Krzysztof}~\snm{Burdzy}\corref{}\thanksref{T1}\ead[label=e2]{burdzy@math.washington.edu}},
\author[A]{\fnms{Soumik}~\snm{Pal}\thanksref{T1}\ead[label=e3]{soumik@math.washington.edu}}\\
\and
\author[B]{\fnms{Bruce E.}~\snm{Sagan}\ead[label=e4]{sagan@math.msu.edu}}
\runauthor{Billey, Burdzy, Pal and Sagan}
\affiliation{University of Washington, University of Washington,
University of Washington and Michigan State University}
\address[A]{S. Billey\\
K. Burdzy\\
S. Pal\\
Department of Mathematics\\
University of Washington\\
Box 354350\\
Seattle, Washington 98195\\
USA\\
\printead{e1}\\
\phantom{E-mail: }\printead*{e2}\\
\phantom{E-mail: }\printead*{e3}} 
\address[B]{B. E. Sagan\\
Department of Mathematics\\
Michigan State University\\
East Lansing, Michigan 48824-1027\\
USA\\
\printead{e4}}
\end{aug}
\thankstext{T1}{Supported in part by NSF Grants DMS-11-01017, DMS-12-06276
and DMS-13-08340.}

\received{\smonth{8} \syear{2013}}
\revised{\smonth{3} \syear{2014}}

%
\begin{abstract}
We study a model of mass redistribution on a finite graph.
We address the questions of convergence to equilibrium and the rate of
convergence. We present theorems on the distribution of empty sites
and the distribution of mass at a fixed vertex.
These distributions are related to random permutations with certain
peak~sets.
\end{abstract}

%
\end{frontmatter}

\section{Introduction}\label{intro}

We study a model of mass redistribution on a finite graph.
A~vertex $x$ of the graph holds mass $M^x_t\geq0$ at time $t$. When a
``meteor hits''
$x$ at time $t$, the mass $M^x_t$ of the soil present at $x$ is
distributed equally among all neighbors of $x$ (added to their masses).
There is no soil (mass) left at $x$ just after a meteor hit. Meteor
hits are modeled as independent Poisson processes, one for each vertex
of the graph.

We will address the following questions about the meteor model.
Does the process converge to equilibrium? If so, at what rate?
Assuming that the mass distribution process is in equilibrium, what is
the distribution of ``meteor craters'' (sites with zero mass) at a
fixed time?
In equilibrium, at a fixed time and vertex, what is the
distribution of soil mass? We will answer some of these questions in
the asymptotic sense, for some families of growing graphs.

We will also study an ``earthworm model'' in which the soil
redistribution events do not occur according to the Poisson arrival
process but along the trajectory of a symmetric random walk on the
graph. See Section~\ref{earth} for the motivation of the ``earthworm'' model.

We will now present sources of inspiration, motivation and possible
applications for our main model.
\begin{longlist}[(iii)]
\item[(i)]
Similar models of mass redistribution appeared in \cite{HW}, but that
paper went in a completely different direction. It was mostly focused
on the limit model when the graph approximates the real line.
Continuous mass redistribution also appeared in a version of the
chip-firing model in \cite{CP}, but the updating mechanism in that
paper is different from ours. Mass redistribution is a part of every
sandpile model, including a ``continuous'' version studied in \cite
{FMQR}. Sandpile models have considerably different structures and
associated questions from ours.
Our model is one of the simplest models for mass redistribution.
Therefore, its analysis is likely to be most complete on the
mathematical side---a program that we only start in this paper. The
elementary character of our model makes it amenable to a variety of
mathematical techniques---something that we demonstrate in this
article. Our model can be easily modified and generalized to
accommodate the needs of applied science.

\item[(ii)] More on the theoretical side, our model is related to products of
random matrices.
Let $G$ be a finite graph with a set of vertices $V=\{ 1,2, \ldots, n\}
$. Let $A$ denote the transpose of the transition probability matrix of
the nearest neighbor simple random walk on this graph, and let $I$ be
the diagonal (identity) matrix. For every $i=1,2,\ldots,n$, let $\hat
A_i$ be the matrix obtained from $A- I$ by zeroing out all but the
$i$th column, and let $A_i = I + \hat A_i$. Consider this collection of
matrices $ \{ A_1, \ldots, A_n \}$. Suppose we generate a
sequence of i.i.d. random variables $\{ I_1, I_2, \ldots\}$ with the
uniform distribution in $V$ and consider the sequence of products of
i.i.d. matrices $A_{I_n} A_{I_{n-1}}\cdots A_{I_2}A_{I_1}$.
It is easy to see that $M_t^x$ is the $x$th coordinate of the product
of a finite sequence of matrices $A_{I_j}$ right-multiplied with the
column vector $M^x_0$.
Does the product $A_{I_n} A_{I_{n-1}}\cdots A_{I_2}A_{I_1}$ have a
limit in some sense as $n$ tends to infinity?
Products of i.i.d. random matrices have been an old and fascinating
subject (see \cite{FK60,DF99}), and several conditions are known for
convergence of distributions of the products. There are also a number
of theorems on the limit distribution. However, the particular class of
products considered here is just beyond the assumptions under which
general results are known to hold. Most of the entries in any $A_i$ are
zero, violating assumptions in \cite{FK60}, equation~(3.1), and these
matrices are not strong contractions in the sense of \cite{DF99},
Theorem~1.1. However, as we will show, the additional graph structure in
the background determines the limit and the rate of convergence of the
products $A_{I_n} A_{I_{n-1}}\cdots A_{I_2}A_{I_1}$.

\item[(iii)] A more recent line of investigation related to our work is on
Markov chains on the space of partitions; see \cite{CL,C14}. One of
the important considerations in this regard is the product of i.i.d.
picks from a probability measure on the space of finite probability
transition matrices. That is, $S_1, S_2, \ldots$ are i.i.d. stochastic
matrices, and one is interested in the \textit{backward} product $S_1
S_2 \cdots S_m$. The knowledge of this limit determines the behavior of
a corresponding Markov chain on the space of finite partitions of
$\mathbb{N}$; see \cite{C14}, equation (5) and Theorem~1.2. The
transpose of each $A_i$ is a stochastic matrix. If we define
$S_i=A_i'$, then $S_1 S_2 \cdots S_m$ is the transpose of $A_m \cdots
A_2 A_1$. Hence, the limits described in this work give explicit
information about certain Markov chains on the space of partitions.
\end{longlist}

In the title of this paper, WIMPs stands for ``weakly interacting
mathematical particles.'' It turns out that one of the main technical
tools in this paper is a pair of ``weakly interacting'' continuous time
symmetric random walks on the graph. If the two random walks are at
different vertices, they move independently. However, if they are at
the same vertex, their next jumps occur at the same time, after an
exponential waiting time, common to both processes. The dependence ends
here---the two processes jump to vertices chosen independently, even
though they jump at the same time. Heuristically, one expects WIMPs to
behave very much like independent continuous time random walks. Showing
this is the heart of a number of arguments but it proves to be harder
than one would expect. In other cases, the slight dependence manifests
itself clearly and generates phenomena that otherwise would be trivial.
WIMPs played an important role in \cite{FF}.

The rest of the paper is organized as follows. Section~\ref{secpre}
contains rigorous definitions of the meteor process and WIMPs.
Section~\ref{secconv} is devoted to basic properties of the meteor
process and convergence to equilibrium. We present three theorems on
convergence. The first one is very abstract and does not provide any
useful information on the rate of convergence. The second one provides
a rate of convergence, but since it applies to all meteor processes, it
is not optimal in specific examples. The third theorem is limited to
tori and gives a sharp estimate for the convergence rate.
Section~\ref{craters} is devoted to the distribution of craters in
circular graphs and is the most combinatorial of all the sections---it
is partly based on results from \cite{BBS}.
We address several questions about craters. The first one is concerned
with the probability of a given pattern of craters.
The second question is about fluctuations of the numbers of patterns
around the mean. We do not provide a standard large deviations result,
but we prove a theorem on the most likely configuration of craters
assuming that there are very few of them.
Section~\ref{secmass} presents results on the mass distribution at a
single vertex or a
family
of vertices, in case of circular graphs.
Section~\ref{secnonc} contains theorems on the mass distribution for
noncircular graphs. The first result is a bound for the variance for a
large class of graphs. The second result is a completely explicit
limiting distribution at a vertex, for the complete graphs, when the
size of the graph grows to infinity. Finally,
Section~\ref{earth} contains the proof of the claim that the earthworm
distributes mass in a torus more or less evenly, on a large scale.

\section{Preliminaries}\label{secpre}

The following setup and notation will be used in most of the
paper. All constants will be assumed to be strictly positive, finite,
real numbers, unless stated otherwise. The notation $|S|$ will be
used for the cardinality of a finite set, $S$.

We will consider only finite connected graphs with no loops and no
multiple edges. We will often denote the chosen graph by $G$ and its
vertex set by $V$. In particular, we often use $k$ for $|V|$. We let
$d_v$ stand for the degree of a vertex $v$, and write
$v\leftrightarrow x$ if
vertices $v$ and $x$ are connected by an edge.

We will write $\calC_k$ to denote the circular graph with $k$
vertices, $k\geq2$. In other words, the vertex set of $\calC_k$ is
$\{1,2,\ldots, k\}$, and the only pairs of vertices joined by edges are
of the form $(j, j+1)$ for $j=1,2, \ldots, k-1$, and $(k,1)$. For
$\calC_k$, all arguments will apply ``mod $k$.'' For example, we will
refer to $k$ as a vertex ``to the left of 1,'' and interpret $j-1$ as
$k$ in the case when $j=1$.

Every vertex $v$ is associated with a Poisson process $N^v$
representing ``arrival times of meteors''
with intensity 1.
We assume that processes $N^v$ are jointly independent. A vertex $v$
holds some ``soil'' with mass equal to $M^v_t\geq0$ at time $t\geq0$.
The processes $M^v$ evolve according to the following scheme.

We assume that $M^v_0 \in[0, \infty)$ for every $v$, a.s.
At the time $t$ of a jump of $N^v$, $M^v$ jumps to 0. At the same time,
the mass $M^v_{t-}$ is ``distributed'' equally among all adjacent
sites, that is, for every vertex $x\leftrightarrow v$,
the process $M^x$ increases by $M^v_{t-}/d_v$; more formally,
$M^x_t = M^x_{t-} + M^v_{t-}/d_v$.
The mass $M^v$ will change only when $N^v$ jumps
and just prior to that time there is positive mass at $v$,
or $N^x$ jumps for
some $x\leftrightarrow v$ and just prior to that time there is
positive mass
at $x$.
We will denote the mass process $\mathcal{M}_t = \{M^{v}_t, v\in V\}$.

We will now define WIMPs (``weakly interacting mathematical
particles'') which will be used in a number of arguments.

%
\begin{definition}\label{wimps}
We will define several processes on the same probability space.
Suppose that two mass processes $\mathcal{M}_0$ and $\widetilde{\mathcal
{M}}_0$ are given, and assume that $a= \sum_{v\in V} M^v_0 =
\sum_{v\in V} \widetilde M^v_0$.

For each $j\geq1$, let $\{Y^j_n, n\geq0\}$ be a discrete time
symmetric random walk on $G$ with the initial distribution
$\mathbb{P}(Y^j_0 = x) = M^x_0/a$ for $x\in V$.
Similarly, let $\{\widetilde Y^j_n, n\geq0\}$, $j\geq1$, be discrete
time symmetric random walks on $G$ with the initial distribution
$\mathbb{P}(\widetilde Y^j_0 = x) = \widetilde M^x_0/a$ for $x\in V$,
$j\geq1$.
We assume that conditional on $ \mathcal{M}_0$ and $\widetilde{\mathcal
{M}}_0$, all processes $\{Y^j_n, n\geq0\}$, $j\geq1$ and $\{
\widetilde Y^j_n, n\geq0\}$, $j\geq1$, are independent.

Recall Poisson processes $N^v$ defined earlier in this section,
and assume that they are independent of $\{Y^j_n, n\geq0\}$, $j\geq1$
and $\{\widetilde Y^j_n, n\geq0\}$, $j\geq1$.
For every $j\geq1$, we define a continuous time Markov process $\{
Z^j_t, t\geq0\}$
by\vspace*{1pt} requiring that the embedded discrete time Markov chain for $Z^j$ is
$Y^j$ and $Z^j$ jumps\vspace*{1pt} at a time $t$ if and only if
$N^v$ jumps at time $t$, where $v=Z^j_{t-}$.
We define $\{\widetilde Z^j_t, t\geq0\}$
in an\vspace*{1pt} analogous way, so that the embedded discrete time Markov chain
for $\widetilde Z^j$ is $\widetilde Y^j$ and $\widetilde Z^j$ jumps at
a time $t$ if and only if
$N^v$ jumps at time $t$, where $v=\widetilde Z^j_{t-}$.
Note that the jump times of all $Z^j$'s and $\widetilde Z^j$'s are
defined by the same family of Poisson processes $\{N^v\}_{v\in V}$.
\end{definition}

%
\begin{remark}\label{remwimps}
The processes $Z^j$ and $\widetilde Z^j$ in Definition~\ref{wimps} are
continuous time nearest neighbor symmetric random walks on $G$ with
exponential holding time with mean 1.

The joint distribution of $(Z^1, Z^2)$ is the following. The state
space for the process $(Z^1, Z^2)$ is $V^2$. If $(Z^1_t, Z^2_t)=(x,y)$
with $x\ne y$, then the process will stay in this state for an
exponential amount of time with mean $1/2$, and at the end of this time
interval, one of the two processes (chosen uniformly) will jump to one
of the nearest neighbors (also chosen uniformly). This behavior is the
same as that of two independent random walks. However, if
$(Z^1_t, Z^2_t)=(x,x)$, then the pair of processes behave in a way that
is different from that of a pair of independent random walks. Namely,
after an exponential waiting time with mean 1 (not $1/2$),
both processes will jump at the same time; each one will jump to one of
the nearest neighbors of $x$ chosen uniformly and independently of the
direction of the jump of the other process.

The same remark applies to any pair of processes
in the family
$\{Z^j, j\geq1\} \cup\{ \widetilde Z^j, j\geq1\}$.
\end{remark}

\section{Basic properties and convergence to equilibrium}\label{secconv}

It will be convenient from the technical point of view to postpone the
presentation of the most elementary properties of the meteor process to
the end of this section. We will start with three theorems on
convergence to the stationary distribution.

%
\begin{remark}\label{j121}
The mass process $\{\mathcal{M}_t, t\geq0\}$ is a somewhat unusual
stochastic process in that its state space can be split into an uncountable
number of disjoint communicating classes. It is easy to see that, due
to the
definition of the evolution of $\{\mathcal{M}_t, t\geq0\}$,
for every time $t\geq0$ and every $v\in V$,
$M^v_t = \sum_{x\in V} a_x M^x_0$, where $a_x$ is a random variable
that depends on
$t$ and $v$. Every $a_x$ has the form $m \prod_{y\in V}d_y^{-j_y}$ for
some\vspace*{1pt} integers
$m$ and $j_y$.
In other words, $a_x$ take values in the ring ${\mathbb Z}[1/e_1, \ldots,
1/e_{i}]$, where $e_{1},\ldots, e_{i}$ is a list of all distinct
values of degrees of vertices in $V$. Therefore, $M^v_t$'s take values
in the free module over the ring ${\mathbb Z}[1/e_1,\ldots, 1/e_{i}]$ spanned
by $\{M^v_ 0,v\in V\}$.

For example, consider the following two initial distributions. Suppose
that $M^v_0 = 1$ for all $v$. Fix some $x\in V$, and let $\widetilde
M^v_0 = 1/\pi$ for all $v\ne x$ and $\widetilde M^x_0 = |V| -
(|V|-1)/\pi$. If $\{\mathcal{M}_t, t\geq0\}$ and $\{\widetilde{\mathcal
{M}}_t, t\geq0\}$
are mass processes with these initial distributions, then
for every $t>0$, the distributions of $\mathcal{M}_t$ and $\widetilde
{\mathcal{M}}_t$
will be mutually singular.

It follows from these observations that proving convergence of $\{
\mathcal{M}_t, t\geq0\}$ to the stationary distribution cannot
proceed along the most classical lines; see \cite{HMS} for the
discussion of this technical issue and a solution. We will follow \cite
{HMS} in spirit although not in all technical details.
\end{remark}

%
\begin{theorem}\label{a297}
Consider the process $\{\mathcal{M}_t, t\geq0\}$ on a graph $G$.
Assume that $|V| =k$ and $\sum_{v\in V} M^v_0 = k$.
When $t\to\infty$, the distribution of
$\mathcal{M}_t$
converges to a distribution $Q$ on $[0,k]^k$. The distribution $Q$ is
the unique stationary distribution for the process
$\{\mathcal{M}_t, t\geq0\}$.
In particular, $Q$ is independent of the initial
distribution of $\mathcal{M}$.
\end{theorem}

\begin{pf}
We\vspace*{1.5pt} will consider a coupling of two copies of the process $\{\mathcal
{M}_t,\break t\geq0\}$. Suppose that $\{\mathcal{M}_t, t\geq0\}$ and $\{
\widetilde{\mathcal{M}}_t, t\geq0\}$
are driven by the same processes $\{N^v_t, t\geq0\}_{v\in V}$ but the
distribution of $\mathcal{M}_ 0$ is not\vspace*{1.5pt} necessarily the same as that
$\widetilde{\mathcal{M}}_ 0$. We do assume that $\sum_{v\in V} M^v_0 =
\sum_{v\in V} \widetilde M^v_0 = k$.

First we will argue that the total variation distance between the
distributions of $\mathcal{M}_t $ and $ \widetilde{\mathcal{M}}_t$,
that is, $D_t:= \sum_{v\in V} |M^v_t - \widetilde M^v_t|$ is a
nonincreasing process, a.s.
Note that since $G$ is finite, the number of jumps of $N^v$'s is finite
on every finite time interval and $D_t$ is constant between any two
jump times.
Suppose that $N^x$ has a jump at a time $T$. Then
\begin{eqnarray*}
D_T - D_{T-} &=& \sum_{v\in V}
\bigl\llvert M^v_T - \widetilde M^v_T
\bigr\rrvert- \sum_{v\in V} \bigl\llvert
M^v_{T-} - \widetilde M^v_{T-}\bigr
\rrvert
\\
&=& -\bigl\llvert M^x_{T-} - \widetilde
M^x_{T-}\bigr\rrvert+ \sum_{v\leftrightarrow x}
\bigl( \bigl\llvert M^v_T - \widetilde
M^v_T\bigr\rrvert- \bigl\llvert M^v_{T-}
- \widetilde M^v_{T-}\bigr\rrvert\bigr)
\\
&=& -\bigl\llvert M^x_{T-} - \widetilde
M^x_{T-}\bigr\rrvert
\\
&&{} + \sum_{v\leftrightarrow x}
\biggl( \biggl\llvert\bigl(M^v_{T-} - \widetilde
M^v_{T-}\bigr) + \frac{M^x_{T-} - \widetilde M^x_{T-}}{d_x}\biggr
\rrvert- \bigl
\llvert M^v_{T-} - \widetilde M^v_{T-}
\bigr\rrvert\biggr)
\\
&\leq&-\bigl\llvert M^x_{T-} - \widetilde
M^x_{T-}\bigr\rrvert
\\
&&{} + \sum_{v\leftrightarrow x}
\biggl( \bigl\llvert M^v_{T-} - \widetilde
M^v_{T-}\bigr\rrvert+ \biggl\llvert\frac{M^x_{T-} - \widetilde M^x_{T-}}{d_x}
\biggr\rrvert- \bigl\llvert M^v_{T-} - \widetilde
M^v_{T-}\bigr\rrvert\biggr)
\\
&=& -\bigl\llvert M^x_{T-} - \widetilde
M^x_{T-}\bigr\rrvert+ \sum_{v\leftrightarrow x}
\biggl\llvert\frac{M^x_{T-} - \widetilde
M^x_{T-}}{d_x}\biggr\rrvert=0.
\end{eqnarray*}
This shows that $D_t$ is nonincreasing.

Recall that $G$ is connected, and fix some vertex $y$. Let $(x_1, x_2,\ldots, x_n)$ be a sequence of vertices of $G$ such that $x_j
\leftrightarrow x_{j+1}$ for $1\leq j \leq n-1$, $x_n \leftrightarrow
y$ and $\{x_1, x_2,\ldots, x_n\} = G \setminus\{y\}$.
The vertices' $x_j$'s are not necessarily distinct.
Recall that $|V|=k$.
Let $d = \max_{x\in V} d_x$.
Let $a = \max_{v\ne y} M^v_0$, $b = \sum_{v\ne y} M^v_0$, and note
that $a\geq b/(k-1)$.
Suppose that the first $n$ meteors hit vertices $x_1, x_2,\ldots,
x_n$, in this order.
During this process, at least $1/d$th part of the mass from any vertex
$x_j$, $j<n$, is pushed to $x_{j+1}$, and at least $1/d$th part of the
mass at $x_n$ is pushed to $y$.
Let $m$ be the smallest integer with the property that $M^{x_m}_0 =a$.
Then at least $ad^{-1}$ of mass will be pushed from $x_m$ to $x_{m+1}$.
This implies that
least $ad^{-2}$ of mass will be pushed from $x_{m+1}$ to $x_{m+2}$. By
induction, at least $ad^{-j}$ of mass will be pushed from $x_{m+j-1}$
to $x_{m+j}$.
Hence, at least $a d^{-n}$ of mass will be added to $y$.
In other words, the mass outside $y$ will be reduced at least by $a
d^{-n} \geq b d^{-n}/(k-1)$. Putting it in a different way, the mass
outside $y$ will be reduced at least by the factor of $ 1-
d^{-n}/(k-1)$.

Consider an arbitrarily small $\varepsilon>0$, and let $m$ be so large
that $(1- d^{-n}/(k-1))^m k\leq\varepsilon$.
If the first $n m$ meteors hit vertices
\[
\underbrace{x_1, x_2,\ldots, x_n,\qquad
x_1, x_2,\ldots, x_n, \qquad
\ldots,\qquad
x_1, x_2,\ldots, x_n}
_{m\ \mathrm{times}},
\]
in this order, then the mass outside $y$ will be reduced to at most
$(1- d^{-n}/(k-1))^m k\leq\varepsilon$. Sooner or later, with
probability 1, there will be a sequence of $nm$ meteor hits described
above, and then the mass outside $y$ will be less than $\varepsilon$.
Hence the mass at $y$ will be between $k-\varepsilon$ and $k$ at the
end of this sequence\vspace*{2pt} of meteor hits. Note that the argument applies
equally to
$\{\mathcal{M}_t, t\geq0\}$ and $\{\widetilde{\mathcal{M}}_t, t\geq
0\}$. Hence, at the end of this sequence of $nm$ hits, the function $D$
will be at most $2\varepsilon$. Since $D_t$ is nonincreasing, we see
that $D_t$ converges to 0, a.s.

For every $t$, the distribution of $\mathcal{M}_t$ is a measure on
$[0,k]^k$, a compact set, so the family of distributions of $\mathcal
{M}_t$, $t\geq0$, is tight. Therefore, there exists a sequence $t_n$
converging to $\infty$ such that the distributions of $\mathcal
{M}_{t_n}$ converge to a distribution $Q$ on $[0,k]^k$, as $n\to\infty$.

Let ${\mathbf d}$ denote the Prokhorov distance (see \cite{Bill},
page~238) between probability measures on $[0,k]^k$, and recall
that convergence in the metric ${\mathbf d}$ is equivalent to the weak
convergence of measures.
By abuse of notation, we will use the same symbol for the Prokhorov
distance between probability measures on $[0,k]^k$ and $\mathbb{R}$.
We will also apply ${\mathbf d}$ to random variables, with the
understanding that it applies to their distributions. Let ${\mathbf0}$
denote the probability distribution on $\mathbb{R}$ concentrated at $0$.
It is easy to see that for every $\delta>0$ there exists $\alpha
(\delta)>0$ such that if ${\mathbf d}(D_t, {\mathbf0}) \leq\alpha
(\delta)$ then ${\mathbf d}(\mathcal{M}_t, \widetilde{\mathcal{M}}_t)
\leq\delta$.

The bounds in our argument showing convergence of $D_t$ to 0 do not
depend on $\mathcal{M}_ 0$ or $\widetilde{\mathcal{M}}_ 0$ so
there exists a deterministic\vspace*{1pt} function $\rho\dvtx  [0,\infty)\to[0,\infty
)$ such that $\lim_{t\to\infty} \rho(t) =0$ and
${\mathbf d}(D_t, {\mathbf0}) \leq\rho(t)$ for any $\mathcal{M}_ 0$
any $\widetilde{\mathcal{M}}_ 0$.

Suppose that there exists a sequence $s_n$ converging to $\infty$ such
that the distributions of $\mathcal{M}_{s_n}$ converge\vspace*{1pt} to a
distribution $Q'$ on $[0,k]^k$, as $n\to\infty$, and $Q' \ne Q$.
Let $\delta= {\mathbf d}(Q, Q') / 2>0$.

Find $u_0$ so large that $\rho(t) < \alpha(\delta)$ for $t\geq u_0$.
Let $t_n$ and $s_m$ be such that $u_0 < t_n < s_m$,
${\mathbf d}(\mathcal{M}_{t_n}, Q) \leq\delta/4$ and ${\mathbf
d}(\mathcal{M}_{s_m}, Q') \leq\delta/4$.
Let $\widetilde{\mathcal{M}}_ 0=\mathcal{M}_ {s_m-t_n}$.\vspace*{2pt} Then
${\mathbf d}(\widetilde{\mathcal{M}}_{t_n}, Q') \leq\delta/4$.
Since $t_n > u_0$, we have $\rho(t_n) < \alpha(\delta)$, so
${\mathbf d}(D_{t_n}, {\mathbf0}) \leq\alpha(\delta)$ and, therefore,
${\mathbf d}(\mathcal{M}_{t_n}, \widetilde{\mathcal{M}}_{t_n}) \leq
\delta$.
By the triangle inequality,
\begin{eqnarray*}
{\mathbf d}\bigl(Q, Q'\bigr) &\leq&{\mathbf d}(
\mathcal{M}_{t_n}, Q) + {\mathbf d}\bigl(\widetilde{
\mathcal{M}}_{t_n}, Q'\bigr) + {\mathbf d}(
\mathcal{M}_{t_n}, \widetilde{\mathcal{M}}_{t_n})
\\
&\leq&\delta/4 + \delta/4 + \delta= 3\delta/2.
\end{eqnarray*}
This contradicts the fact that
${\mathbf d}(Q, Q') = 2 \delta$ and shows that
$\mathcal{M}_t$ converges in distribution to $Q$, as $t\to\infty$.
The fact that $D_t$ converges to 0 shows that $Q$ does not depend on
the distribution of $\mathcal{M}_0$.

Next we will show that the distribution $Q$ is stationary.
It is routine to show that for every $\eta>0$ there exists $\beta
(\eta)>0$ such that
for any distributions $Q$~and~$Q'$ on $[0,k]^k$ with ${\mathbf
d}(Q',Q'')\leq\beta(\eta)$, one can construct $\mathcal{M}_0$ and
$\widetilde{\mathcal{M}}_0$
on the same probability space so that the distribution of $\mathcal
{M}_0$ is~$Q'$, the distribution of $\widetilde{\mathcal{M}}_0$ is~$Q''$, and ${\mathbf d}(D_0, {\mathbf0}) \leq\eta$.

Consider an arbitrarily small $\delta>0$.
Let the distribution of $\mathcal{M}_ 0$ be $Q$ and
find $u_1 $ so large that
${\mathbf d}(\mathcal{M}_{t}, Q) \leq\beta(\alpha(\delta/2)) \land
\delta/2$ for all $t\geq u_1$.
Then we can construct $\mathcal{M}_0$ and $\widetilde{\mathcal{M}}_0$
on the same probability space so that the distribution of $\mathcal
{M}_0$ is $Q$, the distribution of $\widetilde{\mathcal{M}}_0$ is the
same as that of $\mathcal{M}_{u_1}$, and ${\mathbf d}(D_0, {\mathbf0})
\leq\alpha(\delta/2)$.\vspace*{1pt}
Then ${\mathbf d}(D_t, {\mathbf0}) \leq\alpha(\delta/2)$ for all
$t\geq0$ and, therefore, we have
${\mathbf d}(\mathcal{M}_t, \widetilde{\mathcal{M}}_t)
\leq\delta/2$ for any $t\geq0$.
Note that ${\mathbf d}(\widetilde{\mathcal{M}}_{t}, Q) \leq\delta/2$
for all $t\geq0$
because ${\mathbf d}(\mathcal{M}_{t}, Q) \leq\delta/2$ for all
$t\geq u_1$.
We obtain for $t\geq0$,
\begin{eqnarray*}
{\mathbf d}(\mathcal{M}_t, Q) &\leq&{\mathbf d}(\widetilde{
\mathcal{M}}_{t}, Q) + {\mathbf d}(\mathcal{M}_t,
\widetilde{\mathcal{M}}_t) \leq\delta/2 + \delta/2 = \delta.
\end{eqnarray*}
Since $\delta>0$ is arbitrarily small, $Q$ is stationary.
\end{pf}

%
\begin{remark}\label{j132}
In view of Theorem~\ref{a297} and its proof, it is easy to see that
there exists a stationary version of the process $\mathcal{M}_t$ on
the whole real line; that is, there exists a process $\{\mathcal{M}_t,
t\in\mathbb{R}\}$, such that the distribution of $\mathcal{M}_t$ is
the stationary measure $Q$ for each $t\in\mathbb{R}$. Moreover, one
can construct independent Poisson processes $\{N^v_t,t\in\mathbb
{R}\} $, $v\in V$, on the same probability space, such that $\{\mathcal
{M}_t, t\in\mathbb{R}\}$ jumps according to the algorithm described
in Section~\ref{secpre}, relative to these Poisson processes. We set
$N^v_0=0$ for all $v$ for definiteness.
\end{remark}

The next theorem is the only result in our paper that
is proved in a context more general than that in Section~\ref{secpre}.
Consider a graph, and let $\bP= (p_{xy})_{x,y\in V}$ be the
probability transition matrix for a Markov chain on $V$. In this model,
if a meteor hits site $x$, then the mass is distributed to other
vertices in proportion to $p_{xy}$, not necessarily in equal
proportions to all neighbors. We remark parenthetically that, by
convention, we place an edge between two vertices $x$ and $y$ of $G$ if
and only if $p_{xy} + p_{yx} >0$.

%
\begin{theorem}\label{f101}
Consider a graph $G$, and suppose that $U_t$ and $\widetilde U_t$ are
independent continuous time Markov chains with mean 1 exponential
holding times at every vertex and the transition rates
for the embedded discrete time Markov chains given by $\bP$.
Let
%
%
\begin{eqnarray}\label{f122}
\tau_U&=& \inf\{t\geq0\dvtx  U_t =\widetilde
U_t\},
\nonumber\\[-8pt]\\[-8pt]
\alpha(t) &=& \sup_{x,y\in V} \mathbb{P}(\tau_U > t
\mid U_0 = x, \widetilde U_0 = y).\nonumber
\end{eqnarray}
Consider any (possibly random) distributions of mass $ \mathcal{M}_0$
and $\widetilde{\mathcal{M}}_0$; that is, assume that $ M^x_0 \geq0$
and $\widetilde M^x_0 \geq0$
for all $x\in V$ and $\sum_{x\in V} M^x_0 =\sum_{x\in V} \widetilde
M^x_0 = |V|$, a.s.
One can define mass processes $\mathcal{M}_t$
and $\widetilde{\mathcal{M}}_t$ on a common probability space so that
for all $t\geq0$,
%
%
\begin{equation}
\label{f91.3.2} \mathbb{E} \biggl( \sum_{x\in V} \bigl
\llvert M^x_t - \widetilde M^x_t
\bigr\rrvert\biggr) \leq|V| \alpha(t).
\end{equation}
\end{theorem}

%
\begin{remark}
In the notation of Theorem~\ref{f101},
let $T_x = \inf\{t\geq0\dvtx  U_t = x\}$. According to \cite{Ald},
Proposition~1, if $\bP$ represents a reversible Markov chain, then
\[
\sup_{x,y\in V} \mathbb{E}( \tau_U \mid
U_0 = x, \widetilde U_0 = y) \leq c \sup
_{x,y\in V} \mathbb{E}( T_x \mid U_0 = y).
\]
For an arbitrary Markov chain, Conjecture 1 in \cite{Ald}
states that
\[
\sup_{x,y\in V} \mathbb{E}( \tau_U \mid
U_0 = x, \widetilde U_0 = y) \leq c |V| \sup
_{x,y\in V} \mathbb{E}( T_x \mid U_0 = y).
\]
The conjecture remains open at this time, as far as we know.
\end{remark}

\begin{pf*}{Proof of Theorem~\ref{f101}}
Suppose that $Z$~and~$\widetilde Z$ are constructed as $Z^1$~and~$\widetilde Z^1$ in Definition~\ref{wimps}, except that
$\{Y^1_n, n\geq0\}$ and $\{\widetilde Y^1_n, n\geq0\}$ are discrete
time Markov chains with the transition probability matrix $\bP$. The
initial distributions are given by
$\mathbb{P}(Y_0 = x) = M^x_0/|V|$ and $\mathbb{P}(\widetilde Y_0 =
x) = \widetilde M^x_0/|V|$ for $x\in V$.
Let
\begin{eqnarray*}
\tau&=& \inf\{t\geq0\dvtx  Z_t = \widetilde Z_t\},
\\
\widehat Z_t &=& \cases{ \widetilde Z_t, &\quad for $t
\leq\tau$,
\vspace*{3pt}\cr
Z_t, &\quad for $t > \tau$.}
\end{eqnarray*}
The distribution of
$\{\widehat Z_t, t\geq0\}$
is the same as that of
$\{\widetilde Z_t, t\geq0\}$.

Let $\{Z^*_t, t\geq0\}$ have the same distribution as $\{\widetilde
Z_t, t\geq0\}$ and be independent of $\{Z_t, t\geq0\}$, given
$\mathcal{M}_0$ and $\widetilde{\mathcal{M}}_0$. Let $\tau_* = \inf\{
t\geq0\dvtx  Z_t = Z^*_t\}$.
Since the Poisson processes $N^m$ are independent from one another, it
follows easily that the distributions of
\[
\bigl\{\tau, \bigl\{Z_t, t\in[0,\tau]\bigr\}, \bigl\{\widetilde
Z_t, t\in[0,\tau]\bigr\} \bigr\}
\]
and
\[
\bigl\{
\tau_*, \bigl\{Z_t, t\in[0,\tau_*]\bigr\}, \bigl\{
Z^*_t, t\in[0,\tau_*]\bigr\} \bigr\}
\]
are identical. Thus $\tau$ and $\tau_*$ have the same distributions,
and therefore, (\ref{f122}) implies that $\mathbb{P}(\tau>t) \leq
\alpha(t)$.

Let $\mathcal{G}_t = \sigma(\mathcal{M}_s, \widetilde{\mathcal{M}}_s,
0\leq s \leq t)$, and note that $\mathcal{G}_t = \sigma
(\mathcal{M}_0, \widetilde{\mathcal{M}}_0, N^v_s, 0\leq s \leq t, v\in V)$.
The process $Z$
is ``coupled'' with the
processes $N^v$ which determine the motion of mass.
This easily implies that for all $x$ and $t$,
%
%
\begin{equation}
\label{f121} \mathbb{P}(Z_t =x \mid\mathcal{G}_t) =
M^x_t/|V|.
\end{equation}
It is easy to see that the distributions of
\[
\bigl\{ \bigl\{N^v_s, s\in[0,t]\bigr
\}_{v\in V}, \bigl\{\widetilde Z_s, s\in[0,t]\bigr\} \bigr\}
\]
and
\[
\bigl\{ \bigl\{N^v_s, s\in[0,t]
\bigr\}_{v\in V}, \bigl\{\widehat Z_s, s\in[0,t]\bigr\} \bigr
\}
\]
are the same, so we obtain the following formula, analogous to (\ref{f121}),
\[
\mathbb{P}(\widetilde Z_t =x \mid\mathcal{G}_t) =
\mathbb{P}(\widehat Z_t =x \mid\mathcal{G}_t) =
\widetilde M^x_t/|V|.
\]
It follows that
\begin{eqnarray*}
\mathbb{E} \biggl( \sum_{x\in V} \bigl\llvert
M^x_t - \widetilde M^x_t\bigr
\rrvert\biggr) &=& |V|\mathbb{E} \biggl( \sum_{x\in V}
\bigl\llvert\mathbb{P}(Z_t =x \mid\mathcal{G}_t) -
\mathbb{P}(\widehat Z_t =x \mid\mathcal{G}_t)\bigr
\rrvert\biggr)
\\
&=& |V|\mathbb{E} \biggl( \sum_{x\in V} \bigl\llvert
\mathbb{E} (\bone_{\{Z_t =x\}} - \bone_{\{\widehat Z_t =x\}} \mid
\mathcal{G}_t ) \bigr\rrvert\biggr)
\\
&\leq&|V|\mathbb{E} \biggl( \sum_{x\in V} \mathbb{E}
\bigl(\llvert\bone_{\{Z_t =x\}} - \bone_{\{\widehat Z_t =x\}} \rrvert
\mid\mathcal
{G}_t \bigr) \biggr)
\\
&=& |V| \mathbb{E}\mathbb{E} (\bone_{\{ Z_t \ne\widehat Z_t\}} \mid
\mathcal{G}_t
)
\\
&=& |V| \mathbb{E}\mathbb{P} (\tau> t \mid\mathcal{G}_t ) \leq|V|
\alpha(t).
\end{eqnarray*}
This completes the proof.
\end{pf*}

%
\begin{theorem}\label{ma61}
Consider the meteor process on a graph $G=\calC_n^d$ (the product of
$d$ copies of the cycle $\calC_n$).
Consider any distributions (possibly random) of mass $ \mathcal{M}_0$
and $\widetilde{\mathcal{M}}_0$; that is, assume that $ M^x_0 \geq0$
and $\widetilde M^x_0 \geq0$
for all $x\in V$ and $\sum_x M^x_0 = \sum_x \widetilde M^x_0 = |V| =
n^d$, a.s.
There exist constants $c_1, c_2$ and $c_3 $, not depending on $n$ and
$d$, such that\vspace*{1pt}
if $n \geq1 \lor c_1 \sqrt{d \log d}$ and $t \geq c_2 d n^2 $, then
one can define a coupling of processes $\mathcal{M}_t$
and $\widetilde{\mathcal{M}}_t$ on a common probability space so that
%
%
\begin{equation}
\label{f91} \mathbb{E} \biggl(\sum_{x\in V}
\bigl|M^x_t - \widetilde M^x_t\bigr|
\biggr) \leq\exp\bigl(- c_3 t / \bigl(dn^2\bigr)\bigr)
|V|.
\end{equation}
\end{theorem}

\begin{pf}
\textit{Step}~1.
In this step, we will show that
there exist constants $c_1, c_2\in(0,\infty)$ and $c_4 < 2$, not
depending on $n$ and $d$, such that
if $n \geq1 \lor c_1 \sqrt{d \log d}$ and $t \geq c_2 d n^2 $, and
the processes $\mathcal{M}_t$
and $\widetilde{\mathcal{M}}_t$ are independent, then
%
%
\begin{equation}
\label{ma62} \mathbb{E} \biggl(\sum_{x\in V}
\bigl|M^x_t - \widetilde M^x_t\bigr|
\biggr) \leq c_4 |V|.
\end{equation}

Let $Z$ and $\widetilde Z$ be defined as $Z^1$ and $ Z^2$ in Definition
\ref{wimps}.
In particular,
$\mathbb{P}(Z_0 = x) = \mathbb{P}(\widetilde Z_0 = x) = M^x_0/|V|$
for $x\in V$.

Let $Z^*_t = Z_t - \widetilde Z_t$, and note that $Z^*$ is a continuous
time Markov process on~$V$, with the mean holding time equal to $1/2$
at all vertices $x\ne\bzero:=(0,\ldots,0)$. Recall that if
$(Z_t,\widetilde Z_t)=(x,x)$, then after an exponential waiting time
with mean 1 (not~$1/2$),
both processes will jump at the same time. They will jump to one of the
neighbors of $x$ (the same for both processes) with probability
$1/(2d)$. Hence, this jump of $(Z,\widetilde Z)$ will not correspond to
a jump of $Z^*$. It follows that the mean holding time for $Z^*$ at
$\bzero$ is
$\beta:=(1- 1/(2d))^{-1}$. Note that if $Z^*_t=\bzero$, the next jump
it will take will be to a vertex at the distance 2 from $\bzero$. If
$Z^*_t\ne\bzero$, then the next jump
will be to a neighbor of $Z^*_t$.

Let $Z^1_t$ be a continuous time symmetric nearest neighbor random walk
on $V$, with the mean holding time equal to $1/2$ at all vertices $x\ne
\bzero$, and mean holding time at $\bzero$ equal to $\beta$.
The only difference between $Z^1$ and $Z^*$ is that $Z^1$ can jump from
$\bzero$ only to a nearest neighbor while $Z^*$ can jump from $\bzero
$ to some other vertices.

We will construct a coupling of $Z^*$ and $Z^1$ such that $Z^1_0 =
Z^*_0$ and, a.s.,
%
%
\begin{equation}
\label{ma72} \bigl\{t\geq0\dvtx  Z^*_t =\bzero\bigr\} \subset\bigl\{t
\geq0\dvtx  Z^1_t =\bzero\bigr\}.
\end{equation}

We let $Z^1_t = Z^*_t$ for all $t$ less than the time $S_1$ of the
first jump out of $\bzero$.
At the time $S_1$, we let processes $Z^1$ and $Z^*$ make independent
jumps, each one according to its own jump distribution.

Let $T_{Z^1}(\bzero) =\inf\{t\geq0\dvtx  Z^1_t = \bzero\}$, and let
$T_{Z^*}(\bzero)$ have the analogous meaning. Suppose that $x,y \in
V$, $x\leftrightarrow\bzero$ and $y\nleftrightarrow\bzero$. Then
for every $t\geq0$,
\[
\mathbb{P}^x\bigl(T_{Z^1}(\bzero) > t\bigr) \leq
\mathbb{P}^y\bigl(T_{Z^*}(\bzero) > t\bigr)
\]
because $Z^*$ has to pass a neighbor of $\bzero$ on its way to $\bzero
$. Now standard coupling arguments show that we can construct $Z^1$
after $S_1$ in such a way that it hits $\bzero$ at the same time or
earlier than
the hitting time of $\bzero$ by $Z^*$.
Let $S_2$ be the first hitting time of $\bzero$ by $Z^1$ after time
$S_1$. We will consider several cases:
\begin{longlist}[(a)]
\item[(a)] Suppose that $Z^*_{S_2}\nleftrightarrow\bzero$.
We let processes $Z^1$ and $Z^*$ evolve independently after $S_2$ until
the first time $S_3$ such that either
$Z^1_{S_3}\ne\bzero$ or $Z^*_{S_3}\leftrightarrow\bzero$.
\begin{enumerate}[(a2)]
\item[(a1)] Suppose that $Z^1_{S_3}\ne\bzero$. Then
$Z^1_{S_3}\leftrightarrow\bzero$ and
$Z^*_{S_3}\nleftrightarrow\bzero$. Hence, we can couple $Z^1$ and
$Z^*$ after time $S_3$ in such a way that $Z^1$ will hit $\bzero$
before $Z^*$ does. At the time when $Z^1$ hits $\bzero$,
we are back in the case represented by the time $S_2$.

\item[(a2)] Suppose that $Z^*_{S_3}\leftrightarrow\bzero$. Then
$Z^1_{S_3}=\bzero$. We\vspace*{1pt} continue the construction of the processes
after $S_3$ as in case (b) described below.
\end{enumerate}

\item[(b)] Suppose that $Z^*_{S_2}\leftrightarrow\bzero$.
We let processes $Z^1$ and $Z^*$ evolve independently after $S_2$ until
the first time $S_4$ such that
either $Z^1$ or $Z^*$ jumps.
\begin{enumerate}[(b2)]
\item[(b1)] If $Z^*$ jumps at time $S_4$ and
$Z^*_{S_4}\nleftrightarrow\bzero$, then we are back in the case
analogous to (a).

\item[(b2)] If $Z^*$ jumps at time $S_4$ and
$Z^*_{S_4}=\bzero$, then we continue in the same way as after time $0$.

\item[(b3)] If $Z^1$ jumps at time $S_4$, then
$Z^*_{S_4}\leftrightarrow\bzero$ and $Z^1_{S_4}\leftrightarrow
\bzero$. Then we couple $Z^1$ and $Z^*$ after
$S_4$ so that they hit $\bzero$ at the same time. We continue after
this time in the same way as after time $0$.
\end{enumerate}

\item[(c)] Suppose that $Z^*_{S_2}=\bzero$. Then we continue after this time
in the same way as after time $0$.
\end{longlist}

The construction of $Z^1$ can be continued by induction.
This completes the argument justifying the existence of a coupling of
$Z^1$ and $Z^*$ such that $Z^1$ is at $\bzero$ whenever $Z^*$ is at
this point.

It is elementary to check that for some $c_3 $ and all $d$, $n\geq2$,
$j\geq c_3 d n^2$ and $x\in V$, we have
%
%
\begin{equation}
\label{ma21} n^{-d}/2 \leq\mathbb{P}(Y_j \in x)
\leq2n^{-d}.
\end{equation}

Let\vspace*{1pt} $N^*$ be a Poisson process with the mean time between jumps equal
to $\beta$. It is easy to see that there exists $c_4 >0$ such that for
$t \geq2 \beta c_3 d n^2 $,
%
%
\begin{equation}
\label{ma22} \mathbb{P}\bigl(N^*_t \leq c_3 d
n^2\bigr) \leq e^{-c_4 n^2}.
\end{equation}

Let $\widetilde N_t$ be the number of jumps made by $Z^1$ by the time
$t$, and note that $\widetilde N$ is stochastically\vspace*{1pt} minorized by $N^*$.
By (\ref{ma72}), (\ref{ma21}) and (\ref{ma22}), there are $c_5$
and $c_6$ such that for $n \geq c_5 \sqrt{d \log d}$ and $t \geq2
\beta c_3 d n^2 $,
%
%
\begin{eqnarray}
\label{ma32} \mathbb{P}\bigl(Z^*_t = \bzero\bigr) &\leq&\mathbb{P}
\bigl(Z^1_t = \bzero\bigr)\nonumber
\\
&=& \sum_{j=0}^\infty\mathbb{P}
\bigl(Z^1_t = \bzero\mid\widetilde N_t = j
\bigr) \mathbb{P}(\widetilde N_t = j)
\nonumber
\\
&\leq&\mathbb{P}\bigl(\widetilde N_t \leq c_3 d
n^2\bigr) + \sum_{j> c_3 d n^2} \mathbb{P}
\bigl(Z^1_t = \bzero\mid\widetilde N_t = j
\bigr) \mathbb{P}(\widetilde N_t = j)
\nonumber\\[-8pt]\\[-8pt]
&\leq&\mathbb{P}\bigl( N^*_t \leq c_3 d
n^2\bigr) + \sum_{j> c_3 d n^2} 2n^{-d}
\mathbb{P}(\widetilde N_t = j)
\nonumber
\\
&\leq& e^{-c_4 n^2} + 2 n^{-d} \nonumber
\\
&\leq& c_6 n^{-d}. \nonumber
\end{eqnarray}
From now on, we will assume that $n \geq c_5 \sqrt{d \log d}$ and $t
\geq2 \beta c_3 d n^2 $.

Let $\widehat N_t$ be the number of jumps made by $Z$ by the time $t$
and note $\widehat N$ is stochastically minorized by $N^*$.
By (\ref{ma21})--(\ref{ma22}),
for $x\in V$,
%
%
\begin{eqnarray}
\label{ma31} \mathbb{P}(Z_t = x)&=& \sum
_{j=0}^\infty\mathbb{P}(Z_t = x \mid
\widehat N_t = j) \mathbb{P}(\widehat N_t = j)\nonumber
\\
&\geq&\sum_{j> c_3 d n^2} \mathbb{P}(Z_t = x
\mid\widehat N_t = j) \mathbb{P}(\widehat N_t = j)
\nonumber
\\
&\geq&\sum_{j> c_3 d n^2} \bigl(n^{-d}/2\bigr)
\mathbb{P}(\widehat N_t = j)
\nonumber\\[-8pt]\\[-8pt]
&\geq&\bigl(n^{-d}/2\bigr) \mathbb{P}\bigl(\widehat N_t >
c_3 d n^2\bigr)
\nonumber
\\
&\geq&\bigl(n^{-d}/2\bigr) \mathbb{P}\bigl( N^*_t >
c_3 d n^2\bigr)
\nonumber
\\
&\geq& c_7 n^{-d}.
\nonumber
\end{eqnarray}

It follows from (\ref{ma32})
that $\mathbb{P}(Z_t-\widetilde Z_t = \bzero) = \mathbb{P}(Z^*_t =
\bzero)
\leq c_6 n^{-d}$, so for fixed $t$ and $n$, there must exist $V_1
\subset V$ with $|V_1| \geq n^d/2$, such that for all $x\in V_1$,
%
%
\begin{equation}
\label{ma33} \mathbb{P}(Z_t=\widetilde Z_t = x) \leq2
c_6 n^{-2d}.
\end{equation}

Let $\mathcal{G}_t = \sigma(\mathcal{M}_s, 0\leq s \leq t)$.
It follows easily from the definition of $Z$ that for $x\in V$,
\[
\mathbb{P}(Z_t =x \mid\mathcal{G}_t) =
M^x_t/n^d
\]
and, by (\ref{ma31}),
%
%
\begin{eqnarray}
\label{ma34} \mathbb{E}M^x_t &=& n^d
\mathbb{E}\mathbb{P}(Z_t =x \mid\mathcal{G}_t)
=n^d \mathbb{P}(Z_t =x ) \geq c_7.
\end{eqnarray}

The random variables $Z_t$ and $\widetilde Z_t$
are conditionally independent given $\mathcal{G}_t$, so for $x\in V$,
\[
\mathbb{P}(Z_t =\widetilde Z_t =x \mid
\mathcal{G}_t) = \bigl(M^x_t
/n^d\bigr)^2.
\]
Thus, by (\ref{ma33}), for $x\in V_1$,
\begin{eqnarray*}
\mathbb{E}\bigl(M^x_t\bigr)^2 &=&
n^{2d} \mathbb{E}\mathbb{P}(Z_t =\widetilde
Z_t =x \mid\mathcal{G}_t) =n^{2d}
\mathbb{P}(Z_t =\widetilde Z_t =x ) \leq2
c_6.
\end{eqnarray*}
Let $c_8 = \sqrt{2c_6}$. We have for $j\geq1$, $x \in V_1$,
\[
\mathbb{P}\bigl( 2^j c_8 \leq M^x_t
\leq2^{j+1} c_8\bigr) \leq2^{-2j}
c_8^{-2} \mathbb{E}\bigl(M^x_t
\bigr)^2 \leq2^{-2j}.
\]
Let $j_1$ be such that
$\sum_{j\geq j_1} 2^{-j+1} c_8
\leq c_7/2$.
Then, by (\ref{ma34}) and the last estimate,
\begin{eqnarray*}
c_7 &\leq&\mathbb{E}M^x_t
\\
&\leq&(c_7/4) \mathbb{P}\bigl(0<M^x_t
\leq c_7/4\bigr) + 2^{j_1+1} c_8 \mathbb{P}
\bigl(c_7/4\leq M^x_t \leq2^{j_1+1}
c_8\bigr)
\\
&&{}+ \sum_{j\geq j_1} 2^{j+1}
c_8 \mathbb{P}\bigl( 2^j c_8 \leq
M^x_t \leq2^{j+1} c_8\bigr)
\\
&\leq& c_7/4 +2^{j_1+1} c_8 \mathbb{P}
\bigl(c_7/4\leq M^x_t \leq2^{j_1+1}
c_8\bigr) + \sum_{j\geq j_1}
2^{j+1} c_8 2^{-2j}
\\
&\leq& c_7/4 + 2^{j_1+1} c_8 \mathbb{P}
\bigl(c_7/4\leq M^x_t \leq2^{j_1+1}
c_8\bigr) + c_7/2,
\end{eqnarray*}
and, therefore, for $x\in V_1$,
\[
\mathbb{P}\bigl( c_7/4\leq M^x_t
\leq2^{j_1+1} c_8\bigr) \geq c_7
c_8^{-1} 2^{-j_1-3}.
\]
Let $c_9 = c_7 c_8^{-1} 2^{-j_1-3}$. Assume that $\mathcal{M}$ and
$\widetilde{\mathcal{M}}$ are independent. Then, for $x\in V_1$,
%
%
\begin{equation}
\label{ma35} \mathbb{P}\bigl( M^x_t \geq
c_7/4, \widetilde M^x_t \geq
c_7/4\bigr) \geq c_9^2.
\end{equation}
Let $K$ be the number of $x$ such that
$ M^x_t \geq c_7/4$ and $ \widetilde M^x_t \geq c_7/4$.
Then
\[
\sum_{x\in V} \bigl|M^x_t -
\widetilde M^x_t\bigr| \leq\sum_{x\in V}
M^x_t + \sum_{x\in V}
\widetilde M^x_t - K c_7/4 = 2
n^d - K c_7/4.
\]
Recall that $|V_1| \geq n^d/2$. By (\ref{ma35}),
\[
\mathbb{E} \biggl(\sum_{x\in V} \bigl|M^x_t
- \widetilde M^x_t\bigr| \biggr) \leq2 n^d -
\bigl(n^d/2\bigr) c_9^2 c_7/4.
\]
This completes the proof of (\ref{ma62}).

\textit{Step}~2.
In this step, we will show that (\ref{ma62}) holds (with a different
constant) even if $\mathcal{M}_t$
and $\widetilde{\mathcal{M}}_t$ are not independent. More precisely, we
will argue that
there exist constants $c_1, c_2\in(0,\infty)$ and $c_{10} < 2$, not
depending on $G$, such that
if $n \geq1 \lor c_1 \sqrt{d \log d}$ and $t \geq c_2 d n^2 $, then
for some coupling of $\mathcal{M}_t$
and $\widetilde{\mathcal{M}}_t$,
%
%
\begin{equation}
\label{ma73} \mathbb{E} \biggl(\sum_{x\in V}
\bigl|M^x_t - \widetilde M^x_t\bigr|
\biggr) \leq c_{10} |V|.
\end{equation}

We\vspace*{2pt} will employ several families of WIMPs.
Let $\{Z^j_t, t\geq0\}_{j\geq1}$
and $\{\widetilde Z^j_t, t\geq0\}_{j\geq1}$ be as in Definition~\ref
{wimps}. In particular,
the jump times of all $Z^j$'s and $\widetilde Z^j$'s are defined by the
same family of Poisson processes $\{N^v\}_{v\in V}$. Let $\mathcal
{M}_t$ and $\widetilde{\mathcal{M}}_t$ denote the mass processes
corresponding to $\{N^v\}_{v\in V}$.

Let $\{X^j_t, t\geq0\}_{j\geq1}$
be jointly distributed as $\{Z^j_t, t\geq0\}_{j\geq1}$.
Similarly, let
$\{\widetilde X^j_t, t\geq0\}_{j\geq1}$
be jointly distributed as $\{\widetilde Z^j_t, t\geq0\}_{j\geq1}$.
However, we make the family $\{X^j_t, t\geq0\}_{j\geq1}$ independent\vspace*{1pt} of
$\{\widetilde X^j_t, t\geq0\}_{j\geq1}$.
Let $\{\mathcal{R}_t, t\geq0\}$ have the same distribution as
$\{\mathcal{M}_t, t\geq0\}$, and assume that $\mathcal{R}_t$ is
driven by the same family
of Poisson processes as $\{X^j_t, t\geq0\}_{j\geq1}$.
By analogy,
let $\{\widetilde{\mathcal{R}}_t, t\geq0\}$ have the same distribution as
$\{\widetilde{\mathcal{M}}_t, t\geq0\}$, and assume that $\widetilde
{\mathcal{R}}_t$ is driven by the same family
of Poisson\vspace*{2pt} processes as $\{\widetilde X^j_t, t\geq0\}_{j\geq1}$.
The processes $\mathcal{R}_t= \{R^x_t\}_{x\in V}$ and $\widetilde
{\mathcal{R}}_t= \{\widetilde R^x_t\}_{x\in V}$ are independent.

Fix some $t>0$ and integer $m>0$.
We find a maximal matching between (some) $X^j$'s and (some)
$\widetilde X^j$'s; that is, we find an asymmetric\vspace*{2pt} relation $\sim$ (a
subset of $\{1, 2,\ldots, m\}^2$) such that $i \sim j$ only if $X^i_t
= \widetilde X^j_t$.
Moreover $i\sim j_1$ and $i\sim j_2$ implies $j_1= j_2$ and, similarly,
$i_1 \sim j$ and $i_2 \sim j$ implies $i_1=i_2$.
Among all such relations $\sim$ we choose one of those that have the
greatest number of matched pairs. Note that for every $x\in V$, either
all $i$ with $X^i_t=x$ are matched with some $j$, or
all $j$ with $\widetilde X^j_t=x$ are matched with some $i$ (or both).
Recall that $\sim$ depends on $m$ and let $r_m$ be the (random) number
of matched pairs.

By the law of large numbers, a.s., for $x\in V$,
\[
\lim_{m\to\infty} \frac{1}m \sum
_{j=1}^m \bone_{\{X^j_t = x\}} =
R^x_t/|V|, \qquad\lim_{m\to\infty}
\frac{1}m \sum_{j=1}^m
\bone_{\{\widetilde X^j_t = x\}} = \widetilde R^x_t/|V|.
\]
This implies that
\[
\lim_{m\to\infty} \Biggl( \frac{1} m \sum
_{j=1}^m \bone_{\{X^j_t = x\}} - \frac{1}m
\sum_{j=1}^m \bone_{\{\widetilde X^j_t = x\}} \Biggr)
= \frac{1} {|V|}\bigl( R^x_t - \widetilde
R^x_t\bigr),
\]
and, a.s.,
%
%
\begin{equation}
\label{ma82} \lim_{m\to\infty} \biggl( \frac{1} m (2m -
r_m) \biggr) = \frac{1} {|V|} \sum_{x\in V}
\bigl|R^x_t - \widetilde R^x_t\bigr|.
\end{equation}
Hence, a.s.,
%
%
\begin{equation}
\label{ma83} \frac{1} {|V|}\sum_{x\in V}
\bigl|R^x_t - \widetilde R^x_t\bigr| = 2-
\lim_{m\to\infty} \frac{r_m}m.
\end{equation}

Next we define a new relation $\approx$ (a subset of $\{1, 2,\ldots,
m\}^2$). Recall that $t>0$ and $m>0$ are fixed. We will construct
$\approx$ by adding pairs to this relation in a dynamic way. We start
by letting $i \approx j$ if $i \sim j$ at time 0. Informally speaking,
we match $X^i_0$ and $\widetilde X^j_0$ if they are at the same vertex,
and we try to match as many pairs as possible at the initial time. We
wait until the first time $s_1>0$ when there exist $i_1 $ and $j_1$
such that $i_1 \not\approx j$ for all $j$, $i\not\approx j_1$ for all
$i$, and $X^{i_1}_{s_1} = \widetilde X^{j_1}_{s_1}$. We add the pair
$(i_1,j_1)$ to the relation $\approx$. We proceed by induction. Given
$s_{k-1}$, let $s_k>s_{k-1}$ be the first time when there exist $i_k $
and $j_k$ such\vspace*{1pt} that $i_k \not\approx j$ for all~$j$, $i\not\approx
j_k$ for all~$i$ (at times between $s_{k-1}$ and $s_k$), and
$X^{i_k}_{s_k} = \widetilde X^{j_k}_{s_k}$. We add the pair $(i_k,j_k)$
to the relation~$\approx$. We proceed in this way until time $t$. Let
$r^*_m$ be the number of matched pairs at time $t$.

We will find a lower bound for $r^*_m$ in terms of $r_m$.
Suppose that $i_1 \sim j_1$.
This implies that $X^{i_1}_t = \widetilde X^{j_1}_t$.
Hence it is possible that $i_1 \approx j_1$. In this case, a pair
$(i_1, j_1)$ that is in relation $\sim$ is also in relation $\approx$.

If $i_1 \not\approx j_1$, then it must be the case that
in the construction of the relation $\approx$, either $X^{i_1}$ was
matched with some $\widetilde X^{j^-_1}$ before time $t$, or
$\widetilde X^{j_1}$ was matched with some $ X^{i^+_1}$ before time
$t$, or both.
We will write $i\,\dot\approx\, j$ if and only if $j \approx i$.
Let
%
%
\begin{eqnarray}
\label{ma81}
\qquad i^-_{\min} &\sim& j^-_{\min} \cdots
i^-_2 \sim j^-_2 \,\dot\approx\ i^-_1 \sim
j^-_1\,\dot\approx\ i_1 \sim j_1
\,\dot\approx\
i^+_1 \sim j^+_1 \,\dot\approx\ i^+_2\sim
\nonumber\\[-8pt]\\[-8pt]
&\sim& j^+_2 \cdots i^+_{\max} \sim j^+_{\max}\nonumber
\end{eqnarray}
be the maximal chain with the alternating structure that should be
clear from the formula. The chain does not have to end with $j^+_{\max
}$; it could end with
$i^+_{\max}$. A similar remark applies to the left end of the chain.
The minimal ratio of the number of pairs of integers in the chain which
are in relation $\approx$ to the number of pairs of integers in the
chain which are in relation $\sim$ is $1/2$.

Any two chains of the form given in (\ref{ma81}) are either identical
or disjoint.
Recall that if $i_1 \sim j_1$, then either $i_1 \approx j_1$
or the pair $(i_1,j_1)$ is an element of a chain as in~(\ref{ma81}).
It follows that
%
%
\begin{equation}
\label{ma85} r^*_m \geq r_m/2.
\end{equation}

Recall that $t>0$ is fixed.
If $i\approx j$, then
let $\sigma^X(i,j)=\inf\{t\geq0\dvtx  X^i_t=\widetilde X^j_t\}$.
Otherwise, let $\sigma^X(i,j)=t$.

Recall WIMPs $\{Z^j_t, t\geq0\}_{j\geq1}$
and $\{\widetilde Z^j_t, t\geq0\}_{j\geq1}$.
Let a relation $\triangleq$ be defined relative to these WIMPs in
exactly the same manner as $\approx$ was defined for $\{X^j_t, t\geq
0\}_{j\geq1}$
and $\{\widetilde X^j_t, t\geq0\}_{j\geq1}$.
In\vspace*{2pt} other words, $\triangleq$ matches colliding particles of type $Z^i$
with $\widetilde Z^j$ as soon as the collisions occur, with the
restriction that each particle is matched with at most one other particle.
If $i\triangleq j$, then
let $\sigma^Z(i,j)=\inf\{t\geq0\dvtx  Z^i_t=\widetilde Z^j_t\}$.
Otherwise, let $\sigma^Z(i,j)=t$.

Recall that $(i,j)\in\,\approx$ is equivalent to $i\approx j$.
It is easy to see that the distribution of the family
\[
\bigl(\bigl\{\bigl(X^i_s, \widetilde
X^j_s\bigr), 0\leq s \leq\sigma^X(i,j)
\bigr\} _{(i,j) \,\in\,\approx}, \bigl\{\bigl(X^i_s,
\widetilde X^j_s\bigr), 0\leq s \leq t\bigr\}
_{(i,j)\,\notin\,\approx} \bigr)
\]
is the same as that of
\[
\bigl(\bigl\{\bigl(Z^i_s, \widetilde
Z^j_s\bigr), 0\leq s \leq\sigma^Z(i,j)
\bigr\} _{(i,j) \,\in\,\triangleq}, \bigl\{\bigl(Z^i_s,
\widetilde Z^j_s\bigr), 0\leq s \leq t\bigr\}
_{(i,j)\,\notin\,\triangleq
} \bigr)
\]
because the jump times of the processes in each family are
determined by independent Poisson processes at vertices of the graph.
Let $r^Z_m$ be the number of pairs in the relation $\triangleq$. We
see that the distributions of $r^Z_m$ and $r^*_m$ are identical.

If $\sigma^Z(i,j) < t$, then we let $\widehat Z_s = \widetilde Z_s$
for $s\in[0, \sigma^Z(i,j))$ and $\widehat Z_s = Z_s$ for $s\geq
\sigma^Z(i,j)$. Note that the distribution of the family
$\{\widehat Z^j_t, t\geq0\}_{j\geq1}$ is the same as that of the
family $\{\widetilde Z^j_t, t\geq0\}_{j\geq1}$.
If $i \triangleq j$, then $Z^i_t=\widehat Z^j_t$.
We have
\[
\lim_{m\to\infty} \biggl( \frac{1} m \bigl(2m -
r^Z_m\bigr) \biggr) = \frac{1} {|V|} \sum
_{x\in V} \bigl|M^x_t - \widetilde
M^x_t\bigr|,
\]
for the same reason that (\ref{ma82}) holds.
Therefore, using (\ref{ma62}), (\ref{ma83}), (\ref{ma85}) and the
equality of the distributions of $r^Z_m$ and $r^*_m$, we obtain for $n
\geq1 \lor c_1 \sqrt{d \log d}$ and $t \geq c_2 d n^2 $,
%
%
\begin{eqnarray}
\label{ma84} \qquad&&\mathbb{E} \biggl(\frac{1} {|V|}\sum
_{x\in V} \bigl|M^x_t - \widetilde
M^x_t\bigr| \biggr)
\\
&&\qquad = 2 - \mathbb{E} \biggl(\lim_{m\to\infty} \frac{r^Z_m}m \biggr)
\\
&&\qquad \leq2 - \mathbb{E} \biggl(\lim_{m\to\infty} \frac{r_m}{2m}
\biggr)
\nonumber
\\
&&\qquad  = 2 -\mathbb{E} \biggl(\lim_{m\to\infty} \frac{r_m}{m}
\biggr)+ \frac{1}2 \mathbb{E} \biggl(\lim_{m\to\infty}
\frac{r_m}{m} \biggr)
\nonumber
\\
&&\qquad  = \mathbb{E} \biggl(\frac{1} {|V|}\sum_{x\in V}
\bigl|R^x_t - \widetilde R^x_t\bigr|
\biggr) + \frac{1}2 \biggl( 2 - \mathbb{E} \biggl(\frac{1} {|V|}\sum
_{x\in V} \bigl|R^x_t - \widetilde
R^x_t\bigr| \biggr) \biggr)
\nonumber
\\
&&\qquad =1+ \frac{1}2 \mathbb{E} \biggl(\frac{1} {|V|}\sum
_{x\in V} \bigl|R^x_t - \widetilde
R^x_t\bigr| \biggr) \leq1 + c_4/2.
\nonumber
\end{eqnarray}
This proves (\ref{ma73}).

\textit{Step}~3.\vspace*{2pt}
The process $\mathcal{M}_t$ is ``additive'' in the following sense.
Suppose that $\mathcal{M}_t$ and $\widetilde{\mathcal{M}}_t$ are
driven by the same family of Poisson processes $N^v$. Let \mbox{$\widehat
{\mathcal{M}}_0 = \mathcal{M}_0 + \widetilde{\mathcal{M}}_0$,} and\vspace*{1pt}
suppose that $\widehat{\mathcal{M}}_t$ is\vspace*{1pt} also driven by the same
family of Poisson processes $N^v$. Then
$\widehat{\mathcal{M}}_t = \mathcal{M}_t + \widetilde{\mathcal{M}}_t$
for all $t$, a.s.

Fix $t = c_2 d n^2 $ and suppose that $\mathcal{M}_t$ and $\widetilde
{\mathcal{M}}_t$ are driven by the same family of Poisson processes $N^v$.
Let
\begin{eqnarray*}
\mathcal{M}^+_t &=& (\mathcal{M}_t - \widetilde{
\mathcal{M}}_t) \lor0,
\\
\widetilde{\mathcal{M}}^+_t &=& (\widetilde{\mathcal{M}}_t
- \mathcal{M}_t) \lor0,
\\
\mathcal{M}_t^c &=& \bigl(\mathcal{M}_t -
\mathcal{M}^+_t\bigr) = \bigl(\widetilde{\mathcal{M}}_t
-\widetilde{\mathcal{M}}^+_t\bigr).
\end{eqnarray*}
The\vspace*{2pt} process $\mathcal{M}^c_t$ represents the maximum matching mass at
every site, and processes
$\mathcal{M}^+_t$ and $\widetilde{\mathcal{M}}^+_t$ represent the
excesses of $\mathcal{M}_t$ and $\widetilde{\mathcal{M}}_t$ (if any)
above the common mass. Suppose that all these processes are driven by
the same family of Poisson processes $N^v$ after time $t$.
Then for every $s\geq t$,
\[
\mathcal{M}_s = \mathcal{M}^c_s +
\mathcal{M}^+_s, \qquad\widetilde{\mathcal{M}}_s =
\mathcal{M}^c_s + \widetilde{\mathcal{M}}^+_s.
\]
By\vspace*{1pt} the Markov property applied at time $t = c_2 d n^2 $ and (\ref
{ma73}) applied to $\mathcal{M}^+_s$~and~$\widetilde{\mathcal
{M}}^+_s$, we obtain for $s \geq2 c_2 d n^2 $,
\[
\mathbb{E} \biggl(\sum_{x\in V} \bigl|\bigl(M^+
\bigr)^x_s - \bigl(\widetilde M^+\bigr)^x_s\bigr|
\mid\mathcal{F}_t \biggr) \leq c_{10} \frac{1}2
\sum_{x\in V} \bigl(\bigl(M^+\bigr)^x_t
+ \bigl(\widetilde M^+\bigr)^x_t\bigr).
\]
Hence
\begin{eqnarray*}
\mathbb{E} \biggl(\sum_{x\in V} \bigl|\bigl(M^+
\bigr)^x_s - \bigl(\widetilde M^+\bigr)^x_s\bigr|
\biggr) &\leq& (c_{10} /2) \mathbb{E} \biggl(\sum
_{x\in V} \bigl(\bigl(M^+\bigr)^x_t +
\bigl(\widetilde M^+\bigr)^x_t\bigr) \biggr)
\\
&=&(c_{10} /2)\mathbb{E} \biggl(\sum_{x\in V}
\bigl|M^x_t - \widetilde M^x_t\bigr|
\biggr)
\\
&\leq& (c_{10} /2) c_{10} |V|.
\end{eqnarray*}
An inductive argument applied at times $t$ of the form $t = jc_2 d n^2
$, $j \geq2$, yields for $s \geq(j+1) c_2 d n^2 $,
\[
\mathbb{E} \biggl(\sum_{x\in V} \bigl|\bigl(M^+
\bigr)^x_s - \bigl(\widetilde M^+\bigr)^x_s\bigr|
\biggr) \leq(c_{10} /2)^j c_{10} |V|.
\]
This implies (\ref{f91}) and completes the proof.
\end{pf}

%
\begin{remark} (i)
Let $\| \cdot\|_{\TV}$ denote the total variation distance,
and let the mixing time for the random walk on $G$ be defined by
\[
\mathcal{T}= \inf\Bigl\{t\geq0\dvtx  \sup_\mu\|\mu
P_t - \pi\|_{\TV} \leq1/4\Bigr\},
\]
where $\pi$ stands for the stationary distribution, and $P_t$
denotes the transition kernel.
See \cite{LPW}, Chapter~4, for these definitions and
various results on mixing times.

Consider the graph $\calC_n^d$ for some $n,d\geq3$. For this graph,
$\alpha(t) $ defined in (\ref{f122}) is equal to $1/2$ for $t$ of
order $n^d$.
Theorem~\ref{ma61} shows that the left-hand side of~(\ref{f91}) is
bounded by $n^d/2$ for $t$ of order $n^2$, thus greatly improving (\ref
{f91.3.2}) in the case $G=\calC_n^d$.
Since the mixing time for random walk on $\calC_n^d$ is of the order
$n^2$, the bound in Theorem~\ref{ma61} cannot be improved in a
substantial way.
Recall $\bP$ defined before Theorem~\ref{f101}.

\begin{con*}
The mixing time for the random walk
corresponding to $\bP$ should give the optimal bound in (\ref{f91.3.2}).
\end{con*}

A support to our conjecture is lent by the
recent proof (see \cite{CLR}) of the ``Aldous spectral gap
conjecture,'' saying that
the ``interchange process'' and the corresponding random
walk have the same spectral gap.

(ii)
The proof of
Theorem~\ref{ma61} depends on the assumption that $G=\calC_n^d$
only at one point, namely, the estimate
%
%
\begin{eqnarray}
\label{ma71} && \mathbb{P}(Z_t=\widetilde Z_t = x \mid
Z_t = y, \widetilde Z_0 = z)
\nonumber\\[-8pt]\\[-8pt]
&&\qquad \leq c_0
\mathbb{P}(Z_t= x \mid Z_t = y) \mathbb{P}(\widetilde
Z_t = x \mid\widetilde Z_0 = z)\nonumber
\end{eqnarray}
is derived using properties of $\calC_n^d$ in an essential way.
In other words, if a similar estimate can be obtained
for some other family of graphs, the proof of the theorem would
apply in that case. It is not hard to construct
examples showing that there is no universal constant $c_0$
such that (\ref{ma71}) holds for all finite graphs $G$, all \mbox{$x,y,z\in V$}
and all $t>0$. Hence, any generalization of Theorem~\ref{ma61}
has to be limited to a subfamily of finite graphs or come
with a different proof.
\end{remark}

We now present very elementary properties of the meteor process.

%
\begin{proposition}\label{a298}
Let $T^v_t$ denote the time of the last jump of $N^v$ on the interval
$[0,t]$, with the convention that $T^v_t=-1$ if there were no jumps on
this interval. Let $U(v)=\{v\} \cup\{x\in V\dvtx  x\leftrightarrow v\}$.
\begin{longlist}[(iii)]
\item[(i)] Assume that $ M^v_0 + M^x_0 >0$ for a pair of adjacent vertices $v$
and $x$. Then, almost surely, for all $t\geq0$, $ M^v_t + M^x_t >0$.

\item[(ii)]
Let $R_t$ be the number of pairs $(x,v)$ such that $x\leftrightarrow
v$ and
$ M^v_t + M^x_t =0$. The process $R_t$ is nonincreasing, a.s.

\item[(iii)] Assume that $ M^x_0 >0$ for $x\in U(v)\setminus\{v\}$.
Then, a.s., $M^v_t = 0$ if and only if one of the following conditions
holds: \textup{(a)} $T^v_t = \max\{T^x_t\dvtx  x\in U(v)\}>-1$ or \textup{(b)}~$M^v_0 = 0$
and $ \max\{T^x_t\dvtx  x\in U(v)\setminus\{v\}\}=-1$.

\item[(iv)] Suppose that the process $\{\mathcal{M}_t, t\geq0\}$ is in the
stationary regime, that is, its distribution at time 0 is the
stationary distribution $Q$. Then $ M^v_t + M^{x}_t >0$ for all $t\geq
0$ and all pairs of adjacent vertices $v$ and $x$, a.s.

\item[(v)]
Recall from Remark~\ref{j132} the stationary mass process $\{\mathcal
{M}_t, t\in\mathbb{R}\}$ and the corresponding Poisson processes $\{
N^v_t,t\in\mathbb{R}\} $, $v\in V$.
Let $T^v$ denote the time of the last jump of $N^v$ on the interval
$(-\infty, 0]$, and note that $T^v$ is well defined for every $v$
because such a jump exists, a.s.
Then, a.s., $M^v_0 =0$ if an only if $T^v = \max\{T^x\dvtx  x\in U(v)\}$.
\end{longlist}
\end{proposition}

\begin{pf}
(i) Suppose to the contrary that $ M^v_t + M^x_t =0$ for some
$x\leftrightarrow v$ and $t>0$.
The value of $ M^v_t + M^x_t$ can change only when one of the processes
$N^y$, $y \in U_1:= U(x) \cup U(v)$, has a jump.
Note that $U_1$ is a finite set. It follows that the union of jump
times of all processes $N^y$, $y\in U_1$, does not have accumulation
points. Moreover, jumps of different processes $N^y$ in this family
never occur at the same time, a.s. Let $T$ be the infimum of times such
that $ M^v_t + M^x_t =0$. Then $ M^v_s + M^x_s >0$ for all $s< T$ and $
M^v_{T-} + M^x_{T-} >0$. Suppose without loss of generality that
$M^v_{T-} >0$. If $N^v$ has a jump at $T$, then $M^x_T>0$, a
contradiction. If $N^y$ has a jump at $T$ for some $y\in U_1$, then
$N^v$ does not have a jump at $T$ and, therefore, $M^v_T>0$, also a
contradiction. We conclude that the assumption that $ M^v_t + M^x_t =0$
for some $x\leftrightarrow v$ and $t>0$ is false.

(ii)
For $x\leftrightarrow v$, let $R^{x,v}_t$ be 1 if $ M^v_t + M^x_t =0$
and 0 otherwise. This process is nonincreasing, by part (i). Since
$R_t = \sum_{x,v \in V, x\leftrightarrow v} R^{x,v}_t$, it follows
that $R_t$ is nonincreasing.

(iii)
If $M^v_0 = 0$ and $ \max\{T^x_t\dvtx  x\in U(v)\setminus\{v\}\}=-1$, then
processes $N^x$, $x\in U(v)\setminus\{v\}$, do not jump on the
interval $[0,t]$. Hence, $M^v_s = 0 $ for all $s\in[0,t]$. In
particular, $M^v_t = 0$. We will assume that $ \max\{T^x_t\dvtx  x\in
U(v)\setminus\{v\}\}>-1$ in the rest of the proof.

Suppose that $T^v_t = \max\{T^x_t\dvtx  x\in U(v)\}>-1$. Then $M^v_{T^v_t}
= 0$. Since processes $N^x$, $x\leftrightarrow v$, do not have jumps
on the interval $[T^v_t, t]$, we must have $M^v_s = 0$ for all $s\in
[T^v_t, t]$. Hence, $M^v_t = 0$.

Suppose that $T^v_t < \max\{T^x_t\dvtx  x\in U(v)\}$
and let $y$ be such that $T^y_t =  \max\{T^x_t\dvtx  x\in U(v)\}>-1$ and
$y\leftrightarrow v$. By\vspace*{-1pt} part (i), either $M^v_{T^y_t-} >0$ or
$M^y_{T^y_t-} >0$ (or both).
If $M^v_{T^y_t-} >0$, then $M^v_s > 0$ for all $s\in[T^y_t, t]$
because $N^v$ does not jump on this interval.
If $M^y_{T^y_t-} >0$, then $M^v_{T^y_t} >0$ and,\vspace*{-1pt} therefore,
$M^v_s > 0$ for all $s\in[T^y_t, t]$ because $N^v$ does not jump on
this interval.
We see that in either case, $M^v_t > 0$.

(iv)
Since $V$ is finite, there exists a sequence $(x_1, x_2,\ldots, x_n)$
of vertices of $G$ such that $x_j \leftrightarrow x_{j+1}$ for $1\leq
j \leq n-1$, $x_n \leftrightarrow x_1$, and
the sequence contains all vertices in $V$.
The vertices $x_j$'s are not necessarily distinct.

Let $A_i$
be the event that processes $N^v$, $v\in V$, have $2n$ jumps in the
time interval $[i, i+1)$, and the jumps occur at the following vertices
in the following order: $x_1, x_2,\ldots, x_n, x_1, x_2,\ldots, x_n$.
It is easy to see that if $A_i$ occurs, then there is only one vertex
$v$ with $M^v_{i+1}=0$; specifically,\vspace*{1pt} $v=x_n$. Hence if $A_i$ occurs,
then \mbox{$R_{i+1}=0$} and, by part (ii), $R_t = 0$ for $t\geq i+1$. Events
$A_i$ are independent and have positive probability so the probability
of $A_1 \cup A_2\cup\cdots\cup A_m$ is bounded below by $ 1- p^m$ for
some $p<1$. It follows that $\mathbb{P}(R_{m+1} = 0) \geq1 - p^m$
for $m\geq1$. This and stationarity imply that $\mathbb{P}(R_{0} =
0)=\mathbb{P}(R_{m+1} = 0)=1$ for $m\geq1$.

(v) It follows from part (iv) that $ M^v_0 + M^{x}_0 >0$ for all pairs
of adjacent vertices $v$ and $x$, a.s. Hence $ M^v_k + M^{x}_k >0$ for
all $k\in{\mathbb Z}$ and, therefore, $ M^v_t + M^{x}_t >0$ for all
pairs of adjacent vertices $v$ and $x$ and all $t\in\mathbb{R}$, a.s.
Now we can apply the same reasoning as in the proof of case (a) in part (iii).
\end{pf}

\section{Meteor craters in circular graphs}\label{craters}

This section is devoted to meteor processes on circular graphs.
Recall that $\calC_k$ denotes the circular graph with $k$ vertices,
$k\geq2$.

We will say that there is a \emph{crater} at the site $j$ at time $t$
if $M^j_t =0$. Craters are special features of the meteor process for a
number of reasons. First, the mass at a crater has the minimum possible
value. Second, we expect that the distribution of mass $M^j_t$ is a
mixture of an atom at 0 and a distribution with a continuous density.
Third, given the distribution of mass $\mathcal{M}_s$ at all sites at
time $s$ and positions of all craters at times $t\in[s,u]$, we can
determine the mass process $\{\mathcal{M}_t,
t\in[s,u]\}$.
For these reasons, we find it interesting to study the distribution of
craters. An easy argument (see the proof of Theorem~\ref{a285}) shows
that the concept of a crater is essentially equivalent to a \emph
{peak} in a random (uniform) permutation.

The research on peaks and other related permutation statistics, such as
valleys, descents and runs has a very long history. For a review of
some related literature, see the introduction to \cite{CongVis}; the
authors trace the beginning of this line of research to the nineteenth
century. However, the research in this area seems to have a number of
separate lines, because the authors of \cite{CongVis} do not cite
\cite{Chao} or \cite{KMK2,KMK1}.
In view of this disconnected nature of the literature we are not sure
whether we were able to trace all the existing results in the area that
are relevant to our paper.

There are (at least) three natural probabilistic questions that have to
do with craters. The first one is concerned with the probability of a
given pattern of craters. This is equivalent, more or less, to the
question about the asymptotic frequency of a given pattern of craters
in a very large cyclic graph $\calC_k$. We will provide formulas for
two specific crater ``patterns''
in Theorems~\ref{a285} and~\ref{j191}. It is possible that both
results could be derived from \cite{KMK2,KMK1}, but the style of those
old papers may be hard to follow for the modern reader. We will base
our proofs on the combinatorial results in \cite{BBS}. The results in
\cite{BBS} could be used to derive more advanced theorems on craters
that go beyond the scope of this paper.

The second question is about fluctuations of the number of copies of a
pattern. There are a number of combinatorial versions of the
central limit theorem for permutation statistics; see, for example,
\cite{Chao,CongVis} and references therein. We will state a theorem
that appeared in \cite{Bona}, and we will provide a new short proof
based on classical probabilistic tools
and our meteor process.

Finally, there is a question of large deviations for the crater
process. We will not provide a standard large deviations result, but we
will prove a theorem on the most likely configuration of craters
assuming that there are very few of them.

Consider the meteor process on $\calC_k$, and
assume that $\sum_{j=1}^k M^j_0 = k$.
For $n\geq1$ and $k \geq n+4$, let
$F^k_n$ be the event that $M^i_0>0$ for $i=3,4,\ldots, n+2$ and
$M^{2}_0 = M^{n+3}_0 =0$.
In other words, $F^k_n$ is the event that $3$ is the starting point of
a maximal sequence of vertices which are not craters at time $0$ and
that sequence has length~$n$.
We let $F^k_0$ be the event that $M^2_0=0$.

For $n\geq1$ and $k \geq n+3$,
let $\widehat F^k_n$ be the event that $M^i_0>0$ for $i=3,\ldots, n+2$.
In other words, $\widehat F^k_n$ is the event that
sites $3,\ldots, n+2$ are not craters at time 0, but this sequence
does not have to be maximal.

%
\begin{theorem}\label{a285}
Consider the meteor process on $\calC_k$ in the stationary regime;
that is, assume that the distribution of $\mathcal{M}_0$ is the
stationary measure $Q$.
We have
%
%
\begin{eqnarray}
\label{j131} p_0 &:=& \mathbb{P}\bigl(F^k_0
\bigr) = 1/3, \qquad k\geq3,
\\
p_n &:=& \mathbb{P}\bigl(F^k_n\bigr) =
\frac{n (n+3) 2^{n+1}}{(n+4)!}, \qquad n\geq1, k \geq n+4, \label{a296}
\\
\widehat p_n &:=& \mathbb{P}\bigl(\widehat F^k_n
\bigr) = \frac
{2^{n+1}}{(n+2)!}, \qquad n\geq1, k \geq n+3. \label{m31}
\end{eqnarray}
\end{theorem}

\begin{pf}
Recall from Remark~\ref{j132} the stationary mass process $\{\mathcal
{M}_t, t\in\mathbb{R}\}$ and the corresponding Poisson processes $\{
N^j_t,t\in\mathbb{R}\} $, $j=1,\ldots,k$, defined on the whole real
time-line.
As in Proposition~\ref{a298}(v), we
let $T^j$ denote the time of the last jump of $N^j$ on the interval
$(-\infty, 0]$.
According to Proposition~\ref{a298}(v), $M^j_0 =0$ if and only if
$T^{j-1}<T^j>T^{j+1}$.

Note that $T^m \ne T^j$ if $m\ne j$, a.s.
Let $a_1 \cdots a_k$ be the random permutation of $\{1,2,\ldots, k\}$
defined by the condition $a_j < a_m$ if and only if $T^j < T^m$, for
all $j$ and $m$. It is clear that
$a_1 \cdots a_k$ is the uniform
random permutation of $\{1,2,\ldots, k\}$.

We say that $j$ is a \emph{peak} (of the permutation $a_1 \cdots a_k$)
if $ a_{j-1} < a_j > a_{j+1}$.
Hence $M^j_0 = 0$ if and only if $a_j$ is a peak.

By symmetry, any of the random numbers $a_{j-1}, a_j$ and $a_{j+1}$ is
the largest of the three
with the same probability. Hence the probability that $ a_{j-1} <a_j >
a_{j+1}$ is $1/3$.
This proves (\ref{j131}).

The event $F^k_n$ holds if and only if
in the initial part $a_1 \cdots
a_{n+4}$ of the permutation $a_1 \cdots a_k$, there are exactly two
peaks at $2$ and $n+3$.
It is clear that the probability of this event is the same if $a_1
\cdots a_{n+4}$ is a random uniform permutation of $\{1,\ldots, n+4\}$
with the same peak set.
Recall that the number of permutations of $\{1,\ldots, n+4\}$ is $(n+4)!$.
We now see that (\ref{a296}) follows from Theorems 1~and~10 in \cite{BBS}.
Note that we are concerned with permutations of size $n+4$ while the
two cited theorems in \cite{BBS} count permutations of size $n$. This
explains the shift of size 4
in the corresponding formulas in our paper and \cite{BBS}.

Finally, we will prove (\ref{m31}).
The event $\widehat F^k_n$ holds if an only if
there are no peaks in the part $a_2 \cdots a_{n+3}$ of the permutation
$a_1 \cdots a_k$. The probability of this event is the same if $a_2
\cdots a_{n+3}$ is a random uniform permutation of $\{1,\ldots, n+2\}$
with no peaks.
Formula (\ref{m31}) follows from Proposition~2 in \cite{BBS}, with
a shift of size 2 between
the corresponding formulas in our paper and \cite{BBS}.
\end{pf}

The results in \cite{BBS} provide an effective tool for calculating various
distributions related to crater positions. We ask the interested reader
to consult
that paper for the general theory. We will provide here another explicit
probabilistic formula based on combinatorial results from \cite{BBS}.

%
\begin{theorem}\label{j191}
Consider the meteor process on $\calC_k$ in the stationary regime;
that is, assume that the distribution of $\mathcal{M}_0$ is the
stationary measure $Q$.\vspace*{2pt}
For \mbox{$i,j\geq1$} and $k \geq i+j+5$, let
$A^k_{i,j}$ be the event that $M^{2}_0 = M^{i+3}_0 = M^{i+j+4}_0=0$
and
$M^n_0>0$ for $n=3,4,\ldots, i+2, i+4,\ldots, i+j+3$.
In other words, $A^k_{i,j}$ is the event that $2$ is a crater and the
gaps between
this crater and the next two craters have sizes $i$ and $j$.
We have
%
%
\begin{eqnarray}\label{j152}
\mathbb{P}\bigl(A^k_{i,j}\bigr) &=&
\frac{2^{i+j}}{(i+j+5)!}
\biggl[(i+j+4) \biggl(j \pmatrix{i+j+1 \cr i-1}+(j+1) \pmatrix{i+j+1 \cr i}\nonumber
\\
&&\hspace*{115pt}{}+(i+1) \pmatrix{i+j+1 \cr i+1}+i \pmatrix{i+j+1 \cr i+2}
\nonumber\\[-8pt]\\[-8pt]
&&\hspace*{226pt}{} -2 (i+j+1) \biggr)\nonumber
\\
&&\hspace*{220pt} {}+i j \pmatrix{i+j+4 \cr i+2} \biggr].\nonumber
\end{eqnarray}
\end{theorem}

\begin{pf}
The theorem follows from Theorems 9~and~12 of \cite{BBS}.
The argument is totally analogous to that in the proof of Theorem
\ref{a285} so we leave the details to the reader. We just note
that one should take $m=i+3$ and $n=i+j+5$ in Theorem~12 of \cite{BBS}.
\end{pf}

%
\begin{remark}
(i)
If craters occurred in the i.i.d. manner, then
the distribution of the distance between consecutive craters would have
been geometric, with the tail decaying exponentially.
This is not the case.
By the Stirling approximation,
\begin{eqnarray*}
p_n &=& \frac{n (n+3) 2^{n+1}}{(n+4)!} \sim\frac{n (n+3) 2^{n+1}
e^{n+4}}{(n+4)^{n+4}\sqrt{2 \pi(n+4)} }.
\end{eqnarray*}
Hence, $p_n$ converges to 0 at a rate faster than exponential;
specifically, $\log p_n \approx-n\log n$.

(ii) Despite remark (i), the crater process
is ``partly'' memoryless.
Consider the crater distribution at time 0 assuming that the mass
process $\{\mathcal{M}_t, t\in\mathbb{R}\}$ is in the~stationary
regime. The event that there is a crater at site $j$ depends only on
the Poisson processes $N^n$ for $n=j-1, j,j+1$, by Proposition~\ref
{a298}(v). Hence,~the events $\{M^j _0 = 0\}$ for $j =1+ 3m$, $m \in
{\mathbb Z}$, $1\leq j \leq k-2$, form a sequence of Bernoulli trials
(are i.i.d.). It follows that the gap between the first and second
craters in this sequence has an approximately geometric tail, for large
$k$. The same observation holds for two similar sequences of sites,
namely for those indexed by $j =2+ 3m$, $m \in{\mathbb Z}$, $1\leq j
\leq k-2$, and those
indexed by $j = 3m$, $m \in{\mathbb Z}$, $1\leq j \leq k-2$.
However, the three sequences of Bernoulli trials are highly dependent.

(iii)
It is natural to ask for the distribution of the number of
consecutive sites with nonzero mass following a crater.
This somewhat informal statement can be translated into a rigorous
question about the conditional probability of $F^k_n$ given
that $M^2_0 = 0$. The answer is $p_n/p_0 = 3n (n+3) 2^{n+1}/(n+4)! $.
In other words, a crater is followed by exactly $n$ consecutive sites
with nonzero mass
with probability $ 3n (n+3) 2^{n+1}/(n+4)!$.

(iv) Remarks (i) and (ii) make it clear that the we should not expect
independence between the lengths of consecutive stretches of
sites with nonzero mass. More precisely, one can easily check
that, in general, for large $k$,
\[
\frac{1}{p_0} \mathbb{P}\bigl(A^k_{i,j}\bigr) \ne
\frac{1}{p_0} \mathbb{P}\bigl(F^k_i\bigr)
\frac{1}{p_0} \mathbb{P}\bigl(F^k_j\bigr).
\]
Curiously, for $j\geq1$ and $k\geq j+7$,
%
%
\begin{equation}
\label{j151} \frac{1}{p_0} \mathbb{P}\bigl(A^k_{2,j}
\bigr) = \frac{1}{p_0} \mathbb{P}\bigl(F^k_2\bigr)
\frac{1}{p_0} \mathbb{P}\bigl(F^k_j\bigr).
\end{equation}
Hence, if there are exactly two noncraters between two consecutive
craters, then this event
gives no information about the length of the next stretch of
sites with nonzero masses.
Formula (\ref{j151}) follows from
(\ref{a296}) and (\ref{j152}) by direct calculation.
Formula (\ref{j151}) does not seem to hold if 2 is replaced
by any other integer $i\geq1$, $i\ne2$.
We offer an informal explanation of (\ref{j151}).
Suppose that there is a crater at site 5. Then there is no crater at
site 4. The distribution of craters at sites $5,6,\ldots$ is
determined by Poisson processes at sites $4,5,\ldots.$ If we have
extra information that there is a crater at site 2, then this tells us
only that the latest meteor hit among the sites $1,2$ and 3 occurred at
site 2. Since the Poisson processes at sites $1,2$ and~3 are
independent of those at sites $4,5,\ldots,$ the information that 2 is
a crater has no predictive value for craters to the right of 5.

When translated into the language of permutation peaks, the condition
discussed in the last paragraph becomes that there are exactly two
nonpeaks between any two consecutive peaks. Interestingly, exactly the
same condition came up as part of a conjecture in~\cite{BBS} about the
equidistribution of peaks in permutations. This part of the conjecture
was recently proved by Kasraoui in \cite{AK}. Is there some deeper
connection between this result and equation~(\ref{j151})?

(v)
Formula (\ref{m31}) is extremely easy to prove; see the counting
argument in the proof of
Proposition~2 in \cite{BBS}.
We will
derive the harder formula (\ref{a296}) from the easier formula (\ref{m31})
in an informal way. It has been shown in \cite{meteor2}, a follow-up
paper, that
stationary distributions on $C_k$ converge to a stationary distribution
for the meteor process on ${\mathbb Z}$, in an appropriate sense.
It is easy to see that for the meteor process on ${\mathbb Z}$,
\[
\widehat p_n = p_n + 2 p_{n+1} + 3
p_{n+2} + \cdots.
\]
We take the inverse of this linear transformation to see that
\[
p_n = \widehat p_n - 2 \widehat p_{n+1} +
\widehat p_{n+2}.
\]
This and (\ref{m31}) imply that
\[
p_n = \frac{2^{n+1}}{(n+2)!} - 2 \frac{2^{n+2}}{(n+3)!} +
\frac{2^{n+3}}{(n+4)!} = \frac{n (n+3) 2^{n+1}}{(n+4)!}.
\]
\end{remark}

%
\begin{theorem}
We have
%
%
\begin{eqnarray}
\label{a281} \sum_{n=0}^\infty
p_n &=& 2/3,
\\
\sum_{n=0}^\infty n p_n &=& 2/3.
\label{a282}
\end{eqnarray}
\end{theorem}

\begin{pf}
Our argument is based on power series expansions derived by Mathematica
\cite{math}. The following power series
converges for all real $x$,
\[
\sum_{n=1}^\infty\frac{2 n (n+3)}{(n+4)!}
x^{n+4} = \frac{2 x^3}{3}+2 x^2+2 e^x
(x-2)^2-8.
\]
From this, we obtain
\begin{eqnarray*}
\sum_{n=0}^\infty p_n &=& 1/3 +
\sum_{n=1}^\infty p_n =1/3 + \sum
_{n=1}^\infty\frac{ n (n+3)}{(n+4)!}
2^{n+1}
\\
&=& 1/3 + 2^{-4}\sum_{n=1}^\infty
\frac{2 n (n+3)}{(n+4)!} 2^{n+4}
\\
&=&1/3 + 2^{-4} \biggl( \frac{2\cdot2^3}{3}+2 \cdot2^2+2
e^2 (2-2)^2-8 \biggr) = \frac{2}3.
\end{eqnarray*}
A similar calculation yields
\[
\sum_{n=1}^\infty\frac{2 n^2 (n+3)}{(n+4)!}
x^{n+4} = 2 e^x \bigl(x^3-6 x^2+16
x-16 \bigr)-\frac{2}{3} \bigl(x^3+6 x^2-48
\bigr)
\]
and
\begin{eqnarray*}
\sum_{n=0}^\infty n p_n &=& \sum
_{n=1}^\infty\frac{ n^2 (n+3)}{(n+4)!}
2^{n+1} = 2^{-4}\sum_{n=1}^\infty
\frac{2 n^2 (n+3)}{(n+4)!} 2^{n+4}
\\
&=& 2^{-4} \biggl( 2 e^2 \bigl(2^3-6
\cdot2^2+16 \cdot2-16 \bigr)-\frac{2}{3} \bigl(2^3+6
\cdot2^2-48 \bigr) \biggr) = \frac{2}3.
\end{eqnarray*}
This completes the proof.
\end{pf}

%
\begin{remark}
(i) The reader may be puzzled by (\ref{a281}) since the probabilities
do not add up to 1. This sequence of probabilities does not represent
all events in a partition of a probability space. For the meteor
process on ${\mathbb Z}$ constructed in~\cite{meteor2}, the
probabilities $p_n$ represent only the events that a given vertex has
no mass or it is the starting point of a sequence of consecutive
vertices, all with positive masses. It is also possible for a vertex to
be an interior point of
a sequence of consecutive vertices with positive masses. It follows
from (\ref{a281}) that the last event has probability~$1/3$.

(ii)
We will present a simple heuristic proof of (\ref{a281}) and (\ref
{a282}) based on the meteor process on ${\mathbb Z}$ constructed in
\cite{meteor2}.
Recall from (\ref{j131}) that $p_0 = 1/3$. The number of starting
points of sequences of consecutive vertices with positive masses must
be the same as the number of craters, since such vertices are never
adjacent, by Proposition~\ref{a298}(iv). Hence, $\sum_{n=1}^\infty
p_n = 1/3$, implying
(\ref{a281}).
The sum $\sum_{n=1}^\infty n p_n $
represents the proportion of noncraters so it must be equal to $2/3$
because the proportion of craters is $p_0=1/3$.
\end{remark}

One can ask not only how often a given configuration of craters occurs
in a very large circular graph $\calC_k$ but also what the random
fluctuations are. We will prove a central limit theorem to shed some
light on this problem. To match well the existing literature, our
formulation will be more general than necessary for the purpose of
describing the configuration of craters.

Consider the meteor process on $\calC_k$ in the stationary regime;
here and later in this section this means that the distribution of
$\mathcal{M}_0$ is the stationary measure $Q$.
Recall from Remark~\ref{j132} that the stationary mass process
$\{\mathcal{M}_t, t\in\mathbb{R}\}$ and the corresponding Poisson processes
$\{N^m_t,t\in\mathbb{R}\} $, $m=1,\ldots,k$, are defined on the
whole real time-line.
As in Proposition~\ref{a298}(v), we
let $T^m$ denote the time of the last jump of $N^m$ on the interval
$(-\infty, 0]$.
A permutation $\ba= a_1 a_2 \cdots a_n$ of $\{1,\ldots,n\}$ will be
called a \emph{pattern}.
We will denote finite families of patterns by $\calA=\{\ba^1,\ldots,
\ba^m\}$. We will not assume that all patterns in $\calA$ have the
same length.\vspace*{1pt}
We will say that $\calA$ occurs at $j$ if for some $\ba^r = a^r_1
\cdots a^r_{n_r} \in\calA$, we have $T^{j+i-1} < T^{j+m-1}$ if and
only if $a^r_i < a^r_m$ for all $1\leq i,m\leq n_r$.

According to Proposition~\ref{a298}(v), $M^m_0 =0$ if an only if $
T^{m-1} <T^m > T^{m+1}$. Hence, any finite configuration of craters can
be represented as a family of patterns.

%
\begin{theorem}[(\cite{Bona,CGS})]\label{j192}
Consider the meteor process on $\calC_k$ in the stationary regime.
Fix a family of patterns $\calA$, and
let $N$ be the number of sites in $\calC_k$ where $\calA$ occurs.
Then there exist $\mu,\sigma>0$ such that $(N- k \mu)/ \sigma\sqrt
{k}$ converges in distribution to the standard normal random variable
as $k\to\infty$.
\end{theorem}

\begin{pf}
We will supply a proof that is shorter than that in \cite{Bona} or
\cite{CGS}, Example~6.2, and illustrates well the power of the meteor
representation of craters and other patterns.

Let $\{U^j, j \in{\mathbb Z}\}$ be i.i.d. exponential random variables
with mean 1. Note that for any fixed $k$, the distribution of $\{U^j,
1\leq j \leq k\}$ is the same as that of $\{-T^j, 1\leq j \leq k\}$,
where $T^j$'s are defined relative to $\calC_k$.
Let $\xi_m$ be the indicator random variable of the occurrence of
$\calA$ at the $m$th site in
$\{U^j, j \in{\mathbb Z}\}$. In other words,
$\xi_j=1$ if and only if for some $\ba^r = a^r_1 \cdots a^r_{n_r} \in
\calA$, we\vspace*{1pt} have $U^{j+i-1} < U^{j+m-1}$ if and only if $a^r_i < a^r_m$
for all $1\leq i,m\leq n_r$.
Otherwise, $\xi_j=0$.

It is clear that the process $\{\xi_j, j\in{\mathbb Z}\}$ is stationary.

Let $b$ be the length of the longest pattern in $\calA$.
If $|j-m| > b$, then the occurrence of $\calA$ at site $j$ is
independent of the occurrence of $\calA$ at site $m$, since $U^n$'s
are independent. In other words, if $|j-m| > b$, then $\xi_j$ and $\xi
_m$ are independent. This implies that the process $\{\xi_j, j\in
{\mathbb Z}\}$ is $\varphi$-mixing in the sense of \cite{Bill}, Section~20.
The central limit theorem holds for $\sum_{j=1}^{k-b} \xi_j$,
according to \cite{Bill}, Theorem~20.1.
Let $N'$ be the number of sites $1\leq j \leq k-b$ in $\calC_k$ where\vspace*{2pt}
$\calA$ occurs, and
note that $N'$ has the same distribution as $\sum_{j=1}^{k-b} \xi_j$.
Hence, the central limit theorem holds for $N'$. Since $N$ and $N'$
differ by at most $b$,
the theorem follows.
\end{pf}

%
\begin{remark}
Theorem 20.1 of \cite{Bill} not only yields the central limit theorem
for $N$ in Theorem~\ref{j192} but also provides an effective
algorithm for calculating $\mu$ and~$\sigma$. To compute the values
of these parameters, one has to find $\mathbb{E}\xi_1$ and $\mathbb
{E}(\xi_1\xi_m)$ for all $m$. This is equivalent to counting the
corresponding permutations of length at most $2b$ [because we have
$\mathbb{E}(\xi_1\xi_m)=\mathbb{E}\xi_1\mathbb{E}\xi_m =
(\mathbb{E}\xi_1)^2$ for $|1-m| > b$]. For very small $b$, the
counting can be done directly. For moderate $b$, formulas such as those
in \cite{BBS} may be helpful, depending on the family of patterns
$\calA$.
\end{remark}

Craters represent sites that were hit by a meteor more recently than
their nearest neighbors.
We will now state a result about the locations of the
sites such that both of its neighbors were hit by meteors more recently
than the given site. Our result is
partly motivated by a technical application later in this section.

Recall that, according to Proposition~\ref{a298}(v), $m$ is a crater
if an only if $ T^{m-1} < T^m > T^{m+1}$.
We will say that $m$ is a \emph{mound} if an only if $ T^{m-1} > T^m <
T^{m+1}$.
Note that as we move along the graph $\calC_k$, we will encounter an
alternating sequence of craters and mounds, separated by stretches of
sites that are neither. The
craters and mounds correspond to the local maxima and minima of the
function $m\to T^m$.
Craters and mounds correspond to peaks and valleys of permutations.

%
\begin{proposition}\label{j203}
Consider the meteor process on $\calC_k$ in the stationary regime.
Let $B^k_{i,j}$ be the\vspace*{1pt} event that $2$ is a crater followed by a mound
and another crater, with $i$ and $j$ sites, respectively, between the
three distinguished sites. More precisely, for $i,j\geq0$ and $k \geq
i+j+5$, let
$B^k_{i,j}$ be the event that
$2$ and $i+j+4$ are craters, $i+3$ is a mound and $m$ is neither a
crater nor a mound for
$m=3,4,\ldots, i+2, i+4,\ldots, i+j+3$.
\begin{longlist}[(ii)]
\item[(i)] We have
%
%
\begin{eqnarray}
\label{j201} \mathbb{P}\bigl(B^k_{i,j}\bigr) &=&
\frac{2 (i+j+4)}{(i+j+5)!} \biggl[ \pmatrix{i+j+1
\cr
i+1}+(i+1) \pmatrix{i+j+2
\cr
i+2}
\biggr].
\end{eqnarray}

\item[(ii)] Recall events $F^k_n$ from Theorem~\ref{a285}. If $F^k_{n}$
holds, let $R$ denote the position of the unique mound between craters
at $2$ and $n+3$.
For any $\varepsilon>0$ there exist constants $c_1,c_2>0$ such that
for $n\geq1$ and $k\geq n+4$,
%
%
\begin{equation}
\label{j202} \mathbb{P}\bigl(|R/n - 1/2| > \varepsilon\mid F^k_n
\bigr) < c_1 e^{-c_2 n}.
\end{equation}
\end{longlist}
\end{proposition}

\begin{pf}
(i) This part follows from Proposition~23 of \cite{BBS}.
The argument is totally analogous to that in the proof of Theorem
\ref{a285}, so we leave the details to the reader. We just note
that one should take $m=i+3$ and $n=i+j+5$ in Proposition~23 of \cite{BBS}.

(ii)
The function $H(x):= - x \log x - (1-x) \log(1-x)$ is smooth on
$(0,1)$. It is elementary to check that it is increasing on $(0,1/2)$
and decreasing on $(1/2,1)$. Hence, for some $c_3,c_4>0$ and all $x\in(0,1)$,
%
%
\begin{equation}
\label{ag41} H(x) \leq H(1/2) - c_3 |x-1/2|^2 \leq\log
\biggl(\frac{ 2 }{1+ c_4
|x-1/2|^2} \biggr).
\end{equation}

By the Stirling approximation, for any $c_5 < 1 < c_6$, some $m_1$ and
all $m\geq m_1$, we have $c_5 m \log m < \log(m!) < c_6 m \log m$.
Fix any $\varepsilon>0$, and let $c_7>1$ be so small that $c_7\log
(2/(1+ c_4 \varepsilon^2)):= c_8
< \log2$.
For some $m_1$ and $r_1$, all $m\geq m_1$ and $r\geq r_1$ such that
$m-r \geq r_1$, we have
\begin{eqnarray*}
\log\pmatrix{m
\cr
r} &=& \log\biggl(\frac{m!}{r! (m-r)!} \biggr) \leq
c_7\bigl( m \log m - r \log r - (m-r) \log(m-r)\bigr)
\\
&=& c_7 m \biggl( - \frac{r} m \log\frac{r} m -
\biggl(1 - \frac{r}m \biggr) \log\biggl(1 - \frac{r} m \biggr)
\biggr).
\end{eqnarray*}
This and (\ref{ag41}) imply that if $m\geq m_1$, $r\geq r_1$, $m-r
\geq r_1$, $\varepsilon>0$ and
$|r/m - 1/2| > \varepsilon$,
\[
\log\pmatrix{m
\cr
r} \leq c_7 m\log\biggl(\frac{ 2 }{1+ c_4
|r/m-1/2|^2}
\biggr) \leq c_7 m \log\biggl(\frac{ 2 }{1+ c_4 \varepsilon^2} \biggr) =
c_8 m.
\]
If we take $m = i+j$, $r = i$ and we assume that
$|i - (i+j)/2| > \varepsilon(i+j)/2$, then the last estimate yields
for $i+j\geq m_1$ and $i,j\geq r_1$,
\[
\log\pmatrix{i+j
\cr
i} \leq c_8(i+j),
\]
and, therefore,
%
%
\begin{equation}
\label{ag42} \pmatrix{i+j
\cr
i} \leq e^{c_8(i+j)}.
\end{equation}

Note that for some polynomial $q_1$,
\[
\frac{2 (i+j+4)}{(i+j+5)!} \biggl[ \pmatrix{i+j+1
\cr
i+1}+(i+1) \pmatrix{i+j+2
\cr
i+2}
\biggr] \leq\frac{q_1(i+j)}{(i+j)!} \pmatrix{i+j
\cr
i}.
\]
This, (\ref{j201}) and (\ref{ag42}) give for
$i$ and $j$ satisfying $|i - (i+j)/2| > \varepsilon(i+j)/2$, $i+j\geq
m_1$ and $i,j\geq r_1$,
%
%
\begin{equation}
\label{j205} \mathbb{P}\bigl(B^k_{i,j}\bigr) \leq
\frac{ q_1(i+j) e^{c_8(i+j)}}{(i+j)! }.
\end{equation}
By changing the polynomial $q_1$, if necessary, we can drop the
assumptions that
$i+j\geq m_1$ and $i,j\geq r_1$.

Let
\begin{eqnarray*}
\Lambda(n,\varepsilon) &=& \bigl\{(i,j) \in{\mathbb Z}\dvtx  i,j \geq0, i+j+1
= n, \bigl|i
- (i+j)/2\bigr| > \varepsilon(i+j)/2\bigr\}.
\end{eqnarray*}
Recall that $c_8 < \log2$.
We obtain from (\ref{a296}) and (\ref{j205}) that for some \mbox{$c_1,c_2>0$},
\begin{eqnarray*}
\mathbb{P}\bigl(|R/n - 1/2| > \varepsilon\mid F^k_n
\bigr) &=& \bigl(\mathbb{P}\bigl(F^k_n\bigr)
\bigr)^{-1} \sum_{(i,j) \in\Lambda(n,\varepsilon)} \mathbb{P}
\bigl(B^k_{i,j}\bigr)
\\
&\leq&\frac{(n+4)!}{n (n+3) 2^{n+1}} \sum_{(i,j) \in\Lambda(n,\varepsilon)}
\frac{ q_1(i+j) e^{c_8(i+j)}}{(i+j)! }
\\
&=& \frac{(n+4)!}{n (n+3) 2^{n+1}} \sum_{(i,j) \in\Lambda(n,\varepsilon
)} \frac{ q_1(n-1) e^{c_8(n-1)}}{(n-1)! }
\\
&\leq&\frac{(n+4)!}{n (n+3) 2^{n+1}} n \frac{ q_1(n-1)
e^{c_8(n-1)}}{(n-1)! }
\\
& \leq& c_1
e^{-c_2n}.
\end{eqnarray*}
This completes the proof.
\end{pf}

The results and remarks presented so far in this section indicate
clearly that the crater process does not behave like a Poisson point
process on $\calC_k$. There are many ways to make this intuition
precise. Our next result shows that if there are very few craters, then
their positions are not approximately distributed as independent
uniform random variables on $\calC_k$, unlike in the case of a Poisson
point process.
We will prove that craters have a tendency to repel each other. This
``repulsion'' phenomenon is known in some other contexts; for example,
it applies to eigenvalues of random matrices \cite{Dyson} and other
determinantal processes \cite{HKPV}.

%
\begin{theorem}\label{j171}
Consider the meteor process on a circular graph $\calC_k$ with \mbox{$k\geq
3$}, and assume that the mass process $\{\mathcal{M}_t, t\in\mathbb
{R}\}$ is in the stationary regime.
Let $\mathcal{G}_1$ be the family of adjacent craters, that is, $(i,j)
\in\mathcal{G}_1$ if an only if there are craters at $i$ and $j$, and
there are no craters between $i$ and $j$.
We define $\mathcal{G}_2$ as the family of pairs $(i,j)$ such that
there is a crater at $i$ and a mound at $j$, or there is a mound at $i$
and a crater at $j$, and there are neither craters nor mounds between
$i$ and $j$.
For $r>2$, let
\begin{eqnarray*}
A^1_r &=& \biggl\{ \frac{ \max_{(i,j) \in\mathcal{G}_1} |i-j|}{\min
_{(i,j) \in\mathcal{G}_1} |i-j|} \leq r \biggr\},
\\
A^2_r &=& \biggl\{ \frac{ \max_{(i,j) \in\mathcal{G}_2} |i-j|}{\min
_{(i,j) \in\mathcal{G}_2} |i-j|} \leq r \biggr\}.
\end{eqnarray*}

\begin{longlist}[(ii)]
\item[(i)]
Let $H^1_n$ be the event that there are exactly $n$ craters at time 0.
For every \mbox{$n\geq2$}, $p<1$ and $r>2$ there exists $k_1 < \infty$ such
that for all $k\geq k_1$,
\mbox{$\mathbb{P}(A^1_r \mid H^1_n) > p$}.\vspace*{1pt}

\item[(ii)]
Let $H^2_n$ be the event that there are exactly $n$ craters and mounds
at time 0.
For every $n\geq2$, $p<1$ and $r>2$ there exists $k_1 < \infty$ such
that for all $k\geq k_1$, we have
$\mathbb{P}( A^2_r \mid H^2_n) > p$.
\end{longlist}
\end{theorem}

%
\begin{remark}
A combinatorial result in \cite{EM}, Theorem~6.1, shows that, assuming
that there are exactly $n$ craters, their most likely
configuration makes them equidistant from each other.
See also \cite{AK} for a closely related result.
These results are not equivalent to Theorem~\ref{j171} because the
probability that one of the most likely
configurations will occur does not have to be high.
\end{remark}

\begin{pf*}{Proof of Theorem~\ref{j171}}
Recall that, as we move around the graph $\calC_k$, we will encounter
an alternating sequence of craters and mounds, separated by stretches
of sites that are neither. Hence it is easy to see that part (ii)
implies (i). It remains to prove (ii).

Let $G_{i,j}$ be the event that there are craters at sites $i$ and $j$,
and there are no craters between these two sites. Given this event, let
$R_{i,j}$ be the distance from the unique mound between $i$ and $j$ to
the closest of these vertices. We define
$\widehat G_{i,j}$ and $\widehat R_{i,j}$ in an analogous way,
reversing the roles of craters and mounds.

Fix an $n\geq2$.
It is elementary to see that for every $r>2$ there exist $\varepsilon
>0$ and $c_1>0$ such that if $ A^2_r$ fails to hold, then the following
event must occur:
%
%
\begin{eqnarray}\label{j206}
&& \mathop{\bigcup_{i,j \in V}}_{|i-j| > c_1 k+1}
\bigl(G_{i,j} \cap\bigl\{\bigl|R_{i,j}/|i-j-1| - 1/2\bigr| > \varepsilon\bigr\}
\bigr)
\nonumber\\[-8pt]\\[-8pt]
&&\qquad \cup\mathop{\bigcup_{i,j \in V}}_{|i-j| > c_1 k+1} \bigl(
\widehat G_{i,j} \cap\bigl\{\bigl|\widehat R_{i,j}/|i-j-1| - 1/2\bigr| >
\varepsilon\bigr\} \bigr).\nonumber
\end{eqnarray}
If $G_{i,j}$ holds, then the value of $R_{i,j}$ does not depend on the
positions of craters and mounds outside the interval between $i$ and
$j$. Hence
\begin{eqnarray*}
&& \mathbb{P} \bigl(G_{i,j} \cap\bigl\{\bigl|R_{i,j}/|i-j-1| - 1/2\bigr| >
\varepsilon\bigr\} \mid H^2_n \bigr)
\\
&&\qquad = \mathbb{P}\bigl(G_{i,j} \mid H^2_n\bigr)
\mathbb{P}\bigl(\bigl|R_{i,j}/|i-j-1| - 1/2\bigr| > \varepsilon\mid
G_{i,j} \cap H^2_n\bigr)
\\
&&\qquad \leq\mathbb{P}\bigl(\bigl|R_{i,j}/|i-j-1| - 1/2\bigr| > \varepsilon\mid
G_{i,j} \cap H^2_n\bigr)
\\
&&\qquad = \mathbb{P}\bigl(\bigl|R_{i,j}/|i-j-1| - 1/2\bigr| > \varepsilon\mid
G_{i,j} \bigr).
\end{eqnarray*}
Proposition~\ref{j203}(ii) yields for some $c_2,c_3>0$ and $i$ and
$j$ such that $|i-j| > c_1 k+1$,
%
%
\begin{eqnarray}
\label{j207} && \mathbb{P} \bigl(G_{i,j} \cap\bigl\{\bigl|R_{i,j}/|i-j-1|
- 1/2\bigr| > \varepsilon\bigr\}\mid H^2_n \bigr)\nonumber
\\
&&\qquad \leq\mathbb{P}\bigl(\bigl|R_{i,j}/|i-j-1| - 1/2\bigr| > \varepsilon\mid
G_{i,j}\bigr)
\\
&&\qquad \leq c_2 e^{-c_3 |i-j-1|} \leq c_2 e^{-c_3 c_1 k}.
\nonumber
\end{eqnarray}
Interchanging the roles of craters and mounds, we obtain for $i$ and
$j$ such that $|i-j| > c_1 k+1$,
\begin{eqnarray*}
&& \mathbb{P} \bigl(\widehat G_{i,j} \cap\bigl\{\bigl|\widehat
R_{i,j}/|i-j-1| - 1/2\bigr| > \varepsilon\bigr\}\mid H^2_n
\bigr) \leq c_2 e^{-c_3 c_1 k}.
\end{eqnarray*}
This, (\ref{j207}) and (\ref{j206}) imply that
\[
\mathbb{P} \bigl( \bigl(A^2_r\bigr) ^c
\mid H^2_n \bigr) \leq2 \mathop{\sum
_{i,j \in V}}_{|i-j| > c_1 k+1} c_2 e^{-c_3 c_1 k} \leq2
k^2 c_2 e^{-c_3 c_1 k}.
\]
The last quantity goes to 0 as $k\to\infty$.
This completes the proof.
\end{pf*}

The last question that we are going to address in this section concerns
the age of the oldest exposed soil. A meteor hit displaces some soil,
and we can imagine that the displaced soil is placed on the top of the
soil already present at the site where it is deposited. Hence the age
of the oldest exposed soil is the minimum over all $n$ of
$\widetilde T^n:= \max( T^{n-1}, T^n, T^{n+1})$.

%
\begin{theorem}
For any $\varepsilon>0$ and $p<1$ there exists $k_1$ such that for
$k\geq k_1$,
\[
\mathbb{P} \Bigl(\Bigl\llvert\min_{1\leq n\leq k} \widetilde
T^n - (1/3)\log k\Bigr\rrvert< \varepsilon\log k \Bigr) >p.
\]
\end{theorem}

\begin{pf}
Consider any $\alpha\in(0, 2/3)$, and let $\beta= 2/3-\alpha>0$.
The probability that $T^n$ is among the $k^\alpha$ lowest values of $
\{T^j, 1\leq j \leq k\}$ is less than $2 k^\alpha/k = 2 k^{\alpha-1}$.
Hence, for a fixed $n$ and large $k$, the\vspace*{1pt} probability that
$\widetilde T^n$ is among the $k^\alpha$ lowest values of $ \{T^j,
1\leq j \leq k\}$ is less than $2( 2 k^{\alpha-1})^3 = 16 k^{3(\alpha
-1)} = 16 k^{-1 - 3\beta}$ (the dependence between the relevant events
is negligible for large $k$). It follows that the probability that
there exists a site $n$ such that $\widetilde T^n$ is among the
$k^\alpha$ lowest values of $ \{T^j, 1\leq j \leq k\}$ is less than $k
16 k^{-1 - 3\beta}=16 k^{ - 3\beta}$. The last quantity goes to 0 as
$k\to\infty$.

Consider any $\gamma\in( 2/3,1)$ and let $\lambda= \gamma-2/3>0$.
The probability that $T^n$ is among the $k^\gamma$ lowest values of $
\{T^j, 1\leq j \leq k\}$ is more than $ k^\gamma/(2k) = (1/2)
k^{\gamma-1}$.
Hence, for a fixed $n$ and large $k$, the probability that
$\widetilde T^n$ is among the $k^\gamma$ lowest values of $ \{T^j,
1\leq j \leq k\}$ is more than $(1/2)( (1/2) k^{\gamma-1})^3 = (1/16)
k^{3(\gamma-1)} = (1/16) k^{-1 + 3\lambda}$. It follows that the
probability that there exists a site $n$ such that $1\leq n= 3 i \leq
k$, $i\in{\mathbb Z}$ and $\widetilde T^n$ is among the $k^\gamma$
lowest values of $ \{T^j, 1\leq j \leq k\}$ is more
than $1- (1- (1/16) k^{-1 + 3\lambda})^{k/6}$. The\vspace*{1pt} last quantity goes
to 1 as $k\to\infty$.

Let $J$ be the rank of $\min_{1\leq j \leq k}\widetilde T^j$ among the
ordered values of $ \{T^j, 1\leq j \leq k\}$.
We have shown that for any $0<\alpha< 2/3 < \gamma<1$, we have
%
%
\begin{equation}
\label{j212} \lim_{k\to\infty}\mathbb{P}\bigl(k^\alpha< J <
k^\gamma\bigr)=1.
\end{equation}

Note that $ \{-T^j, 1\leq j \leq k\}$ are i.i.d., with the exponential
distribution with mean 1. Let $Y_{(n)}$ denote the $n$th order
statistic for $ \{-T^j, 1\leq j \leq k\}$.
It follows from \cite{dHF}, Theorem~2.2.1, that for any fixed $a\in
(0,1)$, random variables $k^{a/2} (Y_{(k- k^a)} - (1-a) \log k) $
converge weakly to the standard normal random variable as $k\to\infty$.
This and (\ref{j212}) easily imply the theorem.
\end{pf}

\section{Mass distribution}\label{secmass}

Section~\ref{craters} was concerned with the distribution of craters,
that is, sites where the mass $M^j$ is 0.
This section will present some results on the mass distribution at all
sites. In other words, we will consider the nondegenerate part of the
mass distribution at a site.

%
\begin{theorem}\label{m61}
Suppose that $d\geq1$, and let
$\{\mathcal{M}_t, t\geq0\} = \{(M^1_t,M^2_t,\ldots,\break M^{k}_t), t\geq
0\} $ be
the mass process on $G=\calC_n^d$ (the product of $d$ copies of the
cycle~$\calC_n$), under the stationary measure $Q_k$ (here $k=nd$).
Assume that \mbox{$\sum_{x\in V} M^x_0 = k$} under $Q_k$. We have
%
%
\begin{eqnarray}
\mathbb{E}_{Q_k} M^x_0 &=& 1, \qquad x\in V,
\label{ap111}
\\
 \lim_{k\to\infty} \operatorname{Var}_{Q_k}
M^x_0 &=& 1, \qquad x\in V, \label{ap112}
\\
\lim_{k\to\infty} \operatorname{Cov}_{Q_k}\bigl(
M^x_0, M^y_0\bigr) &=& -
\frac{1}{2d},\qquad x\leftrightarrow y,\label{ap113}
\\
\lim_{k\to\infty} \operatorname{Cov}_{Q_k}\bigl(
M^x_0, M^y_0\bigr) &=& 0, \qquad
x\ne y\mbox{ and } x\not\leftrightarrow y. \label{ap114}
\end{eqnarray}
\end{theorem}

\begin{pf}
By\vspace*{1pt} symmetry,
$\mathbb{E}_{Q_k} M^x_0 = \mathbb{E}_{Q_k} M^y_0$
for all $x,y\in V$. Since\break $\sum_{x\in V} M^x_0 = k$ under $Q_k$,
we must have $\mathbb{E}_{Q_k} M^x_0=1$ for $x\in V$. This proves
(\ref{ap111}).

We will base our estimates for $\operatorname{Var}_{Q_k} M^x_0 $ and
$\operatorname{Cov}_{Q_k}( M^x_0, M^y_0)$ on a\vspace*{1pt} representation of
$M^x_0$ using
WIMPs. Let $Z$ and $\widetilde Z$ be defined as $Z^1$ and $ Z^2$ in
Definition~\ref{wimps}.
In particular,
$\mathbb{P}(Z_0 = x) = \mathbb{P}(\widetilde Z_0 = x) = M^x_0/k$
for $x\in V$.

Note that since the state space $\calC_k^{2d}$ for the process
$(Z,\widetilde Z)$ is finite, the process has a stationary
distribution. The stationary distribution is unique because all states
communicate.
We will estimate the probability that $Z_t = \widetilde Z_t$ under the
stationary distribution. Let $\overline Z_t=Z_t-\widetilde Z_t$. It is
easy to see that $\overline Z_t$ is a Markov process (although a
function of a Markov process is not necessarily Markov).
The\vspace*{1pt} state space for $\overline Z_t$
may be identified with $V$ in the obvious way.
Let $ \{\pi_x, x\in V\}$ be the set of stationary probabilities for
the discrete time Markov chain $Z^*_j$ embedded in $\overline Z_t$.

First, we will discuss the case $d=1$.
We claim that, in this case, for some $c_1>0$, $\pi_1 = \pi_{n-1}=
c_1/2$ and $\pi_j = c_1$ for $j\ne1,n-1$.
It is easy to check that the following equations define the stationary
probabilities, and these equations are satisfied by the probabilities
specified above:
\begin{eqnarray*}
\pi_0 &=& \tfrac{1}2 \pi_0 + \tfrac{1}2
\pi_1 + \tfrac{1}2 \pi_{n-1},\qquad
\pi_1 = \tfrac{1}2 \pi_2,\qquad\pi_{n-1}
= \tfrac{1}2 \pi_{n-2},
\\
\pi_2 &=& \tfrac{1}2 \pi_3 + \tfrac{1}2
\pi_1 +\tfrac{1}4 \pi_0, \qquad
\pi_{n-2} = \tfrac{1}2 \pi_{n-3} + \tfrac{1}2
\pi_{n-1} +\tfrac{1}4 \pi_0,
\\
\pi_j & =& \tfrac{1}2 \pi_{j-1} + \tfrac{1}2
\pi_{j+1}, \qquad j \ne n-2, n-1, 0,1,2.
\end{eqnarray*}
Of course, $c_1$ is chosen so that $\sum_n \pi_n = 1$.
The mean holding time for $\overline Z_t$ is $1$ in the state $0$ and
it is $1/2$ in all other states.
This and the formulas for $\pi_j$'s imply that
%
%
\begin{eqnarray}
\label{ap122} \lim_{k\to\infty}(k/2)\mathbb{E}_{Q_k}
(Z_t =\widetilde Z_t) &=& \lim_{k\to\infty}(k/2)
\mathbb{E}_{Q_k} (\overline Z_t= 0)=1,\nonumber
\\
\lim_{k\to\infty}2k\mathbb{E}_{Q_k} (Z_t -
\widetilde Z_t=1) &=& \lim_{k\to\infty}2k
\mathbb{E}_{Q_k} (Z_t -\widetilde Z_t=-1)=1,
\\
\lim_{k\to\infty}k\mathbb{E}_{Q_k} (Z_t -
\widetilde Z_t=j) &=& 1, \qquad j \ne-1,0,1.
\nonumber
\end{eqnarray}

The case $d\geq2$ is similar but requires different notation.
Recall that $ \bzero=(0,\ldots,0)$.
Let $\bfa$ be set of all vertices $(a_1,\ldots, a_d) $ such
that $|a_i| = |a_j| = 1$ for some $i$ and $j$, and $a_m = 0$
for all $m\ne i,j$. Let $\bfb$ be set of all vertices $(b_1,\ldots, b_d) $ such that $|b_i| =2$ for some $i$, and $b_m = 0$
for all $m\ne i$. Let $\bfh$ be set of all vertices $(h_1,\ldots, h_d) $ such that $|h_i| =1$ for some $i$, and $h_m = 0$
for all $m\ne i$. Let $\bfg= V \setminus(\{\bzero\}\cup\bfa
\cup\bfb\cup\bfh)$.

We claim that for some $c_1>0$, $\pi_x = (1 - \frac{1}{2d}) c_1$
for all $x\in\bfh$ and $\pi_x = c_1$ for all other $x\in V$.
It is easy to check that the following equations define the
stationary probabilities, and these equations are satisfied by
the probabilities specified above:
\begin{eqnarray*}
\pi_\bzero&=& \frac{1}{2d} \pi_\bzero+ 2d
\frac{1}{2d} \pi_x, \qquad x \in\bfh,
\\
\pi_x &=& \frac{1}{2d} \pi_y + (2d -2)
\frac{1}{2d} \pi_z,\qquad x\in\bfh, y \in\bfb, z\in\bfa,
\\
\pi_x &=& 2 \biggl(\frac{1}{2d} \biggr)^2
\pi_\bzero+ 2\frac{1}{2d} \pi_y + (2d-2)
\frac{1}{2d} \pi_z, \qquad x\in\bfa, y \in\bfh, z\in\bfg,
\\
\pi_x &=& \biggl(\frac{1}{2d} \biggr)^2
\pi_\bzero+ \frac{1}{2d} \pi_y + (2d-1)
\frac{1}{2d} \pi_z, \qquad x\in\bfb, y \in\bfh, z\in\bfg,
\\
\pi_x &=& 2d\frac{1}{2d} \pi_y, \qquad x\in\bfg,
y \in\bfa\cup\bfb\cup\bfg.
\end{eqnarray*}
Recall that $c_1$ is chosen so that $\sum_n \pi_n = 1$,
the mean holding time for $\overline Z_t$ is $1$ in the state $0$ and
it is $1/2$ in all other states.
This and the formulas for $\pi_j$'s imply that
%
%
\begin{eqnarray}
\label{ap123} \lim_{k\to\infty}(k/2)\mathbb{E}_{Q_k}
(Z_t =\widetilde Z_t) &=& \lim_{k\to\infty}(k/2)
\mathbb{E}_{Q_k} (\overline Z_t= \bzero)=1,\nonumber
\\
\lim_{k\to\infty}\frac{2d}{2d-1}k\mathbb{E}_{Q_k}
(Z_t -\widetilde Z_t=x)&=&1, \qquad x \in\bfh,
\\
\lim_{k\to\infty}k\mathbb{E}_{Q_k} (Z_t -
\widetilde Z_t=x) &=&1, \qquad x \notin\{\bzero\} \cup\bfh.
\nonumber
\end{eqnarray}

Let $\alpha_\bzero=2$, $\alpha_x =1-\frac{1} {2d}$ for
$x\in\bfh$ and $\alpha_x = 1$ for all other $x$. By
(\ref{ap122}) and (\ref{ap123}), for any fixed $x\in V$ and
an arbitrarily small $\varepsilon>0$, there exists $k_1$ so large
that for any $k\geq k_1$, the probability that $Z_t - \widetilde Z_t =
x$ under the stationary distribution is in the interval
$((1-\varepsilon)\alpha_x/k, (1+\varepsilon)\alpha_x/k)$. Hence,
for $y\in
V$,
%
%
\begin{equation}
\label{ap121} \mathbb{P}_{Q_k}(Z_0 =y, \widetilde
Z_0 =y+x)  \in\bigl((1-\varepsilon)\alpha_x/k^2,
(1+\varepsilon)\alpha_x/k^2\bigr).
\end{equation}

Let $\mathcal{G}_t = \sigma(\mathcal{M}_s, 0\leq s \leq t)$.
It is easy to see that
\[
\mathbb{P}_{Q_k}(Z_0 =x \mid\mathcal{G}_0) =
M^x_0/k.
\]
The random variables $Z_0$ and $\widetilde Z_0$
are conditionally independent given $\mathcal{G}_0$, so
\[
\mathbb{P}_{Q_k}(Z_0 =y,\widetilde Z_0 =y+x
\mid\mathcal{G}_0) = M^y_0
M^{y+x}_0/k^2.
\]
Thus
\begin{eqnarray*}
\mathbb{E}_{Q_k} \bigl(M^y_0
M^{y+x}_0\bigr) &=& k^2 \mathbb{E}_{Q_k}
\mathbb{P}_{Q_k}(Z_0 =y,\widetilde Z_0 =y+x
\mid\mathcal{G}_0)
\\
&=& k^2\mathbb{P}_{Q_k}(Z_0
=y,\widetilde Z_0 =y+x ).
\end{eqnarray*}
This and (\ref{ap121}) yield,  for $k>k_1$,
\[
(1-\varepsilon)\alpha_x \leq\mathbb{E}_{Q_k}
\bigl(M^y_0 M^{y+x}_0\bigr) \leq(1+
\varepsilon)\alpha_x.
\]
Since $\varepsilon>0$ is arbitrarily small, it follows that
\[
\lim_{k\to\infty} \mathbb{E}_{Q_k} \bigl(M^y_0
M^{y+x}_0\bigr) = \alpha_x.
\]
For $x=\bzero$, we obtain
$\lim_{k\to\infty}
\mathbb{E}_{Q_k} (M^y_0)^2 = 2$.
This and (\ref{ap111}) imply that
\[
\lim_{k\to\infty} \operatorname{Var}_{Q_k}
M^y_0 =1.
\]
For $x\in\bfh$, we have $\lim_{k\to\infty} \mathbb{E}_{Q_k} (M^y_0
M^{y+x}_0) = 1-\frac{1} {2d}$, so, in view of (\ref{ap111}),
\[
\lim_{k\to\infty} \operatorname{Cov}_{Q_k}\bigl(
M^y_0, M^{y+x}_0\bigr) = -
\frac{1} {2d}.
\]
Finally, for $x\notin\{\bzero\} \cup\bfh$, we have
$\lim_{k\to\infty} \mathbb{E}_{Q_k} (M^y_0 M^{y+x}_0) = 1$, and,
therefore,
\[
\lim_{k\to\infty} \operatorname{Cov}_{Q_k}\bigl(
M^y_0, M^{y+x}_0\bigr) = 0.
\]
This completes the proof.
\end{pf}

%
\begin{remark}
(i) It
has been shown in \cite{meteor2} (a follow-up paper) that
the distributions of $M^1_0$ under the stationary measures $Q_k$
converge as $k\to\infty$.
We have neither explicit description nor detailed information about
the limit distribution.
We performed a number of long simulations.
Figure~\ref{fig1} illustrates some of the numerical results.
The figure on the left shows the empirical distribution of masses $\{
M^j_{10{,}000{,}000}, 1\leq j \leq60\mbox{,}000\}$, based\vspace*{1pt} on a single simulation
with ten million jumps (``meteor hits'') for a circular graph $\calC
_{60{,}000}$. The distribution has an atom at 0 of size about $1/3$, as
predicted by Theorem~\ref{a285}. The distribution does not appear to
have any other atoms.
The graph on the right shows the ``Q--Q'' plot (quantile on quantile
plot) for the continuous component of the empirical distribution of
masses versus the best matching gamma density (in the sense of matching
the first two moments), for a simulation on the graph $\calC_{6000}$.
The ``Q--Q'' plot shows convincingly that the distribution is not in the
gamma family.
We will return to this point in part~(iii) of this remark.

%
\begin{figure}[t]

\includegraphics{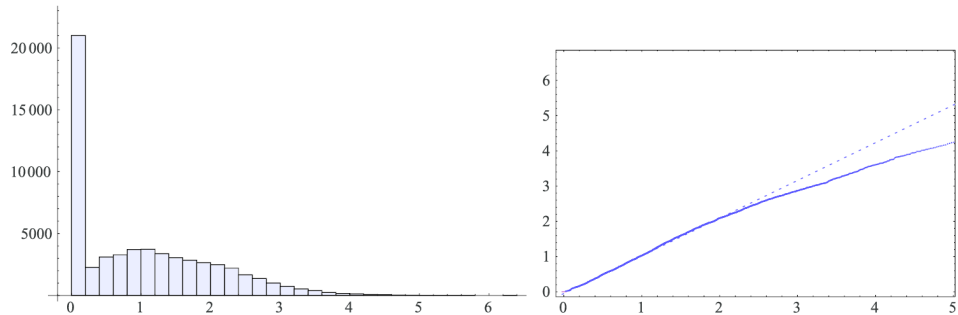}

\caption{The\vspace*{2pt} figure on the left shows the empirical distribution of
masses $\{M^j_{10,000,000},\break  1\leq j \leq60\mbox{,}000\}$, based on a single
simulation with ten million jumps (``meteor hits'') for a circular
graph $\calC_{60,000}$. The distribution has an atom at 0 of
(theoretical) size $1/3$.
The graph on the right shows the ``Q--Q'' plot (quantile on quantile
plot) for the continuous component of the empirical distribution of
masses versus the best matching gamma density (in the sense of matching
the first two moments), for a simulation on the graph $\calC_{6000}$.}\label{fig1}
\end{figure}

(ii) An argument similar to that in the proof of Theorem~\ref{m61}
leads to a (nonasymptotic) formula for the third moment of $M^1_0$,
for a fixed circular graph~$\calC_k$. The calculation is based on the
derivation of the stationary distribution for the Markov process
consisting of three dependent continuous time random walks. The
stationary distribution can be explicitly calculated using computer
algebra for low values of $k$. The values of the third moment of
$M^1_0$ seem to converge to $4.75531$ as $k$ goes to infinity. This
value is consistent with the results of computer simulations.

Calculating the stationary distribution for the Markov chain of three
random walks quickly becomes a time consuming task because the state
space of the Markov chain has $k^3$ elements, assuming that the cycle
has size $k$. To reduce the size of the state space, we collapsed the
states that were images of each other under symmetries of the cycle.
For example, for $k=20$, the state space size was reduced from $20^3 =
8000$ to $44$.

(iii)
It follows from our estimates that the limiting distribution of the
mass size at a given point, after removing the atom at 0, does not
belong to the gamma family. For a gamma random variable $X$ with density
$x^{\alpha- 1} \exp(-x / \beta) / (\Gamma(\alpha)\beta^\alpha)$,
we have $ EX^j = \beta^j \alpha(\alpha+1)\cdots(\alpha+ j -1)$. In
particular, $EX = \beta\alpha$, $EX^2 = \beta^2 \alpha(\alpha+1)$
and $EX^3 = \beta^3 \alpha(\alpha+1) (\alpha+2)$.
Let $W$ be $M^1_0$ conditioned to be nonzero.
In our case, under the stationary distribution $Q_k$, we have $EW =
3/2$, $EW^2 = 3$ and $EW^3 \approx4.755$.
If we have $EX = 3/2$ and $EX^2 = 3$ for a gamma distribution, then
$EX^3 = 7.5\ne4.755$.
There are no values of $\alpha$ and $\beta$ that would make the
moments of $W$ match the moments of a gamma distribution even in an
approximate sense.

(iv) Numerical calculations suggest that $(M^1_0)^2$ and $M^j_0$ are
asymptotically correlated,
when $k\to\infty$. Hence, it appears that $M^1_0$ and $M^j_0$ are
asymptotically dependent,
when $k\to\infty$.
We do not have a heuristic explanation for the lack of asymptotic
correlation of
$M^1_0$ and $M^j_0$, for $j\geq3$,
proved in Theorem~\ref{m61}.

(v) When $k=2$ or $3$, we can provide an explicit description of the
stationary distribution for the mass process $\mathcal{M}_t$ on the
circular graph $\calC_k$.
If $k=2$, then the stationary distribution of $\mathcal{M}_t$ has two
atoms of size $1/2$. One atom is the measure that gives mass 2 to site
1 and mass 0 to site 2. The other atom is the measure that gives mass 2
to site 2 and mass 0 to site 1.

Suppose that $k=3$ and for $j=1,2,3$, let $\mu_j$ be the random
measure which gives mass 0 to site $j$, $\mu_j(j+1)$ is the uniform
random variable on $[0,2]$ and $\mu_j(j+1)= 2-\mu_j(j+2)$. Then the
stationary distribution for $\mathcal{M}_t$ is the mixture, with equal
weights, of $\mu_j$, $j=1,2,3$. It is an elementary exercise to check
that the given measures are stationary.

(vi)
Consider the meteor process on a circular graph $\calC_k$,
and let\break $M^{1,n}_t = \sum_{j=1}^n M^{j\!}_{t\!}$. Then Theorem~\ref{m61}
implies that for any fixed $n$,\break
$\lim_{k\to\infty} \operatorname{Var}_{Q_k\!} M^{1,n}_0 = 1$. In
other words,
although the expected mass in an interval of length $n$, that is,
$\mathbb{E}_{Q_k} M^{1,n}_0 = n$, grows with $n$, the variance of
this mass does not grow (in the limit when $k\to\infty$).

More generally, consider the meteor process on the product $\calC_n^d$
of circular graphs,
and let $k=dn$. For a set $A\subset V$, let $M^A_t = \sum_{x\in A}
M^x_t$. Let $\partial A$ be the number of edges joining two vertices of
which exactly one is in $A$. Then Theorem~\ref{m61}
implies that for any fixed $A$,
$\lim_{k\to\infty} \operatorname{Var}_{Q_k} M^A_0 = |\partial A|/
(2d)$. Obviously, $\mathbb{E}_{Q_k} M^A_0 = |A|$.
The mass enclosed in each of the shapes in Figure~\ref{fig2} has the
same (asymptotic) variance.
%
%
\begin{figure}

\includegraphics{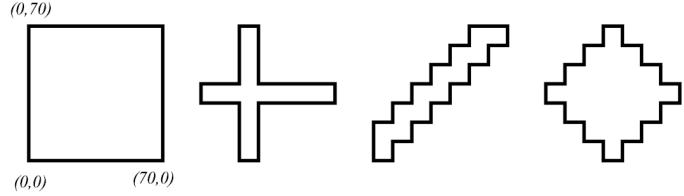}

\caption{All curves have the same height and width. They all have the
same ``boundary length'' $|\partial A|$, where $A$ denotes the set of
vertices inside the given closed curve.
The asymptotic variance
of the mass enclosed by one of these four curves has the same value as
for any other of these curves.}\label{fig2}
\end{figure}
\end{remark}

Consider the meteor process on a circular graph $\calC_k$ with $k\geq
4$ and assume that the mass process $\{\mathcal{M}_t, t\in\mathbb
{R}\}$ is in the stationary regime.
We will estimate the expected value of the height of a crater rim, that
is, the expected value of the mass at a site that is adjacent to a
crater. Note that the expected value of the mass at a uniformly chosen
noncrater is $3/2$ because the expected value of the mass at a site is
1 and the probability that a site is a crater is $1/3$.

%
\begin{proposition}
Consider the meteor process on a circular graph $\calC_k$ with $k\geq
6$, and assume that the mass process $\{\mathcal{M}_t, t\in\mathbb
{R}\}$ is in the stationary regime.
Then
%
%
\begin{equation}
\label{j312} 1.625 = 13/8 < \mathbb{E}_Q\bigl(M^2_0
\mid M^1_0 = 0\bigr) < 5/3 \approx1.667.
\end{equation}
\end{proposition}

\begin{pf}
Recall that $T^j$ denotes the time of the last jump of $N^j$ before 0,
that is, $T^j = \sup\{t\leq0\dvtx  N^j_t \ne N^j_{t-}\}$. The event $A:=\{
M^1_0 = 0\}$ is equivalent to $\{ T^2 < T^1 > T^k\}$. It is easy to see
that the conditional distribution of $T^1$ given $\{ T^2 < T^1 > T^k\}$
is the same as the distribution of $\max(T^k, T^1, T^2)$.
The density of $-\max(T^k, T^1, T^2)$ is $3e^{-3t}$.

The conditional distribution of $\mathcal{M}_{t-}$ given $A\cap\{
T^1=t\}$ is the stationary distribution because the event $A\cap\{
T^1=t\}$ is determined by $\{N^j_s, t\leq s \leq0\}$, and the value of
$\mathcal{M}_{t-}$
is determined by $\{N^j_s, s< t\}$.

Fix some $S>0$, and assume that $A\cap\{T^1=-S\}$ occurred.
For $t\in[0,S]$, let $F^j_t = \mathbb{E}(M^j_{-S+t} \mid A\cap\{
T^1=-S\})$, and note that $F^1_0=0$, $F^2_0=F^k_0 = 3/2$ and \mbox{$F^j_0 =
1$} for all other $j$.
Given $A\cap\{T^1=-S\}$, meteors hit sites $3,4,\ldots, k-1$ at a
constant rate of 1 hit per unit of time
during the time interval $(-S,0)$, so for $t\in[0,S]$,
%
%
\begin{eqnarray}
\label{j311} \frac{d}{dt} F^2_t &=&
\frac{1}2 F^{3}_t,\nonumber
\\
\frac{d}{dt} F^3_t &=& \frac{1}2
F^{4}_t - F^3_t,
\\
\frac{d}{dt} F^j_t &=& \frac{1}2
F^{j-1}_t + \frac{1}2 F^{j+1}_t
- F^j_t, \qquad j = 4,\ldots, k-2.
\nonumber
\end{eqnarray}
These equations and the initial conditions imply that $F^3_t > e^{-t}$
and, therefore, $F^2_t > 3/2 + (1-e^{-t})/2$ for $t\in(0,S]$. It
follows that
\begin{eqnarray*}
\mathbb{E}_Q\bigl(M^2_0 \mid
M^1_0 = 0\bigr) > \int_0^\infty
\bigl(3/2 + \bigl(1-e^{-s}\bigr)/2 \bigr) 3 e^{-3s} \,ds &=& 13/8.
\end{eqnarray*}

We also have $F^3_t < 1$ and, therefore, $F^2_t < 3/2 + t/2 $ for $t\in
(0,S]$. Hence
\[
\mathbb{E}_Q\bigl(M^2_0 \mid
M^1_0 = 0\bigr) < \int_0^\infty
(3/2 + s/2 ) 3 e^{-3s} \,ds = 5/3.
\]
This completes the proof.
\end{pf}

%
\begin{remark}
Computer simulations show that $\mathbb{E}_Q(M^2_0 \mid M^1_0 =
0)\approx1.6443$.
\end{remark}

We note that one can derive sharper estimates for $F^2_t $ using (\ref
{j311}) and hence sharper estimates in (\ref{j312}).

\section{Meteor processes on noncircular graphs}\label{secnonc}

Consider\vspace*{1pt} a circular graph~$\calC_k$, and suppose that the total mass
$\sum_{j=1}^k M^j_0$ is equal to $k$. Then
it is obvious that $\mathbb{E}_{Q_k} M^1_0 =1$
for every $k$, by symmetry. However, the fact that
$\lim_{k\to\infty} \operatorname{Var}_{Q_k} M^1_0 = 1$, proved in
Theorem~\ref{m61}, does not seem to be obvious. We will show that
under some structural assumptions on the graph $G$, the variance of
$M^x_0$ under the stationary distribution cannot be too large. We will
show that the bound for the variance of $M^x_0$ depends mainly on the
degree of the vertex.

A graph is called distance-transitive if for any two vertices $v$ and
$w$ at any distance $i$, and any other two vertices $x$ and $y$ at the
same distance, there is an automorphism of the graph that carries $v$
to $x$ and $w$ to $y$.

%
\begin{theorem}
Assume that $G$ is a distance-transitive $\rho$-regular graph.
Assume that $\sum_{x\in V} M^x_0=|V|$. Then under the stationary
distribution, for any $x\in V$,
\[
\operatorname{Var}_{Q} M^x_0 \leq
\frac{\rho+1 } {\rho-1 + 2\rho/ (|V|-1) }.
\]
\end{theorem}

\begin{pf}
By symmetry, $\mathbb{E}_{Q} M^x_0=1$, for all $x\in V$.

Let $Z$ and $\widetilde Z$ be defined as $Z^1$ and $ Z^2$ in Definition
\ref{wimps}.
In particular,
$\mathbb{P}(Z_0 = x) = \mathbb{P}(\widetilde Z_0 = x) = M^x_0/|V|$
for $x\in V$.

Note that since the state space $V^2$ for the process $(Z,\widetilde
Z)$ is finite, the process has a stationary distribution. The
stationary distribution is unique because all states communicate.
We will estimate the probability that $Z_t = \widetilde Z_t$ under the
stationary distribution.

Fix any vertex and label it $\bzero$. Let $Z^1$ be a continuous time
Markov process on~$V$ defined as follows. We let $Z^1_0$ be a vertex
uniformly chosen from all vertices $x$ with the property that the
distance from $x$ to $\bzero$ is the same as the distance from $Z_0$
to $\widetilde Z_0$. The\vspace*{1pt} process $Z^1$ jumps if an only if $(Z,
\widetilde Z)$ jumps. At a time $t$ of a jump of $(Z, \widetilde Z)$,
the process
$Z^1$ jumps to one of the nearest neighbors of $Z^1_{t-}$, whose
distance from $\bzero$ is the same as the distance between $Z_t$ and
$\widetilde Z_t$.
The process $Z^1$ is a continuous time Markov process on $V$, with the
mean holding time equal to $1/2$ at all vertices $x\ne\bzero$. The
mean holding time for $Z^1$ at $\bzero$ is
$(1- 1/\rho)^{-1}$. If $Z^1_t=\bzero$, the next jump it will take
will be to a vertex at a distance either 1 or 2 from $\bzero$. If
$Z^1_t\ne\bzero$, then the next jump
will be to a neighbor of $Z^1_t$.
Let $Z^2_t$ be a continuous time symmetric nearest neighbor random walk
on $V$, with the mean holding time equal to $1/2$ at all vertices $x\ne
\bzero$, and mean holding time at $\bzero$ equal to $(1- 1/\rho)^{-1}$.
The only difference between $Z^1$ and $Z^2$ is that $Z^2$ can jump from~$\bzero$ only to a nearest neighbor while $Z^1$ can jump from $\bzero
$ to some other vertices.

The long run proportion of time spent by $Z^2$ at $\bzero$ is
\[
\frac{(1- 1/\rho)^{-1}}{(1- 1/\rho)^{-1} + (|V| -1)/2} = \frac{\rho
}{\rho+ (\rho-1)( |V| -1)/2}.
\]
After every jump of $Z^1$ from $\bzero$, this process will take some
time, not necessarily zero, until it hits a neighbor of $\bzero$. Hence,
the long run proportion of time spent by $Z^1$ at $\bzero$ is less
than or equal to
\[
\frac{\rho}{\rho+ (\rho-1)( |V| -1)/2}.
\]
By symmetry, for any $x\in V$, the long run proportion of time spent by
$(Z,\widetilde Z) $ at $(x,x)$ is less than or equal to
\[
\frac{\rho}{|V|(\rho+ (\rho-1)( |V| -1)/2)}.
\]
Hence, for any $x\in V$,
%
%
\begin{equation}
\label{f281} \mathbb{P}_{Q}(Z_0 = \widetilde
Z_0 =x) = \frac{\rho}{|V|(\rho+ (\rho-1)( |V| -1)/2)}.
\end{equation}

Let $\mathcal{G}_t = \sigma(\mathcal{M}_s, 0\leq s \leq t)$.
Then, for $x\in V$,
\[
\mathbb{P}_{Q}(Z_0 =x \mid\mathcal{G}_0) =
M^x_0/|V|.
\]
The random variables $Z_0$ and $\widetilde Z_0$
are conditionally independent given $\mathcal{G}_0$, so
\[
\mathbb{P}_{Q}(Z_0 =\widetilde Z_0 =x \mid
\mathcal{G}_0) = \bigl(M^{x}_0/|V|
\bigr)^2.
\]
Thus
\begin{eqnarray*}
\mathbb{E}_{Q} \bigl(M^x_0
\bigr)^2 &=& |V|^2 \mathbb{E}_{Q}
\mathbb{P}_{Q}(Z_0 =\widetilde Z_0 =x \mid
\mathcal{G}_0) =|V|^2\mathbb{P}_{Q_k}(Z_0
=\widetilde Z_0 =x ).
\end{eqnarray*}
This and (\ref{f281}) yield
\[
\mathbb{E}_{Q} \bigl(M^x_0
\bigr)^2 \leq\frac{|V|^2 \rho}{|V|(\rho+ (\rho-1)( |V| -1)/2)}.
\]
Since $\mathbb{E}_{Q} M^x_0=1$, we obtain
\[
\operatorname{Var}_{Q} M^x_0 \leq
\frac{|V|^2 \rho}{|V|(\rho+ (\rho-1)( |V| -1)/2)} -1 = \frac{\rho+1 }{
\rho-1 + 2\rho/ (|V|-1) }.
\]
This completes the proof.
\end{pf}

Recall that $T^v_t$ denotes the time of the last jump of $N^v$ on the
interval $[0,t]$, with the convention that $T^v_t=-1$ if there were no
jumps on this interval.

%
\begin{theorem}\label{f171}
Suppose that $G$ is a complete graph with $k$ vertices  $\{1,2,\break \ldots,
k\}$, and recall the Poisson processes $N^m$.
Let $Q_k$ be the stationary distribution for the mass process $\mathcal
{M}$. When $k\to\infty$, processes $\{M^1_t -M^1_0 - t +
T^1_t, 
t\geq0\}$, under $Q_k$, converge weakly to the process identically
equal to 0 in the Skorokhod space $D([0,\infty), \mathbb{R})$.
\end{theorem}

%
\begin{corollary}\label{f172}
Under assumptions of Theorem~\ref{f171}, we have the following:
\begin{longlist}[(ii)]
\item[(i)] the distributions of $M^1_0$ under $Q_k$ converge to the
exponential distribution with mean 1, when $k\to\infty$;

\item[(ii)] there is propagation of chaos; that is,
for any finite $n\geq2$, the distributions of $\{M^j_t, t\geq0\}$,
$j=1,\ldots,n$, are asymptotically independent, when $k\to\infty$.
\end{longlist}
\end{corollary}

\begin{pf*}{Proof of Theorem~\ref{f171}}
Let $x^+ = \max(x,0)$.
Let $R^k_{s,t}$ be the mass moved to state 1 during the time interval
$[s,t]$, that is,
$R^k_{s,t}= \sum_{u\in[s,t]} (M^1_u - M^1_{u-})^+$.
If $t_1$~and~$t_2$ are any two consecutive jumps\vspace*{1pt} of $N^1$, then
$M^1_{t_1} = 0$ and $M^1_u - M^1_{t_1} = R^k_{u,t_1}$ for all $u \in
(t_1, t_2)$. Hence it will suffice to prove that for any two fixed
rational numbers
$0<t_1,t_2<\infty$,
$R^k_{t_1,t_2} $ converges to $t_2 - t_1$ weakly, as $k\to\infty$.

Let $Z^j$'s be defined as in Definition~\ref{wimps}.
Let $\mathcal{G}_t = \sigma(\mathcal{M}_s, 0\leq s \leq t)$.
Then for any $t\geq0$, a.s.,
\[
\mathbb{P}_{Q_k}\bigl(Z^j_t =1 \mid
\mathcal{G}_t\bigr) = M^1_t/k.
\]
The processes $\{Z^j_s, s\in[0,t]\}$, $j\geq1$,
are conditionally independent given $\mathcal{G}_t$, so by the law of
large numbers, for every $t\geq0$, a.s.,
%
%
\begin{equation}
\label{f173} \lim_{n\to\infty} \frac{1}n \sum
_{j=1}^n \bone_{\{Z^j_t =1\}} =
M^1_t/k.
\end{equation}
Since the process $M^1$ has only a finite number of jumps on any finite
time interval, the convergence in (\ref{f173}) holds uniformly on
every interval of the form $[t_1, t_2]$, with $0< t_1 < t_2 < \infty$.
Fix any $0< t_1 < t_2 < \infty$
and let
\[
A(k,n) = \frac{k}n \sum_{u \in[t_1,t_2]} \sum
_{j=1}^n \bone_{\{Z^j_u =1, Z^j_{u-} \ne1\}}.
\]
In view of earlier remarks, it will suffice to prove that, in probability,
\[
\lim_{k\to\infty} \lim_{n\to\infty} A(k,n) =
t_2 - t_1.
\]
It will be enough to show that
%
%
\begin{equation}
\label{f174} \lim_{k\to\infty} \lim_{n\to\infty}
\mathbb{E}_{Q_k} A(k,n) = t_2 - t_1
\end{equation}
and
%
%
\begin{equation}
\label{f175} \lim_{k\to\infty} \lim_{n\to\infty}
\operatorname{Var}_{Q_k} A(k,n) = 0.
\end{equation}

Since all $Z^j$'s have the same distribution,
to prove (\ref{f174}), it will suffice to show that
%
%
\begin{equation}
\label{f176} \lim_{k\to\infty} k \mathbb{E}_{Q_k} \sum
_{u \in[t_1,t_2]} \bone_{\{Z^1_u =1, Z^1_{u-} \ne1\}} = t_2 -
t_1.
\end{equation}
By symmetry, $\mathbb{P}_{Q_k}(Z^1_{t_1} = 1) = 1/k$. After the
process $Z^1$ jumps to some other state, it has probability less than
$1- e^{- (t_2 - t_1)/(k-1)}$ of jumping to 1 in the remaining time in
the interval $[t_1, t_2]$.
If it jumps back to 1 and then again to another state, it has, once again,
probability less than $1- e^{- (t_2 - t_1)/(k-1)}$ of jumping to 1 in
the remaining time in the interval $[t_1, t_2]$. A\vspace*{1.5pt} similar argument
applies to further possible jumps to 1.
Hence, if we denote consecutive jumps of $Z^1$ to the state 1 on the
interval $(t_1, t_2]$ by $S_1, S_2,\ldots,$ then
%
%
\begin{eqnarray}
\label{f1712}
\mathbb{P}_{Q_k}\bigl(Z^1_{t_1} =
1, S_m \leq t_2\bigr) &\leq&\frac{1}k \bigl(1-
e^{- (t_2 - t_1)/(k-1)} \bigr)^m
\nonumber\\[-8pt]\\[-8pt]
&\leq& (t_2 - t_1)^m (k-1)^{-m-1},\nonumber
\end{eqnarray}
and, therefore,
\[
\mathbb{E}_{Q_k} \sum_{u \in(t_1,t_2]}
\bone_{\{Z^1_{t_1} = 1\}}\bone_{\{Z^1_u =1, Z^1_{u-} \ne1\}} \leq\sum
_{m\geq1}(t_2
- t_1)^m (k-1)^{-m-1}.
\]
This implies that
%
%
\begin{equation}
\label{f177} \lim_{k\to\infty} k \mathbb{E}_{Q_k} \sum
_{u \in(t_1,t_2]} \bone_{\{Z^1_{t_1} = 1\}} \bone_{\{Z^1_u =1, Z^1_{u-}
\ne1\}} =
0.
\end{equation}

Next consider the case when $Z^1_{t_1} \ne1$. The probability that the
process $Z^1$ jumps to 1 before $t_2$ is equal to $1- e^{- (t_2 -
t_1)/(k-1)}$, so
%
%
\begin{equation}
\label{f1713} \mathbb{P}_{Q_k}\bigl(Z^1_{t_1}
\ne1, S_1 \leq t_2\bigr) = \frac{k-1}k \bigl(1-
e^{- (t_2 - t_1)/(k-1)} \bigr),
\end{equation}
and, consequently,
%
%
\begin{equation}
\label{f178} \lim_{k\to\infty} k \mathbb{E}_{Q_k} (
\bone_{\{Z^1_{t_1} \ne1\}} \bone_{\{S_1 \leq t_2\}} ) = t_2 - t_1.
\end{equation}
By the strong Markov property applied at $S_1$ and (\ref{f177}),
%
%
\begin{equation}
\label{f179} \lim_{k\to\infty} k \mathbb{E}_{Q_k} \sum
_{m\geq2} \bone_{\{Z^1_{t_1} \ne1\}} \bone_{\{S_m \leq t_2\}} =
0.
\end{equation}
We combine (\ref{f177}) and (\ref{f178})--(\ref{f179}) to see
that (\ref{f176}) holds, and therefore, (\ref{f174}) is true.

Given (\ref{f174}),
in order to prove (\ref{f175}), it is necessary and sufficient to
show that
%
%
\begin{equation}
\label{f1710} \lim_{k\to\infty} \lim_{n\to\infty}
\mathbb{E}_{Q_k} A(k,n)^2 \leq(t_2-t_1)^2.
\end{equation}
Let $S^j_1, S^j_2,\ldots$
denote\vspace*{1pt} the consecutive jumps of $Z^j$ to the state 1 on the interval
$(t_1, t_2]$.
We have
%
%
\begin{eqnarray}
\label{f1711} A(k,n)^2 &=& \biggl(\frac{k}n
\biggr)^2 \sum_{u,v \in[t_1,t_2]} \sum
_{i,j=1}^n \bone_{\{Z^i_u =1, Z^i_{u-} \ne1, Z^j_v =1, Z^j_{v-} \ne1\}}\nonumber
\\
&=& \biggl(\frac{k}n \biggr)^2 \sum
_{m,r\geq1} \sum_{i,j=1}^n
\bone_{\{S^i_m \leq t_2, S^j_r \leq t_2\}}\nonumber
\\
&=& \biggl(\frac{k}n \biggr)^2 \sum
_{m,r\geq1} \sum_{j=1}^n
\bone_{\{S^j_m \leq t_2, S^j_r \leq t_2\}}
\nonumber
\\[-8pt]
\\[-8pt]
&&{}  + \biggl(\frac{k}n \biggr)^2
\sum _{m,r\geq1} \sum_{i,j=1, i\ne j}^n
\bone_{\{S^i_m \leq t_2, S^j_r \leq t_2\}}
\nonumber
\\
&\leq&\biggl(\frac{k}n \biggr)^2 \sum
_{j=1}^n \bone_{\{S^j_1 \leq t_2\}} + \biggl(
\frac{k}n \biggr)^2 \sum_{m\geq2}
2m \sum_{j=1}^n \bone_{\{S^j_m \leq t_2\}}
\nonumber
\\
&&{}+ \biggl(\frac{k}n \biggr)^2 \sum
_{m,r\geq1} \sum_{i,j=1, i\ne j}^n
\bone_{\{S^i_m \leq t_2, S^j_r \leq t_2\}}.
\nonumber
\end{eqnarray}

Combining (\ref{f1712}) and (\ref{f1713}), we obtain
\[
\mathbb{P}_{Q_k}\bigl( S^1_1 \leq
t_2\bigr) \leq(t_2 - t_1)
(k-1)^{-2} + \frac{k-1}k \bigl(1- e^{- (t_2 - t_1)/(k-1)} \bigr),
\]
which has a finite value for each $k$,
so
%
%
\begin{equation}
\label{f1714} \lim_{k\to\infty} \lim_{n\to\infty}
\mathbb{E}_{Q_k} \Biggl[ \biggl(\frac{k}n \biggr)^2
\sum_{j=1}^n \bone_{\{S^j_1 \leq t_2\}}
\Biggr] = 0.
\end{equation}

By (\ref{f1712}) and the strong Markov property applied at $S^j_1$,
we have for $m\geq2$,
\[
\mathbb{P}_{Q_k}\bigl( S^j_m \leq
t_2\bigr) \leq(t_2 - t_1)^{m-1}
(k-1)^{-m},
\]
so
%
%
\begin{eqnarray}
\label{f1715}
\quad &&\lim_{k\to\infty} \lim_{n\to\infty}
\mathbb{E}_{Q_k} \Biggl[ \biggl(\frac{k}n \biggr)^2
\sum_{m\geq2} 2m \sum_{j=1}^n
\bone_{\{S^j_m \leq t_2\}} \Biggr]
\nonumber\\[-8pt]\\[-8pt]
&&\qquad \leq \lim_{k\to\infty} \lim_{n\to\infty} \Biggl[
\biggl(\frac{k}n \biggr)^2 \sum_{m\geq2}
2m \sum_{j=1}^n (t_2 -
t_1)^{m-1} (k-1)^{-m} \Biggr] = 0.\nonumber
\end{eqnarray}

In view of (\ref{f1711})--(\ref{f1715}), to complete the proof of
(\ref{f1710}),
it remains to show that
%
%
\begin{equation}
\label{f181} \qquad\lim_{k\to\infty} \lim_{n\to\infty}
\mathbb{E}_{Q_k} \Biggl[ \biggl(\frac{k}n \biggr)^2
\sum_{m,r\geq1} \sum_{i,j=1, i\ne j}^n
\bone_{\{S^i_m \leq t_2, S^j_r \leq t_2\}} \Biggr] \leq(t_2 - t_1)^2.
\end{equation}
Since the joint distribution of $(Z^i, Z^j)$ does not depend on $i$ and
$j$ as long as $i\ne j$, (\ref{f181}) will follow once we prove
%
%
\begin{equation}
\label{f182} \lim_{k\to\infty} \mathbb{E}_{Q_k} \biggl[
k^2 \sum_{m,r\geq1} \bone_{\{S^1_m \leq t_2, S^2_r \leq t_2\}}
\biggr] \leq(t_2 - t_1)^2.
\end{equation}

We will estimate the proportion of time
that $Z^1$ and $Z^2$ spend in the same state. After the two processes
meet, they spend an exponential amount of time together, with mean one,
and then they jump at the same time. They jump to the same state with
probability $1/(k-1)$ and if they do, they spend another period of
exponential length in the same state. The sequence of jumps to the same
state has geometric length with expectation $(k-1)/(k-2)$, so the total
time the processes spend together before they separate has expectation
$(k-1)/(k-2)$. When the processes travel through separate states, each
one jumps to the state occupied by the
other process at the rate $1/(k-1)$, so the waiting time for the next
meeting at some state is exponential with mean $(k-1)/2$. It follows
that in the long run, the proportion of time the two processes are in
the same state is
%
%
\begin{equation}
\label{f1810} \frac{(k-1)/(k-2)}{(k-1)/(k-2) + (k-1)/2} = \frac{2}k.
\end{equation}
By symmetry,
the proportion of time spent by the two processes in state 1 is
$2/k^2$, so
\[
\mathbb{P}_{Q_k}\bigl(Z^1_{t_1} =1,
Z^2_{t_1} =1\bigr) = 2/k^2.
\]
This and the argument given in support of (\ref{f1712}) can be
combined to see that
\begin{eqnarray*}
\mathbb{P}_{Q_k}\bigl(Z^1_{t_1} =
Z^2_{t_1} =1, S^1_m \leq
t_2\bigr)
&\leq& \frac{2}{k^2} \bigl(1- e^{- (t_2 - t_1)/(k-1)}
\bigr)^m
\\
& \leq& 2(t_2 - t_1)^m
(k-1)^{-m-2},
\end{eqnarray*}
and, therefore,
\begin{eqnarray*}
&& \mathbb{E}_{Q_k} \biggl[ k^2 \sum
_{m,r\geq1} \bone_{\{Z^1_{t_1} = Z^2_{t_1} =1,
S^1_m \leq t_2, S^2_r \leq t_2\}} \biggr] 
\\
&&\qquad \leq 2
\mathbb{E}_{Q_k} \biggl[ k^2 \sum
_{ 1\leq r \leq m } \bone_{\{Z^1_{t_1} = Z^2_{t_1} =1,
S^1_m \leq t_2, S^2_r \leq t_2\}} \biggr]
\\
&&\qquad \leq 2 \mathbb{E}_{Q_k} \biggl[ m k^2 \sum
_{ m\geq1 } \bone_{\{Z^1_{t_1} = Z^2_{t_1} =1,
S^1_m \leq t_2\}} \biggr]
\\
&&\qquad \leq 4 k^2
\sum_{m\geq1} m (t_2 - t_1)^m
(k-1)^{-m-2}.
\end{eqnarray*}
This implies that
%
%
\begin{equation}
\label{f183} \lim_{k\to\infty} \mathbb{E}_{Q_k} \biggl[
k^2 \sum_{m,r\geq1} \bone_{\{Z^1_{t_1} = Z^2_{t_1} =1,
S^1_m \leq t_2, S^2_r \leq t_2\}}
\biggr] = 0.
\end{equation}

We will now estimate
$\mathbb{P}_{Q_k} (Z^1_{t_1} =1,
S^1_m \leq t_2, S^2_r \leq t_2 )$. By symmetry,\break $\mathbb
{P}_{Q_k} (Z^1_{t_1} =1 )=1/k$. Consider the case $m=r=1$,
and suppose that $Z^1_{t_1}=1$.
After the process $Z^1$ jumps to some other state, it has probability
less than $1- e^{- (t_2 - t_1)/(k-1)}$ of jumping\vspace*{1pt} to 1 in the remaining
time in the interval $[t_1, t_2]$. The probability that $Z^2$ will jump
to 1 from some other state during $[t_1, t_2]$ is bounded by $1- e^{-
(t_2 - t_1)/(k-1)}$, no matter where $Z^2$ is at the time $t_1$. Hence,
the probability that at least one of the processes $Z^1$ or $Z^2$ jumps
to 1 from some other state during $[t_1, t_2]$ is bounded by $2(1- e^{-
(t_2 - t_1)/(k-1)})$. Now we consider two cases. The first one is that
at the time of the first jump of $Z^1$ or $Z^2$ to 1 from some other
state; the other process jumps as well. The conditional probability
that the second one will also jump to 1 is $1/(k-1)$. The second case
is that the other process does not jump at the same time. The\vspace*{1pt}
probability that it will jump to 1 in the remaining time in $[t_1,
t_2]$ is bounded by $1- e^{- (t_2 - t_1)/(k-1)}$. Altogether,
%
%
\begin{eqnarray}\label{f184}
&& \mathbb{P}_{Q_k} \bigl(Z^1_{t_1} =1,
S^1_1 \leq t_2, S^2_1
\leq t_2 \bigr)\nonumber
\\
&&\qquad \leq\frac{1}k\cdot2\bigl(1-
e^{- (t_2 - t_1)/(k-1)}\bigr) \biggl( \frac{1} {k-1} + 1- e^{- (t_2 -
t_1)/(k-1)} \biggr)
\\
&&\qquad \leq  c(t_1,t_2) \frac{1} {(k-1)^3}.\nonumber
\end{eqnarray}
The same argument that proves (\ref{f1712}) gives for any $n$ and
$m\geq1$,
%
%
\begin{eqnarray}\label{f185}
\mathbb{P}_{Q_k}\bigl(Z^1_{t_1}
=1, S^1_m \leq t_2 \bigr) &\leq&
\frac{1}k \bigl(1- e^{- (t_2 - t_1)/(k-1)} \bigr)^m
\nonumber\\[-8pt]\\[-8pt]
&\leq&(t_2 - t_1)^m (k-1)^{-m-1}\nonumber
\end{eqnarray}
and
%
%
\begin{eqnarray}\label{f186}
\mathbb{P}_{Q_k}\bigl(Z^1_{t_1}
=1, S^2_m \leq t_2 \bigr) &\leq&
\frac{1}k \bigl(1- e^{- (t_2 - t_1)/(k-1)} \bigr)^m
\nonumber\\[-8pt]\\[-8pt]
&\leq& (t_2 - t_1)^m (k-1)^{-m-1}.\nonumber
\end{eqnarray}
We combine (\ref{f184})--(\ref{f186}) to see that
\begin{eqnarray*}
&&\mathbb{E}_{Q_k} \biggl[ k^2 \sum
_{m,r\geq1} \bone_{\{Z^1_{t_1} =1,
S^1_m \leq t_2, S^2_r \leq t_2\}} \biggr]
\\
&&\qquad  \leq\mathbb{E}_{Q_k} \bigl[ k^2 \bone_{\{Z^1_{t_1} =1,
S^1_1 \leq t_2, S^2_1 \leq t_2\}}
\bigr] + \mathbb{E}_{Q_k} \biggl[ k^2 \sum
_{ 1\leq r \leq m, m\geq2 } \bone_{\{Z^1_{t_1} =1,
S^1_m \leq t_2, S^2_r \leq t_2\}} \biggr]
\\
&&\quad\qquad{} + \mathbb{E}_{Q_k} \biggl[ k^2 \sum
_{ 1\leq m \leq r, r\geq2 } \bone_{\{Z^1_{t_1} =1,
S^1_m \leq t_2, S^2_r \leq t_2\}} \biggr]
\\
&&\qquad \leq\mathbb{E}_{Q_k} \bigl[ k^2 \bone_{\{Z^1_{t_1} =1,
S^1_1 \leq t_2, S^2_1 \leq t_2\}}
\bigr] + \mathbb{E}_{Q_k} \biggl[ m k^2 \sum
_{ m\geq2 } \bone_{\{Z^1_{t_1} =1,
S^1_m \leq t_2\}} \biggr]
\\
&&\qquad\quad{} + \mathbb{E}_{Q_k} \biggl[ rk^2 \sum
_{ r\geq2 } \bone_{\{Z^1_{t_1} =1, S^2_r \leq t_2\}} \biggr]
\\
&&\qquad \leq k^2 c(t_1,t_2) \frac{1} {(k-1)^3} + 2
k^2 \sum_{m\geq2} m (t_2 -
t_1)^m (k-1)^{-m-1}.
\end{eqnarray*}
This implies that
%
%
\begin{equation}
\label{f187} \lim_{k\to\infty} \mathbb{E}_{Q_k} \biggl[
k^2 \sum_{m,r\geq1} \bone_{\{Z^1_{t_1} =1,
S^1_m \leq t_2, S^2_r \leq t_2\}}
\biggr] = 0.
\end{equation}
By symmetry,
%
%
\begin{equation}
\label{f188} \lim_{k\to\infty} \mathbb{E}_{Q_k} \biggl[
k^2 \sum_{m,r\geq1} \bone_{\{Z^2_{t_1} =1,
S^1_m \leq t_2, S^2_r \leq t_2\}}
\biggr] = 0.
\end{equation}

It follows from (\ref{f1810}) that
$\mathbb{P}_{Q_k} (Z^1_{t_1} =Z^2_{t_1} \ne1
)=2(k-1)/k^2$. The reasoning completely analogous to that given in the
case when $Z^1_{t_1} = 1$ yields
%
%
\begin{equation}
\label{f1811} \lim_{k\to\infty} \mathbb{E}_{Q_k} \biggl[
k^2 \sum_{m,r\geq1} \bone_{\{Z^1_{t_1} =Z^2_{t_1} \ne1,
S^1_m \leq t_2, S^2_r \leq t_2\}}
\biggr] = 0.
\end{equation}

Finally,
consider the event $F:=\{Z^1_{t_1} \ne Z^2_{t_1}, Z^1_{t_1}\ne1,
Z^2_{t_1}\ne1\}$. The probability that $Z^1$ will jump to 1 during
$[t_1, t_2]$ is equal to $1- e^{- (t_2 - t_1)/(k-1)}$. Let $\tau= \inf
\{t\geq t_1\dvtx  Z^1_t = Z^2_t\}$ ($\tau= t_2$ if the two processes do not
meet before $t_2$).
We have
%
%
\begin{eqnarray}\label{f1813}
&& \mathbb{P}_{Q_k}  \bigl(F, S^1_1 \leq t_2, S^2_1 \leq t_2 \bigr)\nonumber
\\
&&\qquad = \mathbb{P}_{Q_k} \bigl(F, S^1_1 < \tau<
S^2_1 \leq t_2 \bigr) +\mathbb{P}_{Q_k}
\bigl(F, S^2_1 < \tau< S^1_1 \leq
t_2 \bigr)
\nonumber\\[-8pt]\\[-8pt]
&&\quad\qquad{} +\mathbb{P}_{Q_k} \bigl(F, \tau< S^1_1 \leq
t_2, \tau< S^2_1 \leq t_2 \bigr)\nonumber
\\
&&\quad\qquad{} +\mathbb{P}_{Q_k} \bigl(F, S^1_1 \leq
t_2, S^2_1 \leq\tau\land t_2\bigr).
\nonumber
\end{eqnarray}
Our usual estimates give
\[
\mathbb{P}_{Q_k} \bigl(F, S^2_1 < \tau<
S^1_1 \leq t_2 \bigr) \leq\bigl(1-
e^{- (t_2 - t_1)/(k-1)} \bigr)^3 \leq(t_2 -
t_1)^3(k-1)^{-3},
\]
so
%
%
\begin{equation}
\label{f191} \qquad\lim_{k\to\infty} \mathbb{E}_{Q_k} \bigl[
k^2 \bone_{F\cap\{S^2_1 < \tau< S^1_1 \leq t_2\}} \bigr] \leq\lim_{k\to
\infty}
k^2 (t_2 - t_1)^3(k-1)^{-3}
= 0,
\end{equation}
and, by symmetry,
%
%
\begin{equation}
\label{ma131} \lim_{k\to\infty} \mathbb{E}_{Q_k} \bigl[
k^2 \bone_{F\cap\{S^1_1 < \tau< S^2_1 \leq t_2\}} \bigr] = 0.
\end{equation}

Given\vspace*{1pt} $\{Z^1_t, t\in[t_1, t_2]\}$, the conditional probability that
$Z^2$ jumps to 1 before or at time $\tau\land t_2$ is bounded by
$1- e^{- (t_2 - t_1)/(k-1)}$. It follows that
\[
\mathbb{P}_{Q_k} \bigl(F, S^1_1 \leq
t_2, S^2_1 \leq\tau\land t_2
\bigr) \leq\bigl(1- e^{- (t_2 - t_1)/(k-1)} \bigr)^2 \leq(t_2 -
t_1)^2(k-1)^{-2}
\]
and
%
%
\begin{eqnarray}
\label{f192}
\lim_{k\to\infty} \mathbb{E}_{Q_k} \bigl[
k^2 \bone_{F\cap\{S^1_1 \leq t_2, S^2_1 \leq\tau\land t_2\}} \bigr]
&\leq& \lim_{k\to\infty}
k^2 (t_2 - t_1)^2(k-1)^{-2}
\nonumber\\[-8pt]\\[-8pt]
&=& (t_2 - t_1)^2.\nonumber
\end{eqnarray}

The probability that, given $F$, the coupling time $\tau$ will occur
before $t_2$ is bounded by
$1- e^{- (t_2 - t_1)/(k-1)} \leq(t_2 - t_1)/(k-1)$, so\vspace*{1pt} using the
strong Markov property at $\tau$,
the case of $F\cap\{\tau< S^1_1 \leq t_2, \tau< S^2_1 \leq t_2\}$
is reduced to that in~(\ref{f1811}), and
we obtain the following bound:
%
%
\begin{equation}
\label{f193} \lim_{k\to\infty} \mathbb{E}_{Q_k} \bigl[
k^2 \bone_{F\cap\{\tau< S^1_1 \leq t_2, \tau< S^2_1 \leq t_2\}} \bigr] =0.
\end{equation}

In view of (\ref{f1813}), estimates (\ref{f191})--(\ref{f193}) yield
%
%
\begin{equation}
\label{f194} \lim_{k\to\infty} \mathbb{E}_{Q_k} \bigl[
k^2 \bone_{F\cap\{ S^1_1 \leq t_2, S^2_1 \leq t_2\}} \bigr] \leq(t_2 -
t_1)^2.
\end{equation}
A similar analysis gives
%
%
\begin{equation}
\label{f195} \lim_{k\to\infty} \mathbb{E}_{Q_k} \bigl[
k^2 \bone_{F\cap\{ S^1_2 \leq t_2, S^2_1 \leq t_2\}} \bigr] =0
\end{equation}
and
%
%
\begin{equation}
\label{f196} \lim_{k\to\infty} \mathbb{E}_{Q_k} \bigl[
k^2 \bone_{F\cap\{ S^1_1 \leq t_2, S^2_2 \leq t_2\}} \bigr] =0.
\end{equation}

Our usual arguments give for $m\geq0$,
\[
\mathbb{P}_{Q_k}\bigl( S^1_m \leq
t_2\mid F\bigr) \leq\bigl(1- e^{- (t_2 - t_1)/(k-1)} \bigr)^m
\leq2(t_2 - t_1)^m (k-1)^{-m},
\]
so
\begin{eqnarray*}
\mathbb{E}_{Q_k} \biggl[ k^2 \sum
_{m,r\geq3} \bone_{F\cap\{
S^1_m \leq t_2, S^2_r \leq t_2\}} \biggr] &\leq&
\mathbb{E}_{Q_k} \biggl[ m k^2 \sum
_{ m\geq3 } \bone_{F\cap\{
S^1_m \leq t_2\}} \biggr]
\\
&\leq& m k^2 \sum_{ m\geq3 }
2(t_2 - t_1)^m (k-1)^{-m},
\end{eqnarray*}
and, therefore,
\[
\lim_{k\to\infty} \mathbb{E}_{Q_k} \biggl[ k^2
\sum_{m,r\geq3} \bone_{F\cap\{
S^1_m \leq t_2, S^2_r \leq t_2\}} \biggr] =0.
\]
This and (\ref{f194})--(\ref{f196}) give
\[
\lim_{k\to\infty} \mathbb{E}_{Q_k} \biggl[ k^2
\sum_{m,r\geq1} \bone_{F\cap\{
S^1_m \leq t_2, S^2_r \leq t_2\}} \biggr]
\leq(t_2 - t_1)^2.
\]
We deduce (\ref{f182}) from the last estimate, (\ref{f183}) and
(\ref{f187})--(\ref{f1811}). This completes the proof.
\end{pf*}

\begin{pf*}{Proof of Corollary~\ref{f172}}
(i)
The process $\{N^1_t, t\geq0\}$ is Poisson with rate one.
It is routine to check that the exponential distribution with mean 1 is
the stationary distribution for the process $ t\to t -T^1_t$. This
easily implies
part (i) of the corollary.

(ii)
Processes $N^j$, $j = 1,\ldots, n$, are independent, so processes $\{
t - T^j_t, t\geq0\}$, $j = 1,\ldots, n$, are independent. This and
Theorem~\ref{f171} imply part (ii).
\end{pf*}

\section{Earthworm effect}\label{earth}

An ``earthworm'' model was introduced in \cite{BCP}. The model
involves a ball moving in a Euclidean torus which pushes ``soil
particles'' aside. The motion of the center of the ball is that of
Brownian motion.
The paper~\cite{BCP} contains a result which suggests that in
dimensions 3 and higher, the ``spherical earthworm''
does not compactify the soil on a global scale, assuming that the torus
diameter is much larger than that of the ball (the result is
asymptotic, in other words). The result in \cite{BCP} does not answer
a number of conjectures stated in that paper. Finishing that research
program appears to involve major technical challenges. In this article,
we will present a discrete version of the earthworm model and a result
that is closer\vspace*{1pt} to the conjectures stated in \cite{BCP}, at least at
the heuristic level.
We will show that if $G = \calC_n^d$ is a torus with a large diameter,
then in the long run, the soil will be uniformly distributed over $G$,
in an appropriate sense, as a result of earthworm's stirring action.

We now present the rigorous version of the ``earthworm'' model.
Given a graph~$G$ with\vspace*{1pt} a vertex set $V$, we will define the mass
process $\mathcal{M}_t = (M^{v_1}_t, M^{v_2}_t,\ldots)$, with an
evolution different than that in the previous sections of the paper.
Suppose that $B_t$ (the ``earthworm'') is a simple random walk on $G$,
that is, $B_t$ is a Markov process which takes values in $V$, stays
constant for an exponential (mean 1) amount of time, and jumps to a
uniformly chosen nearest neighbor at the end of the exponential holding
time. At the time $t$ of a jump of $B$, $M^{B_{t}}$ jumps to~0. At the
same time, the mass $M^{B_{t}}_{t-}$ is ``distributed'' to all adjacent
sites, that is, for every vertex $x$ connected to $B_{t}$ by an edge,
the process $M^x$ increases by $M^{B_{t}}_{t-}/d_v$, that is,
$M^x_t = M^x_{t-} + M^{B_{t}}_{t-}/d_v$. The processes $M^v$ are
constant between the jumps of $B$. The mass $M^v$ can jump only when
$B$ jumps to $v$ or a neighbor of~$v$ in the graph $G$.

Let $\bM$ be the empirical measure for the process
$\{M^v_t\}_{v\in V}$, that is,
$\bM_t = \sum_{v\in V} \delta_{M^v_t}$, where\vspace*{1pt} $\delta_x$ stands for
the measure with a unit atom at $x$ (``Dirac's delta'').

Note that in the following theorem, by the symmetry of the torus, the
initial position of $B$ is irrelevant, so we may assume that $B_0
=\bzero:=(0,\ldots,0)$.

%
\begin{theorem}
Fix $d\geq1$, and let $\bM^n_t$ be the empirical measure process for
the earthworm process on the graph $G=\calC_n^d$. Assume that $M^v_0=
1/n^d$ for $v\in V$ (hence $\sum_{v\in V} M^v_0 =1$).
\begin{longlist}[(ii)]
\item[(i)] For every $n$, the random measures $\bM^n_t$ converge weakly to a
random measure $\bM^n_\infty$, when $t\to\infty$.

\item[(ii)] For\vspace*{2pt} $R \subset\mathbb{R}^d$ and $a\in\mathbb{R}$, let $aR =
\{x\in\mathbb{R}^d\dvtx  x= ay \mbox{ for some } y \in R\}$ and
$\widehat\bM^n_\infty(R) = \bM^n_\infty(nR) $. When $n\to\infty
$, the random measures $\widehat\bM^n_\infty$ converge weakly to the
random measure equal to, a.s., the uniform probability measure on $[0,1]^d$.
\end{longlist}
\end{theorem}

\begin{pf}
(i) The proof of Theorem~\ref{a297}
applies in the present case, with some minor modification accounting
for the fact that the mass redistribution mechanism is given by $B$
rather than Poisson processes $N^x$. Hence, there exists a unique
stationary distribution $Q$ for $(\mathcal{M}_t, B_t)$. Under $Q$,
$\bM^n_t$ has distribution $\bM^n_\infty$.\vspace*{2pt}

(ii)
Let $|R|$ denote the $d$-dimensional Lebesgue measure of $R\subset
\mathbb{R}^d$.
To prove part (ii) of the theorem, it will suffice to show that for
every fixed rectangle $R\subset[0,1]^d$ with rational vertices,
$\lim_{n\to\infty}\widehat\bM^n_\infty(R) = |R|$, in probability.
It will be enough to show that
%
%
\begin{equation}
\label{f1910} \lim_{n\to\infty} \mathbb{E}_Q \widehat
\bM^n_\infty(R) = |R|
\end{equation}
and
%
%
\begin{equation}
\label{f1911} \lim_{n\to\infty} \operatorname{Var}_Q
\widehat\bM^n_\infty(R) = 0.
\end{equation}

By symmetry,
$\mathbb{E}_{Q} M^x_0 = \mathbb{E}_{Q} M^y_0$
for all $x,y \in V$. Since $\sum_{x\in V} M^x_0 = 1$ under $Q$,
we must have $\mathbb{E}_{Q} M^x_0=1/n^d$.
By abuse of notation, we give $| \cdot|$ another meaning---it will
denote the cardinality of an (at most) countable set. We have
\[
\lim_{n\to\infty} \mathbb{E}_Q \widehat
\bM^n_\infty(R) = \lim_{n\to\infty}
\mathbb{E}_Q \bM^n_\infty(nR) = \lim
_{n\to\infty}\frac{1} {n^d} |nR \cap V| = |R|,
\]
and thus (\ref{f1910}) is proved.

Let $Z$ and $\widetilde Z$ be defined as $Z^1$ and $ Z^2$ in Definition
\ref{wimps}.
In particular,
$\mathbb{P}(Z_0 = x) = \mathbb{P}(\widetilde Z_0 = x) = M^x_0$ for
$x\in V$.
However, note that in the present case, the process $\{Z_t, t\geq0\}$
jumps at a time $t$ if and only if
$B$ jumps at the time $t$ and $B_t=Z_{t-}$.
A similar remark applies to $\{\widetilde Z_t, t\geq0\}$.
Note that the jump times of $Z$ and $\widetilde Z$ are defined by the
same process $B$.

The state space for the process $(Z,\widetilde Z)$ is finite, so it has
a stationary distribution. The stationary distribution is unique
because all states communicate.
We will next estimate the stationary probabilities, in the asymptotic
sense, when $n\to\infty$.

Let $\overline Z_t = Z_t - \widetilde Z_t\in V$ (in the sense of group
operations on the Cayley graph). Although $\overline Z_t$ is not a
Markov process (as far as we can tell), it is clear how to define a
discrete time Markov chain $\{U_j, j\geq1\}$ embedded in $\overline Z_t$.

For $x\in V$, let $\calB(x,r)$ denote the closed ball in $V$ with
center $x$ and radius $r$, relative to the graph distance.

Let $\bfh$ be set of all vertices $(f_1,\ldots, f_d) $ such
that $|f_i| =1$ for some $i$, and $f_m = 0$ for all $m\ne i$.
It has been shown in the proof of Theorem~\ref{m61} that the
stationary distribution $\{\pi_x, x\in V\}$ for $U$ has the\vspace*{1pt}
following form. For some normalizing constant $c_1>0$, $\pi_x =
(1 - \frac{1}{2d}) c_1$ for all $x\in\bfh$ and $\pi_x = c_1$
for all other $x\in V$.

Although $\overline Z$ is not a Markov process, $(\overline Z, Z,
\widetilde Z)$ is.
We will consider the process $(\overline Z, Z, \widetilde Z)$ in the
stationary regime.
Let $\overline Z$ have the corresponding marginal distribution.
We will estimate the proportion of time that $\overline Z$ spends in
different states. For each state $x$, we will estimate the product of
$\pi_x$ and the expected amount of time between the time $\tau_1$ of
the first jump of $\overline Z$ to $x$ and the time $\tau_2$ of the
next jump. Let us call the random time between these jumps $\tau_x =
\tau_2-\tau_1$. Hence we will estimate $\mathbb{E}\tau_x$.

Consider $x \in\calB(\bzero, 2)^c$ and any two neighbors $y$ and $z$
of $x$. We have $\pi_y=\pi_z $, so the probability that the process
$\overline Z$ jumps to $x$ from $y$ is equal to the probability that
the process $\overline Z$ jumps to $x$ from $z$. Hence,
$\mathbb{P}(B_{\tau_1} =x)=\mathbb{P}(B_{\tau_1} =y)$
for any neighbors $x$ and $y$ of $Z_{\tau_1}$ and $\widetilde Z_{\tau_1}$.
The time $\tau_x$ is the same as the waiting time for the first
hit of
$\{Z_{\tau_1},\widetilde Z_{\tau_1}\}$ after time $\tau_1$, for $B$.

Let $K$ be the set of all neighbors of $Z_{\tau_1}$ and $\widetilde
Z_{\tau_1}$. We have shown that
the distribution of $B_{\tau_1}$ is uniform on $K$. It follows
from \cite{AF}, Corollary~24, page~21, Chapter~2,
that the
expected time until $B$ hits $\{Z_{\tau_1},\widetilde Z_{\tau_1}\}$
is $|V| /2 -1$. This implies that $
\mathbb{E}\tau_x = |V|/2 -1 $.
Thus for any
$x,y \in\calB(\bzero, n_1)^c$, we have
$\mathbb{E}\tau_x =\mathbb{E}\tau_y $.
This and the fact that $\pi_x=\pi_y $
imply\vspace*{2pt} that, under the stationary distribution, for
$x,y \in\calB(\bzero, 2)^c$,
$\mathbb{P}(\overline Z_0 = x) = \mathbb{P}(\overline Z_0 = y)$.
Therefore, if $x-y \in\calB(\bzero, 2)^c$,
%
%
\begin{equation}
\label{f235} \mathbb{P}_{Q}(Z_0 =x,\widetilde
Z_0 =y ) = \mathbb{P}_{Q}(Z_0 =x)
\mathbb{P}_{Q}(\widetilde Z_0 =y ).
\end{equation}
For $x \in\calB(\bzero, 2)$ we have a rough bound
$\mathbb{E}\tau_x \leq c_1 n^d$, which yields for $x-y \in\calB
(\bzero, 2)$,
%
%
\begin{equation}
\label{f243} \mathbb{P}_{Q}(Z_0 =x,\widetilde
Z_0 =y ) \leq c_2\mathbb{P}_{Q}(Z_0
=x)\mathbb{P}_{Q}(\widetilde Z_0 =y ).
\end{equation}

Let $\mathcal{G}_t = \sigma(\mathcal{M}_s, 0\leq s \leq t)
= \sigma(\mathcal{M}_0,B_s, 0\leq s \leq t)$.
We have for $x\in V$,
\[
\mathbb{P}_{Q}(Z_0 =x \mid\mathcal{G}_0) =
M^x_0.
\]
The processes $Z_t$ and $\widetilde Z_t$
are conditionally independent given $\mathcal{G}_t$, so for $x,y\in V$,
\[
\mathbb{P}_{Q}(Z_t =x,\widetilde Z_t =y
\mid\mathcal{G}_t) = M^x_t
M^y_t.
\]
By stationarity, for $x,y\in V$,
\[
\mathbb{P}_{Q}(Z_0 =x,\widetilde Z_0 =y
\mid\mathcal{G}_0) = M^x_0
M^y_0.
\]
Thus
\begin{eqnarray*}
\mathbb{E}_{Q} \bigl(M^x_0
M^y_0\bigr) &=& \mathbb{E}_{Q}
\mathbb{P}_{Q}(Z_0 =x,\widetilde Z_0 =y \mid
\mathcal{G}_0) =\mathbb{P}_{Q}(Z_0 =x,
\widetilde Z_0 =y ).
\end{eqnarray*}
We obtain
%
%
\begin{eqnarray}\label{f241}
\qquad\mathbb{E}_Q \bigl(\bM^n_\infty(nR)
\bigr)^2 &=& \sum_{x,y \in nR}
\mathbb{E}_{Q} \bigl(M^x_0
M^y_0\bigr)
\nonumber\\[-8pt]\\[-8pt]
&=& \mathop{\sum_{x,y \in nR}}_{x-y \in\calB(\bzero, 2)^c} \mathbb
{E}_{Q} \bigl(M^x_0 M^y_0
\bigr) + \mathop{\sum_{x,y \in nR}}_{x-y \in\calB(\bzero, 2)} \mathbb
{E}_{Q} \bigl(M^x_0 M^y_0
\bigr).\nonumber
\end{eqnarray}
It follows from (\ref{f243}) that
%
%
\begin{eqnarray}\label{f242}
\mathop{\sum_{x,y \in nR}}_{x-y \in\calB(\bzero, 2)}
\mathbb{E}_{Q} \bigl(M^x_0
M^y_0\bigr) &=& \mathop{\sum_{x,y \in nR}}_{x-y \in\calB(\bzero, 2)}
\mathbb{P}_{Q}(Z_0 =x,\widetilde Z_0 =y )\nonumber
\\
&\leq&\mathop{\sum_{x,y \in nR}}_{x-y \in\calB(\bzero, 2)}
c_2\mathbb{P}_{Q}(Z_0 =x)
\mathbb{P}_{Q}(\widetilde Z_0 =y )
\\
&\leq&  c_3 |R|/n^d.\nonumber
\end{eqnarray}
We use (\ref{f235}) to see that
\begin{eqnarray*}
\mathop{\sum_{x,y \in nR}}_{x-y \in\calB(\bzero, 2)^c} \mathbb
{E}_{Q} \bigl(M^x_0 M^y_0
\bigr) &=& \mathop{\sum_{x,y \in nR}}_{x-y \in\calB(\bzero, 2)^c} \mathbb
{P}_{Q}(Z_0 =x,\widetilde Z_0 =y )
\\
&\leq&\mathop{\sum_{x,y \in nR}}_{x-y \in\calB(\bzero, 2)^c}
\mathbb{P}_{Q}(Z_0 =x)\mathbb{P}_{Q}(
\widetilde Z_0 =y ) \leq|R|^2.
\end{eqnarray*}
This, (\ref{f241}) and (\ref{f242}) give
\[
\lim_{n\to\infty} \mathbb{E}_Q \bigl(
\bM^n_\infty(nR)\bigr)^2 \leq\lim
_{n\to\infty} \bigl( c_3 |R|/n^d +
|R|^2\bigr) = |R|^2.
\]
We obtain,
$\lim_{n\to\infty} \mathbb{E}_Q (\bM^n_\infty(nR))^2 \leq
|R|^2$. This shows (\ref{f1911}), thus completing the proof.
\end{pf}

\section*{Acknowledgments}
We are grateful to Mikl\'os B\'ona, Harry Crane, Persi Diaconis,
Jason Fulman, Ron Irving, Svante Janson, David Levin, Yuval Peres and
Jon Wellner for the most useful advice.
We thank the referee for very careful reading of the paper and many
suggestions for improvement.





\printaddresses

\begin{thebibliography}{29}
\bibitem{AF}
\begin{bmisc}[auto:STB|2014/05/28|10:36:42]
\bauthor{\bsnm{Aldous},~\bfnm{David}\binits{D.}} \AND
\bauthor{\bsnm{Fill},~\bfnm{James}\binits{J.}}
(\byear{2014}).
\bhowpublished{\textit{Reversible Markov Chains and Random Walks on Graphs}.
Book in preparation. Available at \url{http://www.stat.berkeley.edu/\textasciitilde aldous/RWG/book.html}.}
\end{bmisc}
\bptok{imsref}%
\endbibitem

\bibitem{Ald}
\begin{barticle}[mr]
\bauthor{\bsnm{Aldous},~\bfnm{David~J.}\binits{D.~J.}}
(\byear{1991}).
\btitle{Meeting times for independent {M}arkov chains}.
\bjournal{Stochastic Process. Appl.}
\bvolume{38}
\bpages{185--193}.
\bid{doi={10.1016/0304-4149(91)90090-Y}, issn={0304-4149}, mr={1119980}}
\end{barticle}
\bptok{imsref}%
\endbibitem

\bibitem{BBS}
\begin{barticle}[mr]
\bauthor{\bsnm{Billey},~\bfnm{Sara}\binits{S.}},
\bauthor{\bsnm{Burdzy},~\bfnm{Krzysztof}\binits{K.}} \AND
\bauthor{\bsnm{Sagan},~\bfnm{Bruce~E.}\binits{B.~E.}}
(\byear{2013}).
\btitle{Permutations with given peak set}.
\bjournal{J.~Integer Seq.}
\bvolume{16}
\bpages{Article 13.6.1, 18}.
\bid{issn={1530-7638}, mr={3083179}}
\end{barticle}
\bptok{imsref}%
\endbibitem

\bibitem{Bill}
\begin{bbook}[mr]
\bauthor{\bsnm{Billingsley},~\bfnm{Patrick}\binits{P.}}
(\byear{1968}).
\btitle{Convergence of Probability Measures}.
\bpublisher{Wiley},
\blocation{New York}.
\bid{mr={0233396}}
\end{bbook}
\bptok{imsref}%
\endbibitem

\bibitem{Bona}
\begin{bmisc}[auto:STB|2014/05/28|10:36:42]
\bauthor{\bsnm{B{\'o}na},~\bfnm{Mikl{\'o}s}\binits{M.}}
(\byear{2007}).
\bhowpublished{The copies of any permutation pattern are asymptotically normal.
Available at \arxivurl{arXiv:0712.2792}.}
\end{bmisc}
\bptok{imsref}%
\endbibitem

\bibitem{meteor2}
\begin{bmisc}[auto:STB|2014/05/28|10:36:42]
\bauthor{\bsnm{Burdzy},~\bfnm{Krzysztof}\binits{K.}}
(\byear{2013}).
\bhowpublished{Meteor process on $\mathbb{Z}^d$.
Available at \arxivurl{arXiv:1312.6865}.}
\end{bmisc}
\bptok{imsref}%
\endbibitem

\bibitem{BCP}
\begin{barticle}[mr]
\bauthor{\bsnm{Burdzy},~\bfnm{Krzysztof}\binits{K.}},
\bauthor{\bsnm{Chen},~\bfnm{Zhen-Qing}\binits{Z.-Q.}} \AND
\bauthor{\bsnm{Pal},~\bfnm{Soumik}\binits{S.}}
(\byear{2013}).
\btitle{Brownian earthworm}.
\bjournal{Ann. Probab.}
\bvolume{41}
\bpages{4002--4049}.
\bid{doi={10.1214/12-AOP831}, issn={0091-1798}, mr={3161468}}
\end{barticle}
\bptok{imsref}%
\endbibitem

\bibitem{CLR}
\begin{barticle}[mr]
\bauthor{\bsnm{Caputo},~\bfnm{Pietro}\binits{P.}},
\bauthor{\bsnm{Liggett},~\bfnm{Thomas~M.}\binits{T.~M.}} \AND
\bauthor{\bsnm{Richthammer},~\bfnm{Thomas}\binits{T.}}
(\byear{2010}).
\btitle{Proof of {A}ldous' spectral gap conjecture}.
\bjournal{J. Amer. Math. Soc.}
\bvolume{23}
\bpages{831--851}.
\bid{doi={10.1090/S0894-0347-10-00659-4}, issn={0894-0347}, mr={2629990}}
\end{barticle}
\bptok{imsref}%
\endbibitem

\bibitem{CP}
\begin{barticle}[mr]
\bauthor{\bsnm{Chan},~\bfnm{O.-Yeat}\binits{O.-Y.}} \AND
\bauthor{\bsnm{Pra{\l}at},~\bfnm{Pawe{\l}}\binits{P.}}
(\byear{2012}).
\btitle{Chipping away at the edges: How long does it take?}
\bjournal{J.~Comb.}
\bvolume{3}
\bpages{101--121}.
\bid{doi={10.4310/JOC.2012.v3.n1.a5}, issn={2156-3527}, mr={2975324}}
\end{barticle}
\bptok{imsref}%
\endbibitem

\bibitem{Chao}
\begin{barticle}[mr]
\bauthor{\bsnm{Chao},~\bfnm{Chern-Ching}\binits{C.-C.}}
(\byear{1997}).
\btitle{A note on applications of the martingale central limit theorem to random permutations}.
\bjournal{Random Structures Algorithms}
\bvolume{10}
\bpages{323--332}.
\bid{doi={10.1002/(SICI)1098-2418(199705)10:3<323::AID-RSA2>3.3.CO;2-R}, issn={1042-9832}, mr={1606222}}
\end{barticle}
\bptok{imsref}%
\endbibitem

\bibitem{CGS}
\begin{bbook}[mr]
\bauthor{\bsnm{Chen},~\bfnm{Louis~H.~Y.}\binits{L.~H.~Y.}},
\bauthor{\bsnm{Goldstein},~\bfnm{Larry}\binits{L.}} \AND
\bauthor{\bsnm{Shao},~\bfnm{Qi-Man}\binits{Q.-M.}}
(\byear{2011}).
\btitle{Normal Approximation by {S}tein's Method}.
\bpublisher{Springer},
\blocation{Heidelberg}.
\bid{doi={10.1007/978-3-642-15007-4}, mr={2732624}}
\end{bbook}
\bptok{imsref}%
\endbibitem

\bibitem{CongVis}
\begin{barticle}[mr]
\bauthor{\bsnm{Conger},~\bfnm{Mark}\binits{M.}} \AND
\bauthor{\bsnm{Viswanath},~\bfnm{D.}\binits{D.}}
(\byear{2007}).
\btitle{Normal approximations for descents and inversions of permutations of multisets}.
\bjournal{J.~Theoret. Probab.}
\bvolume{20}
\bpages{309--325}.
\bid{doi={10.1007/s10959-007-0070-5}, issn={0894-9840}, mr={2324533}}
\end{barticle}
\bptok{imsref}%
\endbibitem

\bibitem{C14}
\begin{barticle}[auto:STB|2014/05/28|10:36:42]
\bauthor{\bsnm{Crane},~\bfnm{Harry}\binits{H.}}
(\byear{2014}).
\btitle{The cut-and-paste process}.
\bjournal{Ann. Probab.}
\bvolume{42}
\bpages{1952--1979}.
\bid{mr={3262496}}
\end{barticle}
\bptok{imsref}%
\endbibitem

\bibitem{CL}
\begin{barticle}[mr]
\bauthor{\bsnm{Crane},~\bfnm{Harry}\binits{H.}} \AND
\bauthor{\bsnm{Lalley},~\bfnm{Steven~P.}\binits{S.~P.}}
(\byear{2013}).
\btitle{Convergence rates of {M}arkov chains on spaces of partitions}.
\bjournal{Electron. J. Probab.}
\bvolume{18}
\bpages{1--23}.
\bid{doi={10.1214/EJP.v18-2389}, issn={1083-6489}, mr={3078020}}
\end{barticle}
\bptok{imsref}%
\endbibitem

\bibitem{dHF}
\begin{bbook}[mr]
\bauthor{\bparticle{de} \bsnm{Haan},~\bfnm{Laurens}\binits{L.}} \AND
\bauthor{\bsnm{Ferreira},~\bfnm{Ana}\binits{A.}}
(\byear{2006}).
\btitle{Extreme Value Theory: An Introduction}.
\bpublisher{Springer},
\blocation{New York}.
\bid{mr={2234156}}
\end{bbook}
\bptok{imsref}%
\endbibitem

\bibitem{DF99}
\begin{barticle}[mr]
\bauthor{\bsnm{Diaconis},~\bfnm{Persi}\binits{P.}} \AND
\bauthor{\bsnm{Freedman},~\bfnm{David}\binits{D.}}
(\byear{1999}).
\btitle{Iterated random functions}.
\bjournal{SIAM Rev.}
\bvolume{41}
\bpages{45--76}.
\bid{doi={10.1137/S0036144598338446}, issn={0036-1445}, mr={1669737}}
\end{barticle}
\bptok{imsref}%
\endbibitem

\bibitem{Dyson}
\begin{barticle}[mr]
\bauthor{\bsnm{Dyson},~\bfnm{Freeman~J.}\binits{F.~J.}}
(\byear{1962}).
\btitle{A {B}rownian-motion model for the eigenvalues of a random matrix}.
\bjournal{J.~Math. Phys.}
\bvolume{3}
\bpages{1191--1198}.
\bid{issn={0022-2488}, mr={0148397}}
\end{barticle}
\bptok{imsref}%
\endbibitem

\bibitem{EM}
\begin{barticle}[mr]
\bauthor{\bsnm{Ehrenborg},~\bfnm{Richard}\binits{R.}} \AND
\bauthor{\bsnm{Mahajan},~\bfnm{Swapneel}\binits{S.}}
(\byear{1998}).
\btitle{Maximizing the descent statistic}.
\bjournal{Ann. Comb.}
\bvolume{2}
\bpages{111--129}.
\bid{doi={10.1007/BF01608482}, issn={0218-0006}, mr={1682923}}
\end{barticle}
\bptok{imsref}%
\endbibitem

\bibitem{FF}
\begin{barticle}[mr]
\bauthor{\bsnm{Ferrari},~\bfnm{P.~A.}\binits{P.~A.}} \AND
\bauthor{\bsnm{Fontes},~\bfnm{L.~R.~G.}\binits{L.~R.~G.}}
(\byear{1998}).
\btitle{Fluctuations of a surface submitted to a random average process}.
\bjournal{Electron. J. Probab.}
\bvolume{3}
\bpages{34 pp. (electronic)}.
\bid{doi={10.1214/EJP.v3-28}, issn={1083-6489}, mr={1624854}}
\end{barticle}
\bptok{imsref}%
\endbibitem

\bibitem{FMQR}
\begin{barticle}[mr]
\bauthor{\bsnm{Fey-den Boer},~\bfnm{Anne}\binits{A.}},
\bauthor{\bsnm{Meester},~\bfnm{Ronald}\binits{R.}},
\bauthor{\bsnm{Quant},~\bfnm{Corrie}\binits{C.}} \AND
\bauthor{\bsnm{Redig},~\bfnm{Frank}\binits{F.}}
(\byear{2008}).
\btitle{A probabilistic approach to {Z}hang's sandpile model}.
\bjournal{Comm. Math. Phys.}
\bvolume{280}
\bpages{351--388}.
\bid{doi={10.1007/s00220-008-0470-0}, issn={0010-3616}, mr={2395474}}
\end{barticle}
\bptok{imsref}%
\endbibitem

\bibitem{FK60}
\begin{barticle}[mr]
\bauthor{\bsnm{Furstenberg},~\bfnm{H.}\binits{H.}} \AND
\bauthor{\bsnm{Kesten},~\bfnm{H.}\binits{H.}}
(\byear{1960}).
\btitle{Products of random matrices}.
\bjournal{Ann. Math. Statist.}
\bvolume{31}
\bpages{457--469}.
\bid{issn={0003-4851}, mr={0121828}}
\end{barticle}
\bptok{imsref}%
\endbibitem

\bibitem{HMS}
\begin{barticle}[mr]
\bauthor{\bsnm{Hairer},~\bfnm{M.}\binits{M.}},
\bauthor{\bsnm{Mattingly},~\bfnm{J.~C.}\binits{J.~C.}} \AND
\bauthor{\bsnm{Scheutzow},~\bfnm{M.}\binits{M.}}
(\byear{2011}).
\btitle{Asymptotic coupling and a general form of {H}arris' theorem with applications to stochastic delay equations}.
\bjournal{Probab. Theory Related Fields}
\bvolume{149}
\bpages{223--259}.
\bid{doi={10.1007/s00440-009-0250-6}, issn={0178-8051}, mr={2773030}}
\end{barticle}
\bptok{imsref}%
\endbibitem

\bibitem{HKPV}
\begin{bbook}[mr]
\bauthor{\bsnm{Hough},~\bfnm{J.~Ben}\binits{J.~B.}},
\bauthor{\bsnm{Krishnapur},~\bfnm{Manjunath}\binits{M.}},
\bauthor{\bsnm{Peres},~\bfnm{Yuval}\binits{Y.}} \AND
\bauthor{\bsnm{Vir{\'a}g},~\bfnm{B{\'a}lint}\binits{B.}}
(\byear{2009}).
\btitle{Zeros of {G}aussian Analytic Functions and Determinantal Point Processes}.
\bseries{University Lecture Series}
\bvolume{51}.
\bpublisher{Amer. Math. Soc.},
\blocation{Providence, RI}.
\bid{mr={2552864}}
\end{bbook}
\bptok{imsref}%
\endbibitem

\bibitem{HW}
\begin{barticle}[mr]
\bauthor{\bsnm{Howitt},~\bfnm{Chris}\binits{C.}} \AND
\bauthor{\bsnm{Warren},~\bfnm{Jon}\binits{J.}}
(\byear{2009}).
\btitle{Consistent families of {B}rownian motions and stochastic flows of kernels}.
\bjournal{Ann. Probab.}
\bvolume{37}
\bpages{1237--1272}.
\bid{doi={10.1214/08-AOP431}, issn={0091-1798}, mr={2546745}}
\end{barticle}
\bptok{imsref}%
\endbibitem

\bibitem{AK}
\begin{bmisc}[auto:STB|2014/05/28|10:36:42]
\bauthor{\bsnm{Kasraoui},~\bfnm{Anisse}\binits{A.}}
(\byear{2012}).
\bhowpublished{The most frequent peak set of a random permutation.
Available at \arxivurl{arXiv:1210.5869}.}
\end{bmisc}
\bptok{imsref}%
\endbibitem

\bibitem{KMK2}
\begin{barticle}[auto:STB|2014/05/28|10:36:42]
\bauthor{\bsnm{Kermack},~\bfnm{W.~O.}\binits{W.~O.}} \AND
\bauthor{\bsnm{McKendrick},~\bfnm{A.~G.}\binits{A.~G.}}
(\byear{1937}).
\btitle{Some distributions associated with a randomly arranged set of numbers}.
\bjournal{Proc. Roy. Soc. Edinburgh}
\bvolume{57}
\bpages{332--376}.
\end{barticle}
\bptok{imsref}%
\endbibitem

\bibitem{KMK1}
\begin{barticle}[auto:STB|2014/05/28|10:36:42]
\bauthor{\bsnm{Kermack},~\bfnm{W.~O.}\binits{W.~O.}} \AND
\bauthor{\bsnm{McKendrick},~\bfnm{A.~G.}\binits{A.~G.}}
(\byear{1937}).
\btitle{Tests for randomness in a series of numerical observations}.
\bjournal{Proc. Roy. Soc. Edinburgh}
\bvolume{57}
\bpages{228--240}.
\end{barticle}
\bptok{imsref}%
\endbibitem

\bibitem{LPW}
\begin{bbook}[mr]
\bauthor{\bsnm{Levin},~\bfnm{David~A.}\binits{D.~A.}},
\bauthor{\bsnm{Peres},~\bfnm{Yuval}\binits{Y.}} \AND
\bauthor{\bsnm{Wilmer},~\bfnm{Elizabeth~L.}\binits{E.~L.}}
(\byear{2009}).
\btitle{Markov Chains and Mixing Times}.
\bpublisher{Amer. Math. Soc.},
\blocation{Providence, RI}.
\bid{mr={2466937}}
\end{bbook}
\bptok{imsref}%
\endbibitem

\bibitem{math}
\begin{bbook}[auto]
\bauthor{Wolfram Research}
(\byear{2010}).
\btitle{Mathematica. {V}ersion 8.0}.
\bpublisher{Wolfram Research}, \blocation{Champaign, IL}.
\end{bbook}
\bptok{imsref}%
\endbibitem

\end{thebibliography}
\end{document}